\documentclass[11pt]{elsarticle}
\usepackage{graphicx}
\usepackage{caption}
\usepackage{subcaption}
\usepackage{amsmath}
\usepackage{amsfonts}
\usepackage{amssymb}
\usepackage{amsthm}
\usepackage[margin=1in]{geometry}
\usepackage{enumitem}
\usepackage{float}
\usepackage{algorithm}
\usepackage{algpseudocode}
\usepackage{multirow}
\usepackage{colortbl}
\usepackage{hhline}
\usepackage{makecell}
\usepackage{arydshln}
\usepackage[table]{xcolor}
\usepackage{thmtools}
\usepackage{thm-restate}
\usepackage{circledsteps}
\numberwithin{equation}{section}
\usepackage{hyperref}

\biboptions{sort&compress}

\declaretheorem[name=Theorem]{theorem}
\declaretheorem[sibling=theorem, name=Proposition]{proposition}
\declaretheorem[sibling=theorem, name=Lemma]{lemma}
\declaretheorem[sibling=theorem, style=remark, name=Remark]{remark}

\newcommand{\manifold}{\Omega}

\newcommand{\Diff}{\mathcal{D}(\Omega)}
\newcommand{\SDiff}{\mathcal{D}_\mu(\Omega)}
\newcommand{\Prob}{\mathcal{P}(\Omega)}
\newcommand{\Vol}{\mathcal{V}(\Omega)}
\newcommand{\Hor}{\operatorname{Hor}}
\newcommand{\Ver}{\operatorname{Ver}}

\newcommand{\forwardMap}{\varphi}
\newcommand{\map}{\mathcal{X}}
\newcommand{\numericalMap}{X}
\newcommand{\jacobian}{j}
\newcommand{\numericalJacobian}{J}
\newcommand{\mass}{m}
\newcommand{\numericalMass}{M}
\newcommand{\jacobianOperator}{\operatorname{Jac}}
\newcommand{\liftingOperator}{\mathcal{K}}
\newcommand{\solution}{\phi}
\newcommand{\jsource}{q}
\newcommand{\semidirectGroup}{\mathcal{G}}
\newcommand{\groupElement}{\boldsymbol{\mathcal{X}}}
\newcommand{\numericalGroupElement}{\boldsymbol{X}}
\newcommand{\integrator}{\mathcal{S}}
\newcommand{\intMap}{\integrator_\map}
\newcommand{\intJac}{\integrator_\jacobian}
\newcommand{\intMass}{\integrator_\mass}
\newcommand{\slPropagator}{\mathcal{T}}

\newcommand{\interp}{\mathcal{I}}
\newcommand{\grpInterp}{\boldsymbol{\mathcal{I}}}

\newcommand{\productGrp}{\mathcal{G}}

\newcommand{\isoSubGp}{\mathcal{H}}

\newcommand{\fwdGrpElt}{\boldsymbol{\forwardMap}}

\newcommand{\metric}{\mathbb{G}}

\algnewcommand{\LineComment}[1]{\State \(\triangleright\) #1}

\graphicspath{{figures/}}

\makeatletter
\newcommand{\doublehhline}[1]{%
  \hhline{#1}%
  \noalign{%
    \adl@hline\relax{\dimexpr2\arrayrulewidth+\doublerulesep\relax}%
  }%
}
\makeatother
\title{A characteristic mapping method for dissipative evolution equations}
\date{}
\author[1]{William Holman--Bissegger}
\author[2]{Seth Taylor}
\author[1]{Jean--Christophe Nave}

\fntext[1]{Department of Mathematics and Statistics, McGill University, Montr\'{e}al, Qu\'{e}bec, Canada.}
\fntext[2]{Department of Computer Science, University of Saskatchewan, Saskatoon, Saskatchewan, Canada.}

\begin{document}

\begin{abstract}
The Characteristic Mapping (CM) method is a transport scheme which discretizes the Lagrangian flow map as a composition of spatial interpolants. This adaptive approximation technique has demonstrated novel resolution properties and the ability to efficiently simulate the fine scales of inviscid turbulent flows. In this work, we construct an extension of this approach to the solution of dissipative problems, including the Navier--Stokes equations. The extension lifts a general scalar evolution equation to a path in a group of diffeomorphisms, where the semigroup approach of the CM method can be exploited. We give an interpretation of this lift through the geometric framework of optimal transport. We then introduce a semidirect product evolution technique to perform the time stepping of the map components. The map is evolved based on the semi-Lagrangian CM approach, and its log-Jacobian is evolved explicitly by an exponential integrator. The numerical scheme is demonstrated to be high-order accurate and is tested on a range of canonical nonlinear dissipative equations. By retaining the resolution properties of the CM method, the approach shows a promising ability to simulate turbulent fluid flows across all Reynolds numbers.
\end{abstract}

\maketitle

\section{Introduction}

The direct numerical simulation of turbulent systems is fundamentally limited by the resolution required to capture dynamical interactions over spatiotemporal scales spanning several orders of magnitude. Traditional fixed-grid methods constrain the finest representable scale, and information is lost as soon as energy cascades down below it. Many existing methods view this numerical challenge as a nonlinear approximation problem, focusing on the design of adaptive approximations that allocate computational resources across space, time, and scale to regions of increased complexity. Transport-dominated systems however already possess a natural nonlinear description associated, in part, to an evolution on the space of diffeomorphisms of the domain.

The Characteristic Mapping (CM) method \cite{mercier_characteristic_2020, yin_characteristic_2021, yin_characteristic_2023, taylor_projection-based_2023, taylor_characteristic_2023, krah_characteristic_2024, yin_characteristic_2026, bergmann_singularity_2024, taylor_diffeomorphic_2025} leverages this structure for the solution of nonlinear transport equations. By discretizing in the space of maps, the method utilizes the group structure to adaptively split a path of diffeomorphisms into a composition of shorter time deformations. The approach has been shown to exhibit exponential spatial resolution using a linear amount of memory, and to recover the extremely-fine scales in turbulent flows associated to the incompressible Euler equations \cite{yin_characteristic_2021, yin_characteristic_2023, taylor_characteristic_2023}, Vlasov-Poisson equations \cite{krah_characteristic_2024}, and ideal magnetohydrodynamics with source terms \cite{yin_characteristic_2026}. It has been further applied to the study of singularity formation \cite{bergmann_singularity_2024} and as a geometric operator learning method \cite{taylor_diffeomorphic_2025}.

The CM method is well-established for equations in explicit transport form. The aim of this work is to extend the method to the solution of a general scalar evolution equation. We consider problems of the form
\begin{equation}
    \label{eq:main_pde}
    \partial_t \solution + u \cdot \nabla \phi =  F(\solution, t), \quad\quad \solution(0) = \solution_0, 
\end{equation}
for a scalar field $\solution(t) \in C^\infty (\manifold)$, $t \in [0, T]$, over a closed Riemannian manifold $\manifold$. Here, $F$ is a (possibly nonlinear) differential operator and $u$ is a divergence-free velocity field which may depend on the solution, $u = A(\solution, t)$.

\subsection{The Characteristic Mapping Method}
\label{section:cmm_intro}

Let $(\manifold, g)$ be a closed Riemannian manifold and let $v(t) \in \mathfrak{X}(\manifold)$ be a time-dependent vector field. The forward flow map $\varphi_{[t_0, t]} : \manifold \rightarrow \manifold$ of $v$ maps particles transported by $v$ from their starting locations at time $t_0$ to their positions at time $t$, and is defined by the ordinary differential equation
\begin{equation*}
    \partial_t \forwardMap_{[t_0, t]} = v(t) \circ \forwardMap_{[t_0, t]}, \quad\quad \forwardMap_{[t_0, t_0]} = \operatorname{id}.
\end{equation*}
The central object of approximation in the CM method is the inverse flow map $\map_{[t, t_0]} := \forwardMap_{[t_0, t]}^{-1}$, which sends particles at time $t$ back to their initial locations at $t_0$. The inverse map satisfies the transport equation
\begin{equation*}
    \partial_t \map_{[t, t_0]} + D \map_{[t, t_0]} (v(t)) = 0, \quad\quad \map_{[t_0, t_0]} = \operatorname{id},
\end{equation*}
describing the advection of the particle labels by the flow. It defines a solution operator to the continuity equation
\begin{equation}
    \label{eq:continuity_equation}
    \partial_t \rho + \operatorname{div}(\rho v) = 0
\end{equation}
by the pullback action on densities,
\begin{equation}
    \label{eq:def_pullback_in_coordinates}
    \rho(t) = (\rho(t_0) \circ \map_{[t, t_0]}) \jacobianOperator(\map_{[t, t_0]}).
\end{equation}
Here, $\jacobianOperator$ denotes the Jacobian determinant with respect to the Riemannian volume form $\mu$, which is defined by $\jacobianOperator(\map) \mu := \map^\ast \mu$.

By approximating the solution operator instead of the transported quantity directly, only the regularity of the velocity, not that of the density, bears a constraint on the evolution. This property has demonstrated the ability to evolve densities of arbitrary complexity on coarse grids \cite{mercier_characteristic_2020}. It has also demonstrated the ability to capture dynamically generated fine scale features \cite{yin_characteristic_2021}. The latter is enabled by the semigroup property of the map: for arbitrary times $\tau_0 < \tau_1 < \cdots < \tau_m$, we have the decomposition
\begin{equation}
    \label{eq:classical_semigroup_decomposition}
    \map_{[\tau_m, \tau_0]} = \map_{[\tau_1, \tau_0]} \circ \map_{[\tau_2, \tau_1]} \circ \cdots \circ \map_{[\tau_m, \tau_{m - 1}]}.
\end{equation}
This induces an analogous decomposition of the solution operator \eqref{eq:def_pullback_in_coordinates}. The set of remapping times $\tau_i$ may be selected either a priori or dynamically. In a CM method, each submap $\map_{[t, \tau_k]}$ is evolved on a grid during the time interval $[\tau_k, \tau_{k + 1}]$ and then stored. The map \eqref{eq:classical_semigroup_decomposition} then lives in a nonlinear approximation space: the representable scale gets finer with each new submap. This approach allows for representing exponentially fine spatial scales over time with only a linear increase in memory.

\subsection{A Dissipative Extension of the CM Method}
\label{section:dissipative_extension}

The proposed dissipative extension of the characteristic mapping method follows from a particular transformation of the equation \eqref{eq:main_pde} into a continuity equation of the form \eqref{eq:continuity_equation}. In particular, we apply a non-local transformation of the PDE \eqref{eq:main_pde} that can be written formally as
\begin{equation*}
    \partial_t \solution + \operatorname{div} \left( \solution \cdot \left(u -  \frac{1}{\solution} \nabla \Delta^{-1} F(\solution) \right)\right) = 0.
\end{equation*}
Here, $\Delta^{-1}$ is the inverse of the Laplacian restricted to the mean-zero subspace. From this continuity form of the equation, we can then solve for the inverse flow map of the velocity field $v= u -\frac{1}{\solution} \nabla \Delta^{-1} F(\solution)$, which retains the semigroup decomposition of the map and allow for an expression of the $\solution$ as the pullback \eqref{eq:def_pullback_in_coordinates}. \par

A simple change of variables makes this transformation well-defined. Namely, it suffices to shift $\solution$ away from zero, and to normalize $\solution$ so that its rate of change has mean zero, and hence lives in the image of the Laplacian operator. To this end, we fix a constant $c \in \mathbb{R}$ such that $\solution + c$ is bounded away from zero for all time, and define $\rho = \solution + c$. We also define the mass ratio
\begin{equation}
    \label{eq:mass_definition}
    \mass_{[t, t_0]} = \frac{\int_\manifold \rho(t) \mu}{\int_\manifold \rho(t_0) \mu},
\end{equation}
that we track separately. Given a general scalar evolution equation of the form \eqref{eq:main_pde}, we define the resulting \emph{lifting operator} as
\begin{equation*}
    \liftingOperator : C^\infty (\manifold, \mathbb{R}_+) \times [0, T] \rightarrow \mathfrak{X}(\manifold) \times \mathbb{R}
\end{equation*}
by $\liftingOperator(\rho, t) = (v, \xi)$ where
\begin{equation}
    \label{eq:lifting_operator_def}
    v = u + \frac{1}{\rho} \nabla \Delta^{-1} (\xi \rho - F(\rho - c, t)), \quad\quad \xi = \frac{\int_\manifold F(\rho - c, t) \mu}{\int_\manifold \rho \mu},
\end{equation}
and $u = A(\rho - c, t)$. Here, $v$ is the map velocity and $\xi$ is the rate of change in the mass \eqref{eq:mass_definition}. The lifting operator \eqref{eq:lifting_operator_def} forms the basis for our transformation of the original evolution equation into an equation on an extended space of diffeomorphisms $\mathcal{D}(\manifold) \times \mathbb{R}_+$. 

\begin{theorem}
    \label{thm:lift}
    The solution of the PDE \eqref{eq:main_pde} can be recovered through the combined mass-map action
    \begin{equation}
        \label{eq:mass_pullback}
        \rho(t) = \mass_{[t, t_0]} (\rho(t_0) \circ \map_{[t, t_0]}) \jacobianOperator (\map_{[t, t_0]}),
    \end{equation}
    where $(\map, \mass)_{[t, t_0]} \in \Diff \times \mathbb{R}_+$ satisfy
    \begin{equation}
        \label{eq:main_system_on_group}
        \partial_t \map_{[t, t_0]} + D\map_{[t, t_0]}(v) = 0, \quad\quad \partial_t \mass_{[t, t_0]} = \xi \mass_{[t, t_0]},
    \end{equation}
    with velocities given by the lifting operator \eqref{eq:lifting_operator_def}, $(v, \xi) = \liftingOperator(\rho(t), t)$, and $(\map, \mass)_{[t_0, t_0]} = (\mathrm{id}, 1)$.
\end{theorem}

\begin{proof}
    Inserting $\rho = \solution + c$ into the PDE \eqref{eq:main_pde} for $\solution$ and using $\operatorname{div} u = 0$, we obtain 
    \begin{equation}
        \label{eq:rho_almost_in_conservation_form}
        \partial_t \rho + \operatorname{div}(\rho u) = F(\rho - c, t).
    \end{equation}
    Integrating both sides gives $\partial_t \mass = \xi \mass$, with $\mass$ as defined in \eqref{eq:mass_definition} and $\xi$ as defined in \eqref{eq:lifting_operator_def}. Now, the quantity $\overline{\rho}(t) := \rho(t)/\mass_{[t, t_0]}$ has constant total mass, and thus can be written as the solution to a continuity equation using the manipulation described above. Indeed, differentiating $\overline{\rho} = \rho/\mass$ in time and using \eqref{eq:rho_almost_in_conservation_form} gives
    \begin{equation}
        \label{eq:rho_almost_in_conservation_form_2}
        \partial_t \overline{\rho} + \operatorname{div}(\overline{\rho} u) = \frac{1}{m} (F(\rho - c, t) - \xi \rho).
    \end{equation}
    using the divergence-free property of the velocity field $u$. The right hand side has mean zero by the definition \eqref{eq:lifting_operator_def} of $\xi$, and so we can now legitimately perform the manipulation described above: with the choice of velocity
    \begin{equation*}
        v = u + \frac{1}{m \overline{\rho}} \nabla \Delta^{-1} (\xi \rho - F(\rho - c, t)),
    \end{equation*}
    which is equivalent to \eqref{eq:lifting_operator_def}, we have that \eqref{eq:rho_almost_in_conservation_form_2} is equivalent to $\partial_t \overline{\rho} + \operatorname{div}(\overline{\rho} v) = 0$. The pullback \eqref{eq:def_pullback_in_coordinates} applied to $\overline{\rho}$ gives $\frac{\rho(t)}{\mass_{[t, t_0]}} = \left( \frac{\rho(t_0)}{\mass_{[t_0, t_0]}} \circ \map_{[t, t_0]} \right) \jacobianOperator(\map_{[t, t_0]})$, which is equivalent to \eqref{eq:mass_pullback}.
\end{proof}

\subsubsection*{Geometric Descriptions of Continuum Motion}

The ansatz taken for the transformation of the scalar evolution equation \eqref{eq:main_pde} through the lifting operator \eqref{eq:lifting_operator_def} possesses a geometric interpretation which connects our proposed method to the descriptions of continuum motion developed in geometric hydrodynamics \cite{arnold_1966, ebin1970groups, khesin2021geometric} and optimal transport \cite{mccann_convexity_1997, jordan_variational_1998, benamou_computational_2000, otto2001geometry}. In these settings, the spaces of diffeomorphisms or densities can be equipped with a (weak) Riemannian metric structure from which distances and metric gradients can be defined. These constructions have been shown to connect a wide variety of partial differential equations as geodesics, gradient flows, or Newtonian dynamics derived solely from choice of the metric and potential energy \cite{otto2001geometry, khesin2021geometric}. The action of the diffeomorphism group on a reference density relates these spaces through the geometry of a \emph{principal bundle} and the metric provides a way to lift paths and variations on the space of densities to ones of the space of diffeomorphisms.

Extensions of the framework of optimal transport on the space of probability densities have been proposed to accommodate densities with varying total mass \cite{benamou2003numerical,chizat2018unbalanced, khesin_simple_2024}. This theory provides a natural geometric framework upon which our dissipative extension of the CM method and the construction of the lifting operator are formulated \eqref{eq:lifting_operator_def}. In our approach, the path in the space of densities is prescribed by the solution to the PDE \eqref{eq:main_pde}, while the metric determines how to continuously lift this path into the extended diffeomorphism group. We then utilize this lifted representation to approximate the solution through the extended group action. A similar Lie-theoretic approach to lifting evolution equations was also proposed in \cite{taylor2026analysis}, and this technique relies primarily on transitivity of the action and the definition of a continuous lifting operator. Consequently, the dissipative extension of the CM method retains its geometric formulation without requiring the original equation to arise from a variational formulation or an optimal transport problem. We develop this geometric interpretation further in Section~\ref{section:geometric_framework}.

\subsubsection*{Semidirect Product Evolution Technique}

Using the reformulation of the PDE \eqref{eq:main_pde} as the lift in Theorem \ref{thm:lift}, we can compute the map and exploit the semigroup decomposition. However, the direct evolution of the map through \eqref{eq:main_system_on_group} yields a system with a stiff nonlinearity through the evaluation of the Jacobian $\map \mapsto \jacobianOperator(\map)$, which presents a numerical instability. Indeed, we find that the \emph{transport component} $\rho_0 \circ \map$ and \emph{compression component} $\jacobianOperator(\map)$ of the solution require a distinct numerical treatment. Thus, we promote the log-Jacobian $\jacobian := \log \jacobianOperator(\map)$ of the map to a primary variable for discretization. The extended state $(\map, \jacobian, \mass)$ forms an element of the semidirect product group $\Diff \ltimes C^\infty(\manifold) \times \mathbb{R}_+$, and we collect the necessary ingredients for the explicit evolution of $(\map, \jacobian, \mass)$ in the context of the semigroup decomposition.

We then describe an algorithm for the solution of PDEs on the torus $\mathbb{T}^2$ with a leading-order linear dissipative operator $\mathcal{L}$ and a lower-order nonlinearity $G$,
\begin{equation*}
    \partial_t \solution + u \cdot \nabla \solution = \mathcal{L} \solution + G(\solution, t).
\end{equation*}
In this case, the log-Jacobian $\jacobian$ satisfies an evolution equation with the same leading-order term, $\partial_t \jacobian = \mathcal{L} \jacobian + \mathcal{N}(\jacobian, t)$, lending naturally to a stabilized time integration for $\jacobian$ via an exponential integrator. This splitting then enables the map $\map$ to be stably advanced by the usual semi-Lagrangian strategy of the CM method. Despite the distinct contributions to the solution of the transport and compression components, these combine through the semigroup decomposition and the numerical group action in a manner which retains the resolution benefits of the CM method. This is featured for the Navier--Stokes equations in Figure \ref{fig:navier_stokes_components_intro}.

\newcommand{\intronavierstokescomponentsimagewidthpercent}{0.29}
\begin{figure}[H]
    \centering
    \begin{subfigure}[c]{\intronavierstokescomponentsimagewidthpercent\textwidth}
        \includegraphics[width=\textwidth]{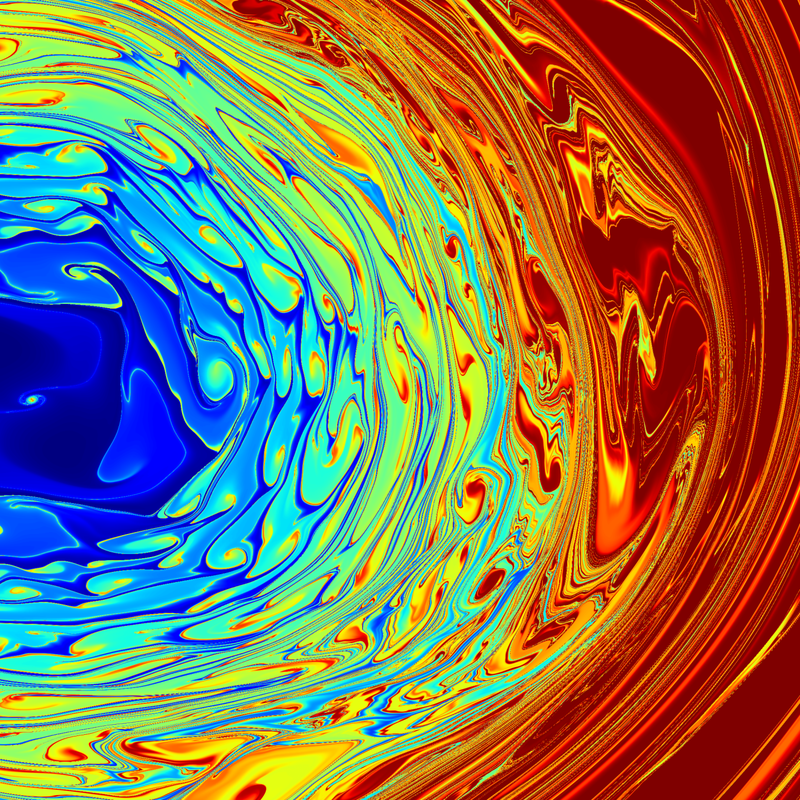}
    \end{subfigure}
    \begin{minipage}[c]{0.045\textwidth}
        \centering
        \Huge $\times$
    \end{minipage}
    \begin{subfigure}[c]{\intronavierstokescomponentsimagewidthpercent\textwidth}
        \includegraphics[width=\textwidth]{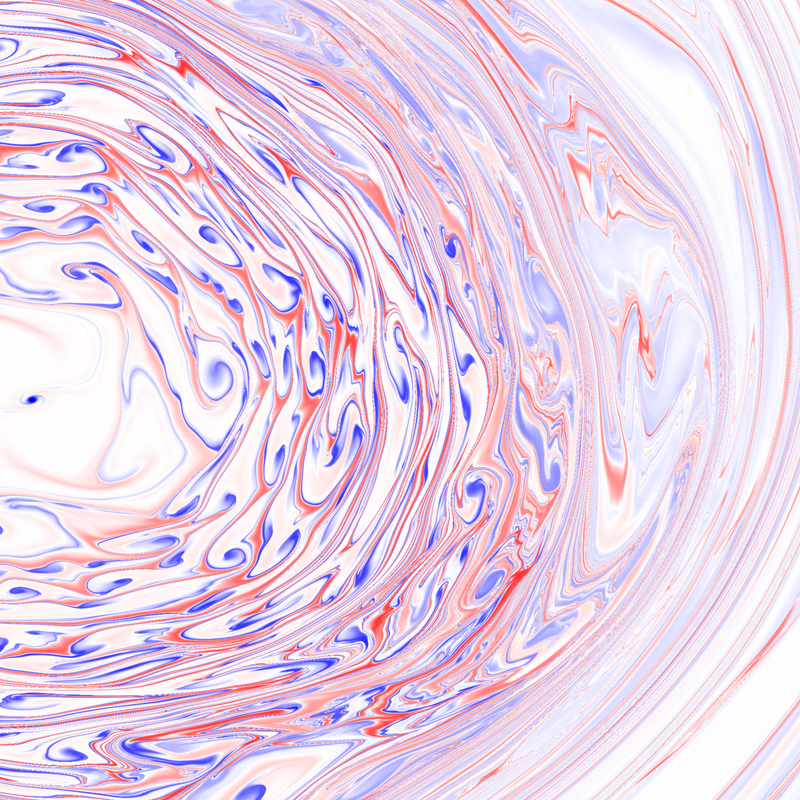}
    \end{subfigure}
    \begin{minipage}[c]{0.045\textwidth}
        \centering
        \Huge $=$
    \end{minipage}
    \begin{subfigure}[c]{\intronavierstokescomponentsimagewidthpercent\textwidth}
        \includegraphics[width=\textwidth]{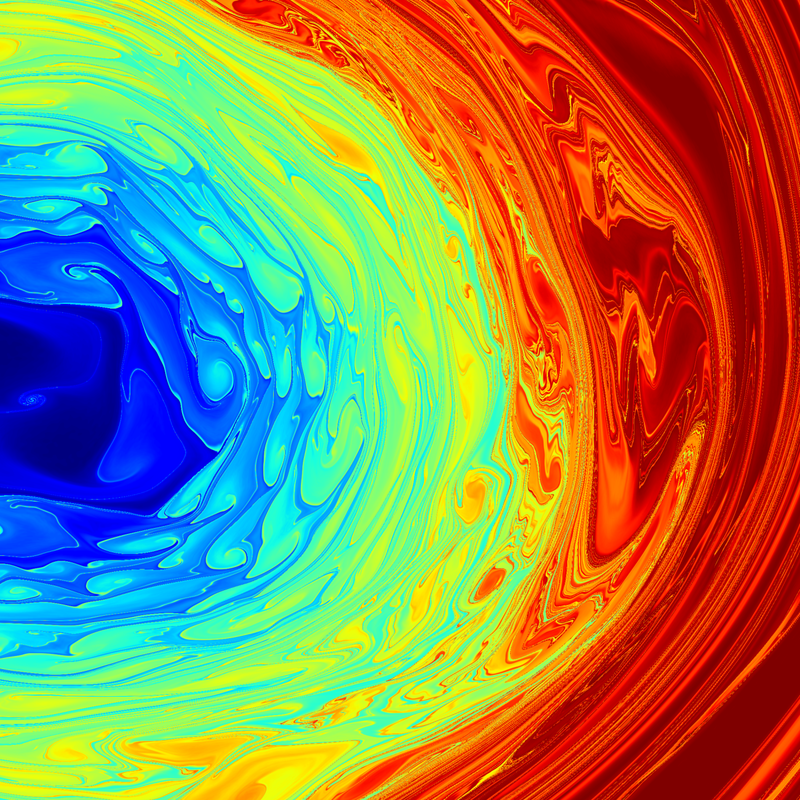}
    \end{subfigure}
    \\[0.3em]
    \makebox[\intronavierstokescomponentsimagewidthpercent\textwidth]{\parbox{\intronavierstokescomponentsimagewidthpercent\textwidth}{\centering $\rho_0 \circ \map$}}
    \hspace*{0.045\textwidth}
    \makebox[\intronavierstokescomponentsimagewidthpercent\textwidth]{\parbox{\intronavierstokescomponentsimagewidthpercent\textwidth}{\centering $\jacobianOperator(\map)$}}
    \hspace*{0.045\textwidth}
    \makebox[\intronavierstokescomponentsimagewidthpercent\textwidth]{\parbox{\intronavierstokescomponentsimagewidthpercent\textwidth}{\centering $\rho$}}
    \caption{Decomposition of the vorticity $\rho$ in the Navier--Stokes vortex merger, $\nu = 10^{-5.5}$, at the final time into the different map components: the transport component $\rho_0 \circ \map$ and the compression component $\jacobianOperator(\map)$.}
    \label{fig:navier_stokes_components_intro}
\end{figure}

\subsection{Related Numerical Works}
\label{section:related_works}

There are several works, particularly in the graphics community, focused on extending traditional map-based methods to problems in hydrodynamics that lack an explicit transport structure.
Some earlier works \cite{tessendorf_characteristic_2011, qu_efficient_2019}, as well as the active body of work on the Particle Flow Maps method \cite{li_impulse-formed_2026, wang_eulerian_2024, li_clebsch_2025, chen_fluid_2025, sun_impulse_2024}, adopt a hybrid map-integral approach. A similar approach is taken in the CM method with source terms \cite{yin_characteristic_2026} where a flow map is computed to account for the transport term in the equation, whereas the contributions of other terms (e.g. viscous and boundary forces) are accumulated through a path integral. In the present work, we instead write the entire evolution as a single action by a map. This systematic approach avoids having to first identify a transport-dominated quantity in the problem and is naturally compatible with a semigroup decomposition, ensuring that the approximation space grows with each submap.

There are also a number of Lagrangian methods for Wasserstein and other gradient flows that have been proposed based on a lifting from the space of densities lifts to a flow on a space of maps \cite{carrillo_chapter_2021}. Utilizing a regularized energy yields a deterministic particle evolution, giving rise to a number of particle \cite{gallouet_convergence_2022, leclerc_lagrangian_2020} and blob \cite{carrillo_blob_2019} methods. Other techniques explicitly approximate the forward flow map and its Jacobian. One family of methods \cite{carrillo_numerical_2010, carrillo_numerical_2016} derives the map evolution equation from a continuous variational formulation before discretizing. Another family \cite{junge_fully_2017, liu_lagrangian_2020, carrillo_lagrangian_2018} first discretizes a variational principle in the map variables, then advances these in time via incremental minimization, in the style of the JKO scheme \cite{jordan_variational_1998}. We note further, that machine learning approaches have also been developed, see for instance \cite{lee_deep_2024}, where the map velocity is parametrized by a neural network. In the present work, we will consider the inverse flow map, as it naturally lives in a fixed Eulerian frame and derive explicit evolution equations for the map and its Jacobian rather than considering a variational formulation. Contrary to the above methods, we consider the geometry of the problem purely for the purpose of lifting the evolution to the diffeomorphism group, where the semigroup decomposition is available. We note that the variational approaches listed above lend naturally to the case of a vanishing density, which we do not consider here.

\subsection{Outline}

This paper is organized as follows. We begin in Section \ref{section:geometric_framework} by giving a geometric interpretation of the lift constructed in Theorem \ref{thm:lift} as a certain horizontal lift in an unbalanced optimal transport setting, after reviewing the relevant geometric notions. This formalization can safely be skipped by readers primarily interested in the numerical aspects of this work. In Section \ref{section:semidirect_product_evolution}, we promote the log-Jacobian of the map to a primary variable for discretization, and introduce the necessary ingredients for its explicit evolution in the context of a semigroup decomposition. In Section \ref{section:implementation_analysis}, we describe an algorithm for the solution of a certain class of PDEs on the torus, and give an asymptotic error estimate for the scheme. Finally, in Section \ref{section:examples}, we consider a number of dissipative equations, compute their lifts, and perform experiments highlighting the roles of the transport and compression components of the map. We also validate the error analysis with the method of manufactured solutions.

\section{Geometric Formulation}
\label{section:geometric_framework}

In this section we provide a geometric interpretation of the extension of the characteristic mapping method to dissipative evolution equations using the theory of unbalanced optimal transport \cite{khesin_simple_2024}. We begin by first motivating the transformation of a general evolution equation into a lifted differential equation on an extension of the diffeomorphism group. We then review the notion of a horizontal lift in a principal bundle, giving the canonical example of the Wasserstein--Otto geometry. We then show that the extension of the CM method described in Theorem \ref{thm:lift} arises precisely as a horizontal lift in this geometry with a particular metric structure. \par

\subsection{Preliminaries}

We will work formally with a differential structure on the space of diffeomorphisms $\Diff$ of a closed Riemannian manifold $(\manifold, g)$. This structure can be made precise in a Fr{\'e}chet topology in the smooth setting or with a metric topology in the Sobolev setting \cite{schmeding2022introduction, ebin1970groups}. We will consider this differential structure only formally, in order to emphasize the construction. A rigorous functional-analytic treatment of this formulation is of interest, yet beyond the scope of the current work. \par

Fix the Riemannian volume form $\mu$ on $\manifold$, and consider the space of all positive volume forms
\begin{equation*}
    \Vol = \left\{ \rho \mu : \rho \in C^\infty (\manifold, \mathbb{R}_+) \right\}
\end{equation*}
as well as the space of normalized densities
\begin{equation*}
    \Prob = \left\{ \varrho \in \mathcal{V}(\manifold): \int_\manifold \varrho = \int_\manifold \mu \right\}.
\end{equation*}
The diffeomorphism group acts on $\Prob$ transitively \cite{moser1965volume} by the pushforward 
\begin{equation}
\label{eq:pullack_on_densities}
\Diff \times \Prob \rightarrow \Prob, \quad\quad (\forwardMap, \varrho) \mapsto \forwardMap_\ast \varrho = \rho \circ \forwardMap^{-1} \jacobianOperator(\forwardMap^{-1}) \mu,
\end{equation}
where $\varrho = \rho \mu \in \Prob$, and $\jacobianOperator$ is the Jacobian determinant with respect to $\mu$. The isotropy group of this action is the space of volume-preserving diffeomorphisms $\SDiff$. \par

In our dissipative extension of the CM method, we consider a central extension of the diffeomorphism group, defined by the Lie group $(\productGrp, \cdot)$ where
\begin{equation*}
    \productGrp = \Diff \times \mathbb{R}_+, \quad\quad (\forwardMap_1, \mass_1) \cdot (\forwardMap_2, \mass_2) = (\forwardMap_1 \circ \forwardMap_2, \mass_1 \mass_2),
\end{equation*}
This group acts on all the space of volume forms via
\begin{equation}
\label{extended_group_action}
    \productGrp \times \Vol \rightarrow \Vol, \quad\quad ((\forwardMap, \mass), \varrho) \mapsto \mass \forwardMap_\ast \varrho.
\end{equation}
This action is transitive, and forms the basis for a simple unbalanced optimal transport \cite{khesin_simple_2024}, with the component $\mass$ serving to keep track of a changing mass. We will assemble group elements into the boldface notation $\fwdGrpElt = (\forwardMap, \mass) \in \productGrp$, and denote elements of the tangent space as
\begin{equation*}
    (v, \xi) \cdot \fwdGrpElt = (v \circ \forwardMap, \xi \mass) \in T_{\fwdGrpElt} \productGrp = T_\forwardMap \Diff \times T_\mass \mathbb{R}_+.
\end{equation*}

The orbit map of the Riemannian volume form $\mu$ relative to this action is defined by
\begin{equation}
    \label{eq:orbit_map_def}
    \pi : \productGrp \rightarrow \Vol, \quad\quad \pi(\fwdGrpElt) = \mass \forwardMap_\ast \mu.
\end{equation}
This map defines a surjection of the extended group onto the space of volume forms due to the transitivity of the action \eqref{extended_group_action}, forming the basis for the lifted evolution equation on $\productGrp$. Its differential $d\pi : T\productGrp \rightarrow T\Vol$ is defined pointwise on the tangent spaces as
\begin{equation}
    \label{eq:orbit_map_differential}
    d \pi_{\fwdGrpElt} : T_{\fwdGrpElt} \productGrp \rightarrow T_{\pi(\fwdGrpElt)} \Vol, \quad\quad (v, \xi) \cdot \fwdGrpElt \mapsto (\xi \rho - \operatorname{div}(\rho v)) \mu,
\end{equation}
see \cite[Lem. 2.1]{khesin_simple_2024}. The isotropy group of the orbit map is the subgroup $\isoSubGp = \pi^{-1} (\mu) = \SDiff \times \{ 1 \}$. The map \eqref{eq:orbit_map_def} gives $\productGrp$ the structure of a \emph{principal $\isoSubGp$-bundle} over $\Vol$, a structure which is central to the lifting construction. 

\subsection{Construction of the Horizontal Lift}
\label{sec:horizontal_lift}

\begin{figure}
    \centering
    \begin{subfigure}[t]{0.58\textwidth}
        %\asyinclude[width=\the\linewidth,keepAspect=true,inline=true]{principal_bundle_figures/pb_figure_1.asy}
        \includegraphics[width=277.4390bp,height=183.5600bp]{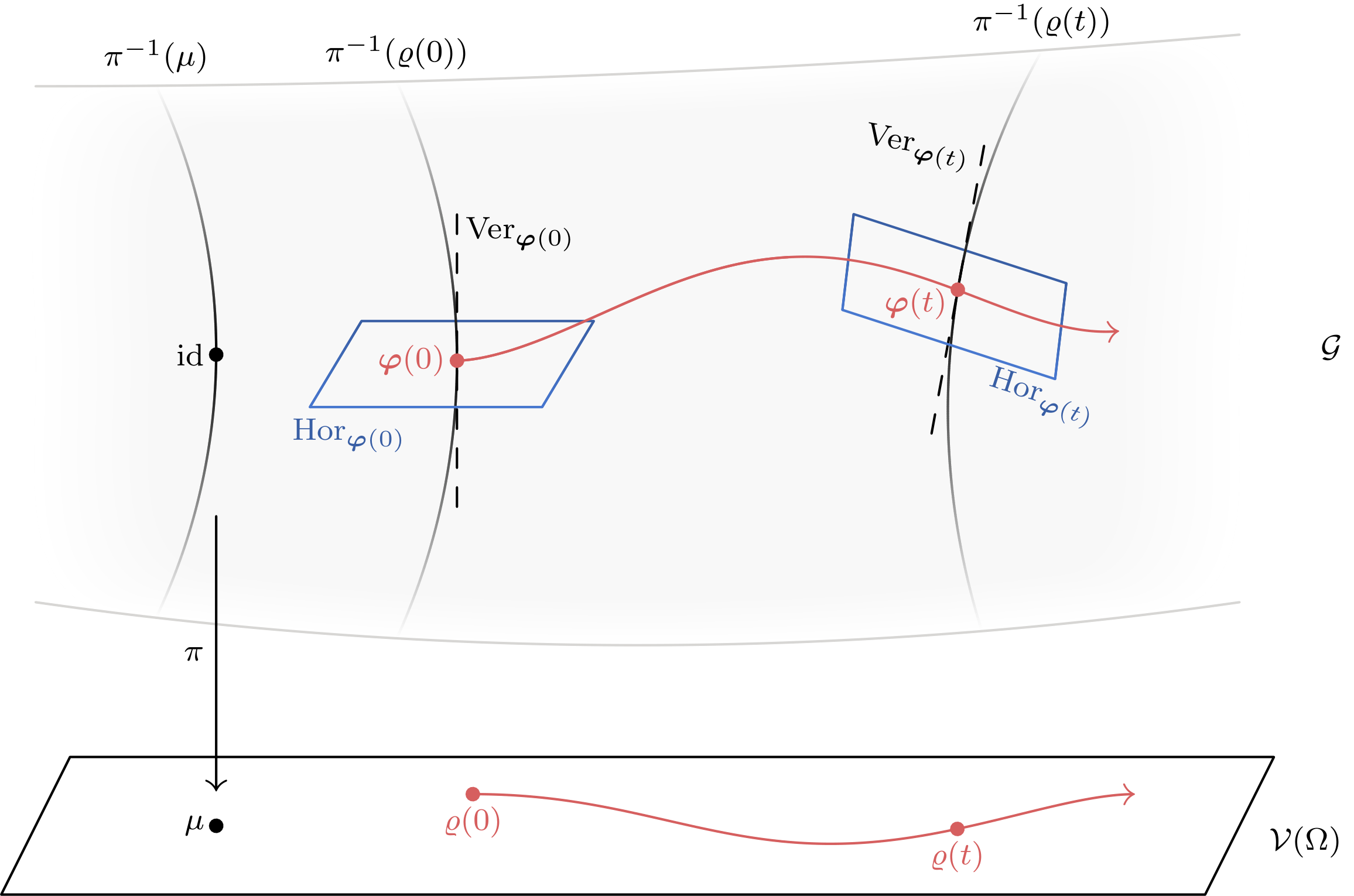}
        \caption{Orthogonal decomposition of tangent spaces $T_{\fwdGrpElt} \mathcal{G} = \Ver_{\fwdGrpElt} \oplus \Hor_{\fwdGrpElt}$, and the horizontal lift $\fwdGrpElt(t) \in \mathcal{G}$ of a curve $\varrho(t) \in \Vol$.}
        \label{fig:pb_figure_1}
    \end{subfigure}
    \hfill
    \begin{subfigure}[t]{0.36\textwidth}
        %\asyinclude[width=\the\linewidth,keepAspect=true,inline=true]{principal_bundle_figures/pb_figure_2.asy}
        \includegraphics[width=174.4800bp,height=187.2600bp]{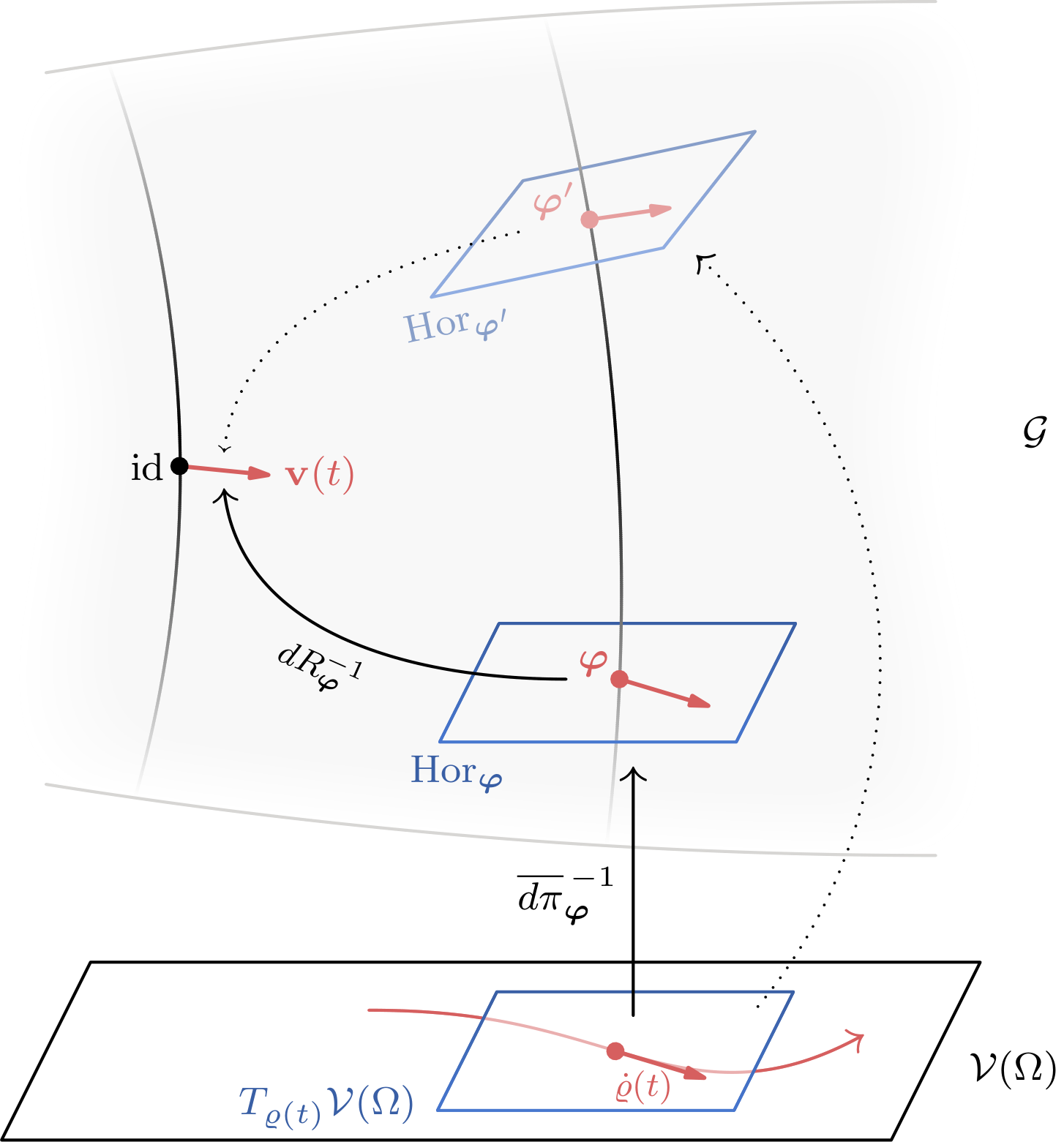}
        \caption{Conditions \ref{item:metric_assumption_complement} and \ref{item:metric_assumption_equivariance} imply that a tangent vector $\dot{\varrho}$ downstairs lifts to a well-defined trivialized tangent vector upstairs.}
        \label{fig:pb_figure_2}
    \end{subfigure}
    \caption{Illustrations of the principal bundle structure $\pi : \mathcal{G} \rightarrow \Vol$ and the horizontal lifting construction.}
    \label{fig:pb_figure}
\end{figure}

Let the path $\phi(t) \in C^\infty(\manifold)$ describe the evolution of a scalar field, such as the solution to the equation \eqref{eq:main_pde}. Using a transformation of the form
\begin{equation*}
    \varrho(t) = \rho(t) \mu := (\solution(t) + c) \mu \in \Vol \,,
\end{equation*}
where $c$ is taken such that $\solution + c > 0$ in $\manifold \times [0, T]$, surjectivity of the orbit map \eqref{eq:orbit_map_def} guarantees the existence of a path in the transformation group $\fwdGrpElt(t) \in \productGrp$ descending to $\varrho(t)$ via
\begin{equation*}
    \pi(\fwdGrpElt(t)) = \varrho(t), \quad\quad t \in [0, T].
\end{equation*}
Uniqueness is obstructed by the kernel of the differential of the orbit map \eqref{eq:orbit_map_differential}, referred to as the \emph{vertical distribution}
\begin{equation*}
    \Ver_{\fwdGrpElt} := \ker d \pi_{\fwdGrpElt} \subset T_{\fwdGrpElt} \mathcal{G}.
\end{equation*}
The construction of a smooth lift necessitates a choice of gauge that allows us to continuously complement the vertical distribution while remaining compatible with the principal bundle structure. A canonical construction to this end is to introduce a (weak) Riemannian metric on the group $\productGrp$,
\begin{equation*}
    \metric_{\fwdGrpElt} : T_{\fwdGrpElt}  \mathcal{G} \times T_{\fwdGrpElt} \mathcal{G} \to \mathbb{R}.
\end{equation*}
We can then define a \emph{horizontal distribution} as the orthogonal complement
\begin{equation*}
    \Hor_{\fwdGrpElt} := \Ver_{\fwdGrpElt}^\perp 
\end{equation*}
with respect to the metric. A uniquely defined inverse to the differential is then obtained by restricting to lifts which are \emph{horizontal,} meaning $\dot{\fwdGrpElt}(t) \in \Hor_{\fwdGrpElt(t)}$ for all $t \in [0, T]$. A horizontal lift is illustrated in Figure \ref{fig:pb_figure_1}.

The metric should be chosen such that the two following conditions hold:
\begin{enumerate}[noitemsep, topsep=0.1cm, label=(\roman*)]
    \item\label{item:metric_assumption_complement} The horizontal distribution complements the vertical distribution, in the sense
    \begin{equation*}
        T_{\fwdGrpElt} \mathcal{G} = \Ver_{\fwdGrpElt} \oplus \Hor_{\fwdGrpElt}.
    \end{equation*}
    This condition is not generically true for orthogonal complements in infinite dimensions.

    \item\label{item:metric_assumption_equivariance} The horizontal distribution satisfies an equivariance condition on the isotropy subgroup,
    \begin{equation}
        \label{eq:equivariance_condition}
        \Hor_{\fwdGrpElt \cdot \boldsymbol{\eta}} = d R_{\boldsymbol{\eta}} (\Hor_{\fwdGrpElt}), \quad\quad \text{for all}\ \boldsymbol{\eta} \in \isoSubGp.
    \end{equation}
    This condition naturally leads to the consideration of \emph{semi-invariant} Riemannian metrics (c.f. \cite{bauer_semi-invariant_2020}), meaning those that are right-invariant only on the isotropy subgroup $\isoSubGp$. We note that this holds automatically when the metric depends only on the basepoint $\fwdGrpElt$ through a function of $\rho \mu = \pi(\fwdGrpElt)$, as $\pi$ is invariant on the fibres.
\end{enumerate}
When these two conditions are satisfied, the following theorem holds as particular case of \cite[§39.1]{kriegl1997convenient}.
\begin{theorem}
    \label{thm:generic_lifting_theorem}
    If conditions \ref{item:metric_assumption_complement} and \ref{item:metric_assumption_equivariance} above are satisfied, then curves $\varrho(t) \in \Vol$ lift to unique horizontal curves $\fwdGrpElt(t) \in \mathcal{G}$ after an initial choice in the fibre $\fwdGrpElt(0) \in \pi^{-1}(\varrho(0))$.
\end{theorem}

The roles of conditions \ref{item:metric_assumption_complement} and \ref{item:metric_assumption_equivariance} in constructing a lift can be interpreted as follows. First, condition \ref{item:metric_assumption_complement} implies that the differential $d \pi$ of the orbit map restricts to an isomorphism on the horizontal distribution,
\begin{equation*}
    \overline{d \pi}_{\fwdGrpElt} := d \pi_{\fwdGrpElt}|_{\Hor_{\fwdGrpElt}} : \Hor_{\fwdGrpElt} \overset{\sim}{\longrightarrow} T_{\pi(\fwdGrpElt)} \Vol.
\end{equation*}
Thus, given a curve $\varrho(t) \in \Vol$, for each $t$ we can lift the tangent vector $\dot{\varrho}(t) \in T_{\varrho(t)} \Vol$ downstairs to a tangent vector $(v(t), \xi(t)) \cdot \fwdGrpElt \in T_{\fwdGrpElt} \productGrp$ upstairs through
\begin{equation*}
    (v(t), \xi(t)) \cdot \fwdGrpElt = \overline{d \pi}_{\fwdGrpElt}^{\, -1} (\dot{\varrho}(t)),
\end{equation*}
after making a choice of element in the fibre $\fwdGrpElt \in \pi^{-1}(\varrho(t))$. Condition \ref{item:metric_assumption_equivariance} then implies that the right-trivialization
\begin{equation}
    \label{eq:right_trivialized_horizontal_lift_well_defined}
    (v(t), \xi(t)) = dR_{\fwdGrpElt}^{-1} \circ \overline{d \pi}_{\fwdGrpElt}^{\,-1}(\dot{\varrho}(t))
\end{equation}
is independent of the choice $\fwdGrpElt \in \pi^{-1}(\varrho)$ in the fibre. We can therefore directly construct the horizontal lift $\fwdGrpElt(t)$ of $\varrho(t)$ as the flow of the $(v(t), \xi(t))$ given by \eqref{eq:right_trivialized_horizontal_lift_well_defined}, after a single choice  $\fwdGrpElt(0) \in \pi^{-1}(\varrho(0))$, without having to make consistent choices of elements in the fibre $\pi^{-1}(\varrho(t))$ for each $t > 0$. This discussion is summarized in Figure \ref{fig:pb_figure_2}.

\subsubsection*{Example: Wasserstein-Otto Geometry}

As an example construction of the horizontal lift, we consider Otto's Riemannian geometric interpretation of dissipative evolution equations \cite{otto2001geometry}. In this setting, we consider the orbit map
\begin{equation*}
    \pi : \Diff \to \Prob, \quad\quad \forwardMap \mapsto \forwardMap_* \mu,
\end{equation*}
where $\mu \in \Prob$. The isotropy group is $\SDiff$ and the vertical distribution has the expression
\begin{equation*}
    \Ver_\varphi := \ker \left( d \pi_\mu \vert_\forwardMap \right) = \left\{ v \circ \forwardMap \in T_\forwardMap \Diff : \operatorname{div}(\rho v ) = 0  \right\}.
\end{equation*}
Otto equipped $\Diff$ with the semi-invariant $L^2$ metric,
\begin{equation}
    \label{eq:otto_metric}
    \metric_\forwardMap (v \circ \forwardMap, v \circ \forwardMap) = \int_\manifold \rho g(v, v) \mu.
\end{equation}
One finds the induced horizontal distribution to consist of the gradient fields
\begin{equation}
    \label{eq:otto_horizontal_dist}
    \Hor_\forwardMap = \nabla C^\infty (\manifold) \circ \forwardMap.
\end{equation}
Condition \ref{item:metric_assumption_complement} is equivalent to the Helmholtz--Hodge decomposition. Condition \ref{item:metric_assumption_equivariance} holds by semi-invariance of the metric, which itself follows from the fact that the expression \eqref{eq:otto_metric} only depends on the basepoint $\forwardMap$ through $\rho = \jacobianOperator(\forwardMap)$. Restricting the differential of the orbit map to \eqref{eq:otto_horizontal_dist} yields a smooth inverse through the solution of a variable-coefficient elliptic problem, defining the lift
\begin{equation}
    \label{eq:otto_lift}
    \dot{\forwardMap} \circ \forwardMap^{-1} = \nabla \psi, \quad\quad -\operatorname{div}(\rho \nabla \psi) = \dot{\rho}.
\end{equation}

\subsection{An Unbalanced, Linearized Wasserstein Geometry}

In our extension of the CM method, we require the construction of a horizontal lift in the extended group $\productGrp = \Diff \times \mathbb{R}_+$. The differential of the orbit map \eqref{eq:orbit_map_differential} has the vertical distribution (c.f. \cite[Lem. 2.1]{khesin_simple_2024})
\begin{equation}
    \label{eq:our_vertical_distribution}
    \Ver_{\fwdGrpElt} = \left\{ (u, 0) \cdot \fwdGrpElt: u \in \mathfrak{X}(\manifold), \ \operatorname{div}(\rho u) = 0 \right\}.
\end{equation}
Our proposed semi-invariant Riemannian metric on $\productGrp$ is defined as
\begin{equation}
    \label{eq:our_metric}
    \metric_{\fwdGrpElt} ((v, \xi) \cdot \fwdGrpElt, (v, \xi) \cdot \fwdGrpElt) = \int_\manifold \rho^2 g(v, v) \mu + \xi^2.
\end{equation}

\begin{proposition}
    \label{prop:our_horizontal_dist_properties}
    The horizontal distribution induced by the metric \eqref{eq:our_metric} satisfies conditions \ref{item:metric_assumption_complement} and \ref{item:metric_assumption_equivariance} above, and has the expression
    \begin{equation}
        \label{eq:our_horizontal_dist}
        \Hor_{\fwdGrpElt} = \left( \frac{1}{\rho} \nabla C^\infty (\manifold) \right) \circ \forwardMap \oplus \mathbb{R} \mass.
    \end{equation}
\end{proposition}

\begin{proof}
    Let $(v, \xi) \cdot \fwdGrpElt \in \Ver_{\fwdGrpElt}^\perp$. Then for all $(u, 0) \cdot \fwdGrpElt \in \Ver_{\fwdGrpElt}$ (recall \eqref{eq:our_vertical_distribution}), we have
    \begin{equation*}
        \metric_{\fwdGrpElt} ((v, \mass) \cdot \fwdGrpElt, (u, 0) \cdot \fwdGrpElt) = \int_\manifold g(\rho v, \rho u) \mu = 0.
    \end{equation*}
    As this holds for all fields $\rho u$ with $\operatorname{div}(\rho u) = 0$, the flat $L^2$ Helmholtz decomposition gives $\rho v = \nabla \psi$ for some $\psi \in C^\infty (\manifold)$. Conversely, if $v$ has the form $v = \frac{1}{\rho} \nabla \psi$, then for all $(u, 0) \cdot \fwdGrpElt \in \Ver_{\fwdGrpElt}$ we have
    \begin{equation*}
        \metric_{\fwdGrpElt} ((v, \xi) \cdot \fwdGrpElt, (u, 0) \cdot \fwdGrpElt) = \int_\manifold g(\nabla \psi, \rho u) \mu = -\int_\manifold \psi \operatorname{div}(\rho u) \mu = 0.
    \end{equation*}
    This establishes the expression for the horizontal distribution \eqref{eq:our_horizontal_dist}.

    We further have that an arbitrary $(v, \xi) \cdot \fwdGrpElt \in T_{\fwdGrpElt} \mathcal{G}$ may be decomposed orthogonally as
    \begin{equation*}
        v = u + \frac{1}{\rho} \nabla \psi \in \Ver_{\fwdGrpElt} \oplus \Hor_{\fwdGrpElt},
    \end{equation*}
    where $\nabla \psi$ is determined by $\Delta \psi = \operatorname{div}(\rho v)$, and then setting $u = v - \frac{1}{\rho} \nabla \psi \in \Ver_{\fwdGrpElt}$. This establishes that $\Ver$ and $\Hor$ together generate the entire tangent space, i.e. shows condition \ref{item:metric_assumption_complement}. The equivariance condition \ref{item:metric_assumption_equivariance} holds from the semi-invariance of the metric, which itself holds again because the metric depends on the basepoint $\fwdGrpElt$ only through $\rho \mu = \pi(\fwdGrpElt)$.
\end{proof}

\begin{proposition}
    \label{prop:our_lift_characterization}
    The unique horizontal lift $\fwdGrpElt(t)$ of the curve $\varrho(t) = \rho(t) \mu$ induced by the metric \eqref{eq:our_metric} is given by the solution to $\dot{\fwdGrpElt} = (v, \xi) \cdot \fwdGrpElt$ where
    \begin{equation*}
        v(t) = \frac{1}{\rho(t)} \nabla \Delta^{-1} \left( \xi(t) \rho(t) - \dot{\rho}(t) \right), \quad\quad \xi(t) = \frac{\int_\manifold \dot{\rho}(t)\mu}{\int_\manifold \rho(t) \mu},
    \end{equation*}
    after an initial choice in the fibre $\fwdGrpElt(0) \in \pi^{-1} (\varrho(0))$.
\end{proposition}

\begin{proof}
Using the characterization of the horizontal distribution \eqref{eq:our_horizontal_dist} we can write $\rho v = \nabla \psi$. Inserting $(v, \xi)$ into the expression for the differential \eqref{eq:orbit_map_differential} gives
\begin{equation}
    \label{eq:inserted_v_xi_into_differential}
    \dot{\rho} \mu = d \pi_{\fwdGrpElt} ((v, \xi) \cdot \fwdGrpElt) = (\xi \rho - \Delta \psi) \mu.
\end{equation}
Integrating \eqref{eq:inserted_v_xi_into_differential} over the domain and using that $\Delta \psi$ has mean zero gives the expression for $\xi$. It then follows that $(\xi \rho - \dot{\rho}) \in \mathrm{im}(\Delta)$ and there exists a unique solution to
\begin{equation*}
\Delta \psi = \xi \rho - \dot{\rho} \quad\implies\quad v = \frac{1}{\rho} \nabla \psi = \frac{1}{\rho} \nabla \Delta^{-1} (\xi \rho - \dot{\rho}),
\end{equation*}
establishing the claim.
\end{proof}

This proposition leads us to the following geometric interpretation of Theorem \ref{thm:lift}.

\begin{theorem}
    \label{eq:section_2_conclusion_theorem}
    Our reformulation of the PDE $\partial_t \solution + u \cdot \nabla \solution = F(\solution, t)$ in Theorem \ref{thm:lift} is equivalent to embedding the solution into the space of volume forms and lifting the evolution to the extended diffeomorphism group $\Diff \times \mathbb{R}_+$. In the lift $(v, \xi) = \liftingOperator(\rho, t)$ of Theorem \ref{thm:lift}, the velocity
    \begin{equation}
        \label{eq:lift_restated_endofsection2}
        v = u + \frac{1}{\rho} \nabla \Delta^{-1} (\xi \rho - F(\rho - c, t))
    \end{equation}
    is precisely the horizontal lift of the dissipative term $F$ in the metric \eqref{eq:our_metric}, translated by the incompressible velocity $u$.
\end{theorem}

We have explicitly made two choices in this construction. First, we consider a different metric than the standard semi-invariant $L^2$ metric \eqref{eq:otto_metric} defining the Otto lift. Second, we keep the incompressible velocity $u$ as-is, without passing through a bundle projection $d \pi^{-1}$. These choices are numerically-motivated. Indeed, had we chosen to employ the Otto lift, the velocity would be given by $v = \nabla \psi$ where $\nabla \psi$ is determined (see equation \eqref{eq:otto_lift}) by
\begin{equation*}
    -\operatorname{div}(\rho \nabla \psi) = -\operatorname{div}(\rho u) + F(\rho - c, t).
\end{equation*}
We consider this choice suboptimal from a numerical point of view: it requires inverting the variable-coefficient operator $\psi \mapsto -\operatorname{div}(\rho \nabla \psi)$, and it projects the velocity $u$ onto a space of gradients, killing incompressibility. In our formulation, the horizontal component requires inverting only a constant-coefficient elliptic operator, and the incompressible field is maintained. We showed in Theorem \ref{eq:section_2_conclusion_theorem} that these requirements can still be met in the standard optimal transport framework.

\begin{remark}
    The lift \eqref{eq:lift_restated_endofsection2} is not in general horizontal when $u \neq 0$. However, it is still given by a horizontal lift \emph{after an incompressible change of coordinates}, which can be seen as follows. In the balanced case $\xi = 0$, let $\varrho(t) = (\solution(t) + c) \mu$ and let $\eta(t) \in \SDiff$ be the flow of the incompressible velocity $u$. One can show
    \begin{equation*}
        \eta_\ast (\partial_t \eta^\ast \varrho) = \partial_t \varrho + \mathcal{L}_u \varrho = F(\rho - c, t) \mu.
    \end{equation*}
    Thus, the compression component $-\frac{1}{\rho} \nabla \Delta^{-1} F(\rho - c, t)$ of the lift arises precisely as the horizontal lift of the conjugation $\eta_\ast (\partial_t \eta^\ast \varrho)$, that is, a horizontal lift after an incompressible change of frame.
\end{remark}

\begin{remark}
    \label{remark:descending_metric}
    The equivariance condition \ref{item:metric_assumption_equivariance} on the isotropy subgroup also implies, by \cite[Prop. 4.2]{Modin2015}, that $\metric$ descends to a metric $\overline{\metric}$ on $\Vol$, making $\pi$ into a Riemannian submersion, as is the case for the $L^2$ metric on $\Diff$ descending to the Wasserstein geometry on $\Prob$ \cite{otto2001geometry}. In this case one can compute the expression for $\overline{\metric}$ to be
    \begin{equation*}
        \overline{\metric}_\varrho (\dot{\varrho}, \dot{\varrho}) = \metric_{\fwdGrpElt} (\dot{\fwdGrpElt}, \dot{\fwdGrpElt}) = \frac{1}{\mass^2} \int_\manifold g(\nabla \Delta^{-1} (\xi \rho - \dot{\rho}), \nabla \Delta^{-1} (\xi \rho - \dot{\rho})) \mu + \xi^2
    \end{equation*}
    for an arbitrary $\dot{\varrho} \in T_\varrho \Vol$, where $\fwdGrpElt \in \pi^{-1}(\varrho)$ and $\dot{\fwdGrpElt}$ is the horizontal lift of $\dot{\varrho}$. In the balanced case $\xi = 0$, one recognizes here the flat $H^{-1}$ metric on $\Prob$, justifying the name ``linearized Wasserstein geometry,''. Further, since $\pi$ is a Riemannian submersion, a horizontal $\dot{\fwdGrpElt}$ minimizes $\lVert \dot{\fwdGrpElt} \rVert^2$ among all vectors descending to $\dot{\varrho}$. Hence in the case $u = 0$, our construction is pointwise minimal with respect to the metric among all lifts.
\end{remark}

\section{Semidirect Product Evolution Technique}
\label{section:semidirect_product_evolution}

The stable evolution of the lifted dynamical system on $\Diff \times \mathbb{R}_+$ described in Theorem \ref{thm:lift} requires a specialized numerical integration strategy that we develop in this section. To this end, we promote the log-Jacobian of the map to be a primary variable for discretization. In this section we motivate this splitting, and describe the resulting semidirect product structure. We then write an explicit evolution of the log-Jacobian, specializing to PDEs \eqref{eq:main_pde} of the form 
\begin{equation}
    \label{eq:pde_with_highest_order_linear_term}
    \partial_t \solution + u \cdot \nabla \solution = \mathcal{L} \solution + G(\solution, t), \quad\quad u = A(\solution, t),
\end{equation}
where $\mathcal{L} = \alpha \Delta$ or $\mathcal{L} = -\alpha \Delta^2$, and $G(\cdot, t) : C^\infty (\manifold) \rightarrow C^\infty (\manifold)$ and $A(\cdot, t) : C^\infty (\manifold) \rightarrow \mathfrak{X}_\mu (\manifold)$ are differential operators of a lower order than $\mathcal{L}$.

\subsection{Numerical Approximation of the Group Action}

Many of the numerical benefits of the CM method stem from approximating the backward characteristic map $\map \in \Diff$ and evaluating its action on the solution analytically. In the balanced case $m \equiv 1$, the density is reconstructed from the backward flow as a multiplication of two components
\begin{equation}
    \rho(t_n) = (\underbrace{\rho_0 \circ \map_{[t_n, 0]}}_{\text{transport}}) \cdot \underbrace{\jacobianOperator(\map_{[t_n, 0]})}_{\text{compression}}.
\end{equation}
These two factors have distinct numerical properties. The transport factor involves only the evaluation of the initial condition, and is controlled by the error in the map values. The compression factor instead depends on the spatial derivatives of the map, and using the identity
\begin{equation}
    \label{eq:jacobian_product_rule}
    \jacobianOperator(\map_1 \circ \map_2) = \jacobianOperator(\map_1) \circ \map_2 \cdot \jacobianOperator(\map_2),
\end{equation}
we see that it evolves multiplicatively via
\begin{equation}
    \label{eq:jacobian_multiplicative_evolution}
    \jacobianOperator(\map_{[t_n, 0]}) = \prod_{i = 1}^n \jacobianOperator(\map_{[t_i, t_{i - 1}]}) \circ \map_{[t_n, t_i]}.
\end{equation}
Errors made in the Jacobians of successive submaps therefore enter multiplicatively along the numerical trajectories. This sensitivity is further exacerbated by the evolution equation \eqref{eq:main_system_on_group}. Even if the derivative errors of the map can be controlled, a high-order differential operator $F(\solution)$ introduces a stiff, nonlinear evolution for $\map$ that can impose severe time-step restrictions on explicit numerical schemes. 

\subsection{Semidirect Product Structure}

These difficulties suggest separate numerical treatments of the transport and compression components of the group action. In particular, we leverage the previous semi-Lagrangian CM methods to compute $\map$ and design a stable numerical integration by promoting the log-Jacobian of the map
\begin{equation*}
    \jacobian_{[t, t_0]} := \log \left( \jacobianOperator(\map_{[t, t_0]}) \right)
\end{equation*}
to be a primary variable for discretization. We then work in the extended variables
\begin{equation*}
    (\map, \jacobian, \mass)_{[t, t_0]} \in \semidirectGroup := \Diff \ltimes C^\infty (\manifold) \times \mathbb{R}_+.
\end{equation*}
For notational convenience, we will assemble the extended variables into the boldface notation
\begin{equation*}
    \groupElement_{[t, t_0]} := (\map, \jacobian, \mass)_{[t, t_0]}.
\end{equation*}
The extended group $\semidirectGroup$ is a semidirect product group. Its multiplication is given by
\begin{equation*}
    (\map_1, \jacobian_1, \mass_1) \cdot (\map_2, \jacobian_2, \mass_2) = (\map_1 \circ \map_2, \jacobian_1 \circ \map_2 + \jacobian_2, \mass_1 \mass_2),
\end{equation*}
which follows from the Jacobian product rule \eqref{eq:jacobian_product_rule}. Given a partition $0 = \tau_0 < \tau_1 < \cdots < \tau_m = T$ of the time interval $[0, T]$, we have the usual semigroup decomposition of a path $\groupElement_{[t, \tau_0]} \in \semidirectGroup$,
\begin{equation*}
    \groupElement_{[t, \tau_0]} = \groupElement_{[\tau_1, \tau_0]} \cdots \groupElement_{[t, \tau_k]},
\end{equation*}
whenever $t \in [\tau_k, \tau_{k + 1}]$. By defining a right action of $\semidirectGroup$ on $C^\infty (\manifold)$,
\begin{equation*}
    \rho \cdot (\map, \jacobian, \mass) := \mass (\rho \circ \map) e^\jacobian,
\end{equation*}
we obtain an expression for the solution in terms of the subinterval group elements
\begin{equation}
    \label{eq:semigroup_decomp_plus_advection}
    \rho(t) = \rho_0 \cdot \groupElement_{[\tau_1, \tau_0]} \cdots \groupElement_{[t, \tau_k]}.
\end{equation}
This group action allows us to explicitly separate the integration of the transport and compression components and forms the basis of the proposed numerical method. In the remainder of this work, we will work at a time $t \in [\tau_k, \tau_{k + 1}]$, evolving the ``current submap'' in this time interval. We will usually reserve the time-indexed notation $\groupElement_{[\tau_{i + 1}, \tau_i]}$ to mean a submap from an earlier time interval, while using the shorthand $\groupElement(t) = \groupElement_{[t, \tau_k]}$ for the current submap, where the initial time $\tau_k$ is understood.

\begin{remark}
    The choice to discretize the log-Jacobian $\jacobian$ over the Jacobian $\jacobianOperator(\map)$ avoids issues related to the loss of positivity of the latter during integration. The additive accumulation over many submaps via $\exp(\jacobian_1 + \cdots + \jacobian_m)$ was also observed to be marginally more numerically stable than the multiplicative counterpart $\jacobianOperator(\map_1) \cdots \jacobianOperator(\map_m)$.
\end{remark}

\subsection{Explicit Log-Jacobian Evolution}

In order to compute the evolution via \eqref{eq:semigroup_decomp_plus_advection}, we require an additional evolution equation for the submap log-Jacobian $\jacobian := \jacobian_{[t, t_0]}$. For a general velocity field $v$, $\jacobian$ satisfies
\begin{equation}
    \label{eq:general_jacobian_evolution_equation}
    \partial_t \jacobian + v \cdot \nabla \jacobian = -\operatorname{div} v, \quad\quad j(\tau_k) = 0.
\end{equation}
Indeed, the Jacobian $\jacobianOperator(\map_{[t, t_0]})$ by definition satisfies $\jacobianOperator(\map_{[t, t_0]}) \mu = \map_{[t, t_0]}^\ast \mu$. Taking the time derivative of this expression and inserting $\jacobianOperator(\map) = e^\jacobian$ gives \eqref{eq:general_jacobian_evolution_equation}.

In our case, the velocity field is determined in terms of $\rho$ by the first component of the lifting operator, $v = \liftingOperator_1(\rho, t)$. When $\rho$ is given in terms of the submaps \eqref{eq:semigroup_decomp_plus_advection}, we see that the resulting evolution \eqref{eq:general_jacobian_evolution_equation} for $\jacobian$ will depend on the current map $\map_{[t_,\tau_k]}$ and mass $\mass_{[t_,\tau_k]}$, as well as the previous submaps $\groupElement_{[\tau_{i + 1}, \tau_i]}$. It is convenient to introduce an additional quantity in order to isolate the contributions of these quantities in the evolution of $\jacobian$. For $k \geq 0$, we define the scalar functions
\begin{equation}
    \label{eq:j_source_term_def}
    \jsource_k : \manifold \times [\tau_k, \tau_{k + 1}] \rightarrow \mathbb{R}, \quad\quad \jsource_k (t) = \rho_0 \cdot \groupElement_{[\tau_1, \tau_0]} \cdots \groupElement_{[\tau_k, \tau_{k - 1}]} \cdot (\map_{[t, \tau_k]}, 0, \mass_{[t, \tau_k]}).
\end{equation}
$\jsource_k$ accounts for the transport and mass factor up to the current time, while omitting the compression associated with the current subinterval.  The full density field is recovered by restoring the current Jacobian factor via $\rho(t) = \jsource_k(t) e^{\jacobian_{[t, \tau_k]}}$, valid for $t \in [\tau_k, \tau_{k + 1}]$. This factorization isolates the explicit dependence on the current log-Jacobian and was found to be necessary for both the stability of the scheme and in the error analysis. From now on, we will write $\jsource := \jsource_k$ when the current time interval $[\tau_k, \tau_{k + 1}]$ is understood. With the introduction of $\jsource$, the evolution \eqref{eq:general_jacobian_evolution_equation} of $\jacobian$ becomes
\begin{equation}
    \label{eq:jacobian_evolution_equation}
    \partial_t \jacobian + v \cdot \nabla \jacobian = -\operatorname{div} v, \quad\quad v = \liftingOperator_1 (\jsource e^\jacobian, t).
\end{equation}
In particular, as $\jacobian$ does not enter $\jsource$, we have a separation between $\jacobian$ and the other submap quantities. The quantity $\jsource$ can be seen as a source or data term in the $\jacobian$-evolution.

For a specific PDE, expanding $v$ in \eqref{eq:jacobian_evolution_equation} often reveals a more structured equation for $\jacobian$, informing a well-suited choice of numerical integrator. We now consider the class of PDEs \eqref{eq:pde_with_highest_order_linear_term}, with a leading order linear term $\mathcal{L} = \alpha \Delta$ or $\mathcal{L} = -\alpha \Delta^2$ and a lower-order nonlinearity $G$.

\begin{proposition}
    \label{prop:jacobian_L_plus_N_form}
    In the case of \eqref{eq:pde_with_highest_order_linear_term}, the evolution equation for $\jacobian$ can be written in the form
    \begin{equation}
        \label{eq:jacobian_L_plus_N_form}
        \partial_t \jacobian = \mathcal{L} \jacobian + \mathcal{N}(j, \jsource, t)
    \end{equation}
    where $\mathcal{N}$ is a differential operator in $\jacobian$ of a lower order than $\mathcal{L}$. In the case of $\mathcal{L} = \alpha \Delta$, we obtain
    \begin{equation}
        \label{eq:jacobian_L_plus_N_form_order2}
        \partial_t \jacobian + (u - \alpha \nabla \log \jsource) \cdot \nabla \jacobian = \mathcal{L} \jacobian + \alpha |\nabla \jacobian|^2 + \mathcal{L} \log \jsource + \frac{1}{\jsource e^\jacobian} \left( G_0 - \nabla \log \jsource \cdot \nabla \Delta^{-1} G_0 \right),
    \end{equation}
    where
    \begin{equation}
        \label{eq:jacobian_L_plus_N_form_quantity_definitions}
        u = A(\jsource e^\jacobian - c, t), \quad\quad G_0 := G(\jsource e^\jacobian - c, t) - \xi \jsource e^\jacobian, \quad\quad \xi = \frac{\int_\manifold G(\jsource e^\jacobian - c, t) \mu}{\int_\manifold \jsource e^\jacobian \mu}.
    \end{equation}
\end{proposition}

We see that the evolution \eqref{eq:jacobian_L_plus_N_form} for $\jacobian$ inherits the stiffness of the original PDE, justifying the splitting of the transport and compression components. Conversely, we see that \eqref{eq:jacobian_L_plus_N_form} retains the leading-order form of the equation, which we exploit for a stable integration in the upcoming section. In the fourth-order case $\mathcal{L} = -\alpha \Delta^2$, we obtain an expression similar to \eqref{eq:jacobian_L_plus_N_form_order2}, only with more mixed lower-order terms. The structure in this case is not particularly revealing, and most importantly we still have the leading-order-linear form \eqref{eq:jacobian_L_plus_N_form}.

\begin{proof}
    Establishing \eqref{eq:jacobian_L_plus_N_form} is straightforward by plugging the expression \eqref{eq:lifting_operator_def} for the lift into the log-Jacobian evolution equation \eqref{eq:jacobian_evolution_equation}. We illustrate the second-order case \eqref{eq:jacobian_L_plus_N_form_order2}; the case $\mathcal{L} = -\alpha \Delta^2$ is similar, only with more mixed lower-order terms. 
    
    In the case $F = \alpha \Delta + G$, the lift $(v, \xi) = \liftingOperator(\rho, t)$, after inserting $\rho = \jsource e^\jacobian$, has the expression
    \begin{equation*}
        v = u + \frac{1}{\rho} \nabla \Delta^{-1} (\xi \rho - \alpha \Delta \rho - G(\rho - c, t)) = u - \alpha \nabla \jacobian - \alpha \nabla \log \jsource - \frac{1}{\jsource e^\jacobian} \nabla \Delta^{-1} G_0,
    \end{equation*}
    where we have used that $\nabla \Delta^{-1} \Delta \rho = \nabla \rho$ since $\nabla$ forgets the constant of integration.  Here $G_0$ is defined in \eqref{eq:jacobian_L_plus_N_form_quantity_definitions}, and $\xi$ is given in \eqref{eq:jacobian_L_plus_N_form_quantity_definitions} after using that $\Delta \rho$ has mean zero. We then compute
    \begin{equation*}
        \begin{aligned}
            \partial_t \jacobian = -\operatorname{div} v - v \cdot \nabla \jacobian = \alpha \Delta \jacobian + \alpha \Delta \log \jsource + \frac{1}{\jsource e^\jacobian} G_0 - \frac{1}{(\jsource e^\jacobian)^2} (e^\jacobian \nabla \jsource + \jsource e^\jacobian \nabla \jacobian) \cdot \nabla \Delta^{-1} G_0 \\
            - u \cdot \nabla \jacobian + \alpha |\nabla \jacobian|^2 + \alpha \nabla \log \jsource \cdot \nabla j + \frac{1}{\jsource e^\jacobian} \nabla \Delta^{-1} G_0 \cdot \nabla \jacobian,
        \end{aligned}
    \end{equation*}
    where we have used $\Delta \Delta^{-1} G_0 = G_0$ since $G_0$ has mean zero. This simplifies to \eqref{eq:jacobian_L_plus_N_form_order2}.
\end{proof}

\section{Numerical Implementation and Analysis}
\label{section:implementation_analysis}

In this section we describe an algorithm for the solution of problems \eqref{eq:pde_with_highest_order_linear_term} with a leading-order linear term. We consider the domain $\mathbb T^2 = \mathbb{R}^2 / 2 \pi \mathbb{Z}^2$, lending an implementation using Fourier and spline discretizations in space. The map and mass are advanced with RK4, and the log-Jacobian is evolved in a stable manner with an exponential integrator. We then prove an asymptotic error estimate, indicating an accuracy of $\mathcal{O}(\Delta t^4 + \Delta x^6 + \Delta t \Delta x^4)$.

We take uniform times $t_0, t_1, \ldots$ with step $\Delta t$, conforming to the remapping times $\tau_k$.

\subsection{Spatial Discretization}

We denote numerical maps, log-Jacobians, and masses by $\numericalMap$, $\numericalJacobian$, and $\numericalMass$, respectively. The numerical maps $\numericalMap$ and log-Jacobians $\numericalJacobian$ are discretized on a uniform periodic grid $\Gamma$ over $[0, 2\pi)^2$, with $N$ points in each dimension. As in the smooth case, we collect the map components into the boldface notation $\numericalGroupElement = (\numericalMap, \numericalJacobian, \numericalMass)$, reserving the time-indexed notation $\numericalGroupElement_{[\tau_{i + 1}, \tau_i]}$ for a completed submap and using the shorthand $\numericalGroupElement^n = \numericalGroupElement_{[t_n, \tau_k]}$ for the current submap.

We make extensive use of spatial interpolation. To this end, we employ uniform periodic spline interpolation of order 6 \cite{schoenberg_cardinal_1973}. The interpolation operator $\interp : C^0 (\mathbb{T}^2) \rightarrow C^0 (\mathbb{T}^2)$ projects to the unique interpolant which is a quintic polynomial in each grid cell and is globally $C^4$. On the torus, since the spline coefficients of the interpolant are related to the grid data by a circular convolution, the interpolant may be constructed by a simple division in Fourier space. The interpolation may then be evaluated in a numerically-stable manner via de Boor's algorithm \cite{lee_simplified_1982}. For $f \in C^6 (\mathbb{T}^2)$, we have the uniform interpolation error \cite{behforooz_overconvergence_1988}
\begin{equation*}
    \lVert \interp[f] - f \rVert_{C^0} = \mathcal{O}(\Delta x^6).
\end{equation*}
We also note that uniform spline interpolation on $\mathbb{T}^2$ is $C^0$-stable: $\lVert \interp \rVert_{C^0 \rightarrow C^0}$ is bounded uniformly in the grid size $N$ \cite{schurer_note_1968}. Finally, we note that the interpolation enjoys a useful superconvergence property of the derivative on the grid nodes: for $f \in C^6 (\mathbb{T}^2)$ we have the full order
\begin{equation*}
    \sup_\Gamma |\nabla \interp[f] - \nabla f| = \mathcal{O}(\Delta x^6),
\end{equation*}
rather than the usual derivative order reduction $\mathcal{O}(\Delta x^5)$ \cite{behforooz_overconvergence_1988}.

We extend the interpolation operator to vector fields and maps in the natural way, overloading the operator $\interp$. For vector fields $v$, the interpolation is defined componentwise, $\interp[v] = (\interp[v_1], \interp[v_2])$. For a map $\numericalMap$ we interpolate the deformation of the identity, $\interp[\numericalMap] := \operatorname{id}{} + \interp[\numericalMap - \operatorname{id}]$. We finally define the componentwise interpolation $\grpInterp[\numericalGroupElement] = (\interp[\numericalMap], \interp[\numericalJacobian], \numericalMass)$.

\subsection{Numerical Integrators for the Submap Components}
\label{section:numerical_propagators}

We now define three numerical integrators $\intMap$, $\intJac$ and $\intMass$ for the map, log-Jacobian and mass, respectively. These propagators take as input numerical approximations $\tilde{v}$, $\tilde{\jsource}$ and $\tilde{\xi}$ of the relevant velocity quantity over the interval $[t_n, t_{n + 1}]$. For example, the map propagator defines an update of the form
\begin{equation*}
    \numericalMap^{n + 1} = \intMap^n (\numericalMap^n, \tilde{v}).
\end{equation*}
The superscript $n$ denotes the dependence of the propagator on $t_n$ and $\Delta t$.

\subsubsection{Map and Mass Integration}

For the Lagrangian quantities, we employ the RK4 integrator. Specifically, let $\Phi_{\Delta t} (x_0, v, t_0)$ denote one RK4 step approximating $x(t_0 + \Delta t)$, where $\dot{x} = v(x, t)$ and $x(t_0) = x_0$. Then the map integrator is defined by the usual semi-Lagrangian step, which consists of integrating the time-$\Delta t$ departure and evaluating the previous map at the departure points via interpolation:
\begin{equation*}
    \intMap^n (\numericalMap, \tilde{v}) := \interp \left[ \numericalMap \right] \circ \numericalMap_{[t_{n + 1}, t_n]}, \quad\quad \numericalMap_{[t_{n + 1}, t_n]} := \Phi_{-\Delta t} (\operatorname{id}{}, \tilde{v}, t_{n + 1}).
\end{equation*}
In terms of the usual semi-Lagrangian propagator $\slPropagator$ on scalar functions
\begin{equation*}
    \slPropagator^n(f) := \interp[f] \circ X_{[t_{n + 1}, t_n]},
\end{equation*}
the map integrator can be written
\begin{equation*}
    \intMap^n (X) = X_{[t_{n + 1}, t_n]} + \slPropagator^n(X - \operatorname{id}).
\end{equation*}
The integrator for the mass is simply taken to be one RK4 step,
\begin{equation*}
    \intMass^n(\numericalMass, \tilde{\xi}) := \Phi_{\Delta t} (\numericalMass, m \mapsto \tilde{\xi} m, t_n).
\end{equation*}

\subsubsection{Log-Jacobian Integration}
\label{section:jacobian_time_stepper}

From Proposition \ref{prop:jacobian_L_plus_N_form}, we have that $\jacobian$ satisfies an evolution of the form
\begin{equation*}
    \partial_t \jacobian = \mathcal{L} \jacobian + \mathcal{N}(\jacobian, \jsource, t),
\end{equation*}
where $\mathcal{N}$ is lower-order than $\mathcal{L}$ in $\jacobian$. There exist a number of integrators which are well-suited to the stable integration of problems of this form. Among them is the class of exponential integrators, which we employ here. For a background, we refer the survey \cite{Hochbruck_Ostermann_2010}, as well as to a numerical application of such methods to the Allen--Cahn and Cahn--Hilliard equations \cite{montanelli_solving_2020}. With Duhamel's formula applied to the PDE \eqref{eq:jacobian_L_plus_N_form} for $\jacobian$, we can write the exact solution over one time step as
\begin{equation*}
    \jacobian(t_{n + 1}) = e^{\Delta t \mathcal{L}} \jacobian(t_n) + \int_{t_n}^{t_{n + 1}} e^{(t_{n + 1} - s) \mathcal{L}} \mathcal{N}(\jacobian(s), \jsource(s), s) d s.
\end{equation*}
Based on this expression, an exponential integrator approximates $\jacobian(t_{n + 1})$ by exponentiating the linear term $\mathcal{L}$ exactly and approximating the integral above by quadrature. As a result, an exponential integrator bypasses the CFL condition associated with the operator $\mathcal{L}$.

In this work we employ the fourth-order, five-stage exponential time differencing (ETD) method of Hochbruck and Ostermann \cite{Hochbruck_Ostermann_2005}. We define the update map $\intJac$ by
\begin{equation*}
    \intJac^n (\numericalJacobian, \tilde{\jsource}) := e^{\Delta t \mathcal{L}} \numericalJacobian + \Delta t \sum_{i = 1}^5 b_i \mathcal{N}_i,
\end{equation*}
where the internal stages
\begin{equation*}
    \begin{aligned}
        k_i &= e^{c_i \Delta t \mathcal{L}} \numericalJacobian + \Delta t \sum_{\ell = 1}^{i - 1} a_{i \ell} \mathcal{N}_\ell, \\
        \mathcal{N}_i &= \mathcal{N}(k_i, \tilde{\jsource}(t_n + c_i \Delta t), t_n + c_i \Delta t).
    \end{aligned}
\end{equation*}
The coefficients $a_{i \ell}$ and $b_i$ are operators, given in terms of the exponential $\varphi$ functions; these are defined recursively by
\begin{equation}
    \label{eq:exponential_phi_functions_def}
    \varphi_0 (z) = e^z, \quad\quad \varphi_{k + 1} (z) = \frac{\varphi_k (z) - \frac{1}{k!}}{z}.
\end{equation}
One then has, for instance, $a_{11} = \frac{1}{2} \varphi_1 \left( \frac{\Delta t}{2} \mathcal{L} \right)$. The full Butcher tableau for the 4th-order, 5-stage ETD method of Hochbruck and Ostermann can be found in \cite[Eq. (5.19)]{Hochbruck_Ostermann_2005}. This ETD method was chosen as it satisfies the stiff order conditions in the weak form to fourth order \cite{Hochbruck_Ostermann_2005}, which in our strongly parabolic regime should be sufficient to ensure convergence.

In our case, the action of the exponential $\varphi$ functions are given by diagonal operators in Fourier space. The expressions \eqref{eq:exponential_phi_functions_def} suffer cancellation errors for small $z$ and cannot be evaluated directly. For our purposes, we precompute $\varphi_k (z)$ near $z = 0$ by a high-degree Taylor expansion in quadruple precision. While we have found this method to be sufficient in our case of a diagonal $\mathcal{L}$, we note that there exist specialized methods for the stable evaluation of the action of $\varphi$ functions in a more general context, see \cite[§4]{Hochbruck_Ostermann_2010}.

\subsection{Algorithm Implementation}

One time step of the proposed algorithm is presented at a high level in Algorithm \ref{alg:time_step}.

\begin{algorithm}[h!]
\caption{One time step of the numerical implementation.}
\label{alg:time_step}

\newcommand{\algstage}[1]{\makebox[4.5em][l]{\textbf{#1}}}

\begin{enumerate}[
    label=\textbf{\arabic*.}, ref=\arabic*,
    leftmargin=*, labelsep=0.6em,
    itemsep=0.7em, topsep=0.4em, parsep=0pt
]
    \item \label{item:alg_step_1} \algstage{Sample}%
    \(\begin{aligned}[t]
        \rho^n &\gets
        \rho_0 \cdot \numericalGroupElement_{[\tau_1,\tau_0]}
        \cdots \numericalGroupElement_{[\tau_k,\tau_{k-1}]}
        \cdot \numericalGroupElement^n, \\
        \jsource^n &\gets
        \rho_0 \cdot \numericalGroupElement_{[\tau_1,\tau_0]}
        \cdots \numericalGroupElement_{[\tau_k,\tau_{k-1}]}
        \cdot (\numericalMap^n,0,\numericalMass^n).
    \end{aligned}\)

    \item \label{item:alg_step_2} \algstage{Lift}%
    \((v^n,\xi^n)\gets\liftingOperator(\rho^n,t_n).\)

    \item \label{item:alg_step_3} \algstage{Flow}%
    \(\begin{aligned}[t]
        &\text{Interpolate } (v^i, q^i, \xi^i) \mapsto (\tilde{v}^n,\tilde{\jsource}^n,
        \tilde{\xi}^n)\text{, then advance:}\\[0.3ex]
        &\begin{aligned}[t]
            \numericalMap^{n+1} \gets \intMap^n(\numericalMap^n,\tilde{v}^n), \quad\quad \numericalJacobian^{n + 1} \gets \intJac^n(\numericalJacobian^n,\tilde{\jsource}^n), \quad\quad \numericalMass^{n + 1} \gets \intMass^n(\numericalMass^n,\tilde{\xi}^n).
        \end{aligned}
    \end{aligned}\)

    \item \label{item:alg_step_4} \algstage{Remap}%
    \(\begin{aligned}[t]
        &\text{If }t_{n+1}=\tau_{k+1}\text{, store and reset:}\\[0.3ex]
        &\numericalGroupElement_{[\tau_{k+1},\tau_k]}
        \gets\numericalGroupElement^{n+1},
        \qquad
        \numericalGroupElement^{n+1}\gets(\operatorname{id},0,1).
    \end{aligned}\)
\end{enumerate}
\end{algorithm}

Each time step, we begin in step \ref{item:alg_step_1} by approximating the state $\rho(t_n)$. Numerically, this is performed by interpolating all of the submaps, sampling the group action on the grid $\Gamma$, and truncating the result to a number $N_\epsilon$ of modes:
\begin{equation}
    \label{eq:computation_of_rho_n_with_truncation}
    \rho^n := T_{N_\epsilon} \left[ \rho_0 \cdot \grpInterp [\numericalGroupElement_{[\tau_1, \tau_0]}] \cdot \ldots \cdot \grpInterp [\numericalGroupElement_{[\tau_k, \tau_{k - 1}]}] \cdot \grpInterp[\numericalGroupElement^n] \right].
\end{equation}
Here we have written $T_{N_\epsilon}$ for the Fourier truncation operator to $N_\epsilon$ modes. This truncation serves a number of purposes. First, it implies that we are effectively solving for a Galerkin truncation of the original PDE, which facilitates the error analysis, see Section \ref{section:error_estimate}. Second, as was observed in the CM method for the incompressible Euler equations \cite{yin_characteristic_2021}, this truncation is necessary to prevent the exponentially-fine details generated by the composition of maps from polluting the velocity, and subsequently reducing the quality of the approximation space. We sample the log-Jacobian data term $\jsource$ similarly, though this time without the truncation:
\begin{equation*}
    \jsource_k^n := \rho_0 \cdot \grpInterp [\numericalGroupElement_{[\tau_1, \tau_0]}] \cdot \ldots \cdot \grpInterp[\numericalGroupElement_{[\tau_k, \tau_{k + 1}]}] \cdot (\interp[\numericalMap^n], 0, \numericalMass^n).
\end{equation*}
For $\jsource$, the truncation was not found to be necessary, neither for stability nor for the error analysis. We note that $\rho^n$ as computed in \eqref{eq:computation_of_rho_n_with_truncation} does not represent the solution of the method. It is an approximation introduced only to advance the time integration. The solution of the method is defined to be the full element of the growing approximation space as
    \begin{equation*}
        \rho(t_n) \approx \rho_0 \cdot \grpInterp[\numericalGroupElement_{[\tau_1, \tau_0]}] \cdot \ldots  \cdot \grpInterp[\numericalGroupElement_{[t_n, \tau_k]}].
    \end{equation*}

In step \ref{item:alg_step_2} of the algorithm we compute the lift. When approximating the lifting operator \eqref{eq:lifting_operator_def}, differential operators are evaluated in Fourier space, and nonlinear operations are dealiased using the two-thirds rule.

Step \ref{item:alg_step_3} propagates the submap quantities over one time step, using the integrators defined in Section \ref{section:numerical_propagators}. These integrators require high-order approximations of their respective velocities $v$, $\jsource$ and $\xi$. To this end, we interpolate the velocity data $v^i$, $\jsource^i$, $\xi^i$ in time over the last few steps using 4th-order Lagrange interpolation. We also interpolate the map velocities $v^i$ in space, as the semi-Lagrangian integration of the map requires evaluation of the velocity at arbitrary locations in the domain. These interpolants are defined as
\begin{equation*}
    \tilde{v}^n (t) := \sum_{i = n - 3}^n \interp[v^i] \ell_{t_i} (t), \quad\quad \tilde{\jsource}^n (t) := \sum_{i = n - 3}^n \jsource^i \ell_{t_i} (t), \quad\quad \tilde{\xi}^n (t) := \sum_{i = n - 3}^n \xi^i \ell_{t_i} (t).
\end{equation*}
Here $\ell_{t_i}$ denotes the cubic Lagrange basis polynomial $\ell_{t_i} (t_j) = \delta_{ij}$ over the nodes $(t_{n - 3}, \ldots, t_n)$. 

\begin{remark}
    For simplicity, we sample the velocity data $(\rho^n, \jsource^n, v^n)$ on the same grid $\Gamma$ as the map $\numericalMap$ and log-Jacobian $\numericalJacobian$. We remark that the sampling space of the velocity data is in principle independent of that of the map, and the velocities could be represented on a finer grid. Such a scale separation has been exploited in previous CM methods \cite{yin_characteristic_2021} to ensure a faithful representation of fine-scale velocities.
\end{remark}

\begin{remark}
    The numerical velocities $(\tilde{v}^n, \tilde{\jsource}^n, \tilde{\xi}^n)$ are constructed in terms of past samples. To initialize the algorithm, for the first few steps, we replace these with high-order approximations of $(v, \jsource, \xi)$, which can be constructed as follows. Using the initial data and the PDE, we can construct a Taylor series in time of the solution $\solution$ about $t = 0$, which gives a high-order approximation to $(v, \xi) = \liftingOperator(\solution + c, t)$. From these quantities we can approximate the flows $(\map, \mass)_{[t, 0]}$, and hence  $\jsource (t) = \rho_0 \cdot (\map, 0, \mass)_{[t, 0]}$, to a high order for small times. 
\end{remark}

\subsection{Error Estimate}
\label{section:error_estimate}

We now give an error estimate for the numerical solution $\rho^n$ computed using one submap. The estimate can be derived under a number of technical assumptions, which we list below. These assumptions pertain to the stability of the scheme and to properties of the chosen numerical integrators. In particular, we assume without justification that the scheme is stable, possibly under a time step restriction. Given the wide range of equations we consider, establishing conditional stability is outside the scope of the current work. However, we will see in the upcoming section that the scheme is energy-stable under a reasonable time step restriction of the form $\Delta t \lesssim \Delta x$. In the following, we denote by $p \in \{2, 4\}$ the order of the operator $\mathcal{L}$.

\paragraph{Error analysis assumptions}
\begin{enumerate}[noitemsep, topsep=0em, label=(\roman*)]
    \item General assumptions. We assume that $A$ and $G$ are locally Lipschitz maps $C^{p - 1} \rightarrow C^0$. We assume that the solution $\groupElement(t)$ and that all the numerical solutions $\numericalGroupElement^n$ remain in a ball in $C^p$, uniformly in $t$ and $n$, for all sufficiently small $\Delta t$.
    
    \item\label{item:error_assumption_trunc} Truncated solution. We compare against a Galerkin truncation of the PDE to a fixed number $N_\epsilon$ of modes, and assume that the truncation satisfies $\lVert T_{N_\epsilon} [\rho] - \rho \rVert_{C^p} \leq \epsilon$ for some $\epsilon$.
    
    \item\label{item:error_assumption_sl} Stability of semi-Lagrangian propagator. We assume that the semi-Lagrangian propagator $\slPropagator$ is $C^0$-stable, meaning
    \begin{equation}
        \label{eq:sl_propagator_c0_stable}
        \lVert \slPropagator^{n - 1} \circ \cdots \circ \slPropagator^{i + 1} \rVert_{C^0 \rightarrow C^0} = \mathcal{O}(1)
    \end{equation}
    uniformly in $i$, $n$, $\Delta t$. This property holds when replacing $\slPropagator$ with the exact transport operator, and is automatic in the weaker $L^2$ norm \cite{ferretti_stability_2019}. We further assume that the displacement of $\slPropagator$ can be measured in a weaker norm: for all $f \in C^0$,
    \begin{equation}
        \label{eq:sl_propagator_displacement_bound}
        \lVert (\slPropagator - I) f \rVert_{C^{-1}} \lesssim \Delta t \lVert f \rVert_{C^0}.
    \end{equation}
    Here, $\lVert \cdot \rVert_{C^{-1}}$ denotes the inhomogeneous Besov norm $B_{\infty, \infty}^{-1}$ \cite[§2.7]{bahouri_fourier_2011}. This property also mirrors the exact transport operator: if $\partial_t f + v \cdot \nabla f = 0$ then
    \begin{equation*}
        \lVert f(t + \Delta t) - f(t) \rVert_{C^{-1}} \lesssim \Delta t \lVert \partial_t f \rVert_{C^{-1}} \lesssim \Delta t \lVert v \cdot \nabla f \rVert_{C^{-1}} \lesssim \Delta t \lVert f \rVert_{C^0}.
    \end{equation*}
    Of course, it suffices to assume \eqref{eq:sl_propagator_c0_stable} and \eqref{eq:sl_propagator_displacement_bound} only on the set of numerical solutions.

    \item\label{item:error_assumption_etd} Error of ETD integrator. The ETD step $\numericalJacobian^{n + 1} = \intJac^n (\numericalJacobian^n, \tilde{\jsource}^n)$ approximates $\tilde{\jmath}^n (t_{n + 1})$, where $\tilde{\jmath}^n$ is the solution of the following problem in $[t_n, t_{n + 1}]$:
    \begin{equation}
        \label{eq:jacobian_modified_pde}
        \partial_t \tilde{\jmath}^n = \mathcal{L} \tilde{\jmath}^n + \mathcal{N}(\tilde{\jmath}^n, \tilde{\jsource}^n, t), \quad\quad \tilde{\jmath}^n (t_n) = \tilde{\jmath}^{n - 1} (t_n),
    \end{equation}
    with $\tilde{\jmath}^0 (0) = 0$. We assume that the ETD integrator is stable and satisfies a global error
    \begin{equation*}
        \lVert \numericalJacobian^n - \tilde{\jmath}^n (t_n) \rVert_{C^0} = \mathcal{O}(\Delta t^4)
    \end{equation*}
    uniformly in $n$ and $\Delta t$, and across all source terms $\tilde{q}^n$, possibly under a time step restriction. This assumption is motivated by the convergence theory of Hochbruck and Ostermann, see \cite[Thm. 2.17]{Hochbruck_Ostermann_2010}, also \cite{Hochbruck_Ostermann_2005}. We note that, for higher-order integrators, the analysis in \cite{Hochbruck_Ostermann_2010, Hochbruck_Ostermann_2005} does not let $\mathcal{N}$ be an unbounded operator. However, there is ample numerical evidence supporting the convergence of higher-order ETD methods in this case despite a missing formal convergence analysis, see e.g. \cite{montanelli_solving_2020} for the Allen--Cahn and Cahn--Hilliard equations.
\end{enumerate}

\begin{theorem}
    \label{thm:error_estimate}
    Under the assumptions above, the global error satisfies the asymptotic estimate
    \begin{equation*}
        \lVert \rho_0 \cdot \numericalGroupElement^n - \rho(t_n) \rVert_{C^0} = \mathcal{O}(\Delta t^4 + \Delta x^6 + \Delta t \Delta x^4 + \epsilon).
    \end{equation*}
\end{theorem}

Theorem \ref{thm:error_estimate} is established in \ref{appendix:proof_of_error_theorem}. In the rest of this section, we give a high-level overview of the proof.

For clarity, we denote the target error as $\delta := \mathcal{O}(\Delta t^4 + \Delta x^6 + \Delta t \Delta x^4 + \epsilon)$.
Recall the log-Jacobian data term $\jsource = \mass \cdot (\rho_0 \circ \map)$. This transported and scaled quantity satisfies $\partial_t \jsource + v \cdot \nabla \jsource = \xi \jsource$, where $(v, \xi) = \liftingOperator(\rho, t)$. At each time step we have the numerical approximations $\jsource^n$, as well as their local time interpolations $\tilde{\jsource}^n$ over $[t_n, t_{n + 1}]$. We note that in general these quantities jump at the nodes, $\tilde{\jsource}^n (t_{n + 1}) \neq \tilde{\jsource}^{n + 1}(t_{n + 1})$. Recall finally from assumption \ref{item:error_assumption_etd} that the log-Jacobian integrator approximates the quantity $\tilde{\jmath}^n$, which satisfies the PDE $\partial_t \tilde{\jmath}^n = \mathcal{L} \tilde{\jmath}^n + \mathcal{N}(\tilde{\jmath}^n, \tilde{\jsource}^n, t)$.

For the purposes of the proof, we consider the quantities $\tilde{\rho}^n := \tilde{\jsource}^n e^{\tilde{\jmath}^n}$ in $[t_n, t_{n + 1}]$. We define their errors relative to $\rho$,
\begin{equation*}
    E^n := \sup_{s \in [t_n, t_{n + 1}]} \lVert \tilde{\rho}^n (s) - \rho(s) \rVert_{C^0}.
\end{equation*}
We have that in the time interval $[t_n, t_{n + 1}]$, $\tilde{\rho}^n$ solves the PDE
\begin{equation}
    \label{eq:rho_tilde_pde_section4}
    \partial_t \tilde{\rho}^n + \tilde{u} \cdot \nabla \tilde{\rho}^n = \mathcal{L} \tilde{\rho}^n + G(\tilde{\rho}^n - c, t) + e^{\tilde{\jmath}^n} \left( \partial_t \tilde{\jsource}^n + \tilde{v} \cdot \nabla \tilde{\jsource}^n - \tilde{\xi} \tilde{\jsource}^n \right)
\end{equation}
where we have defined $\tilde{u} := A(\tilde{\rho}^n - c, t)$ and $(\tilde{v}, \tilde{\xi}) := \liftingOperator(\tilde{\rho}^n, t)$. In other words, we have that $\tilde{\rho}^n$ solves the original PDE \eqref{eq:pde_with_highest_order_linear_term}, up to an error in the action of $\jsource$. Thus, we proceed by first bounding the error in the action of $\jsource$, then bounding the error in $\tilde{\rho}^n$.

We first establish a bound in the pointwise error in $\jsource$,
\begin{equation*}
    \sup_{s \in [t_n, t_{n + 1}]} \lVert \tilde{\jsource}^n (s) - \jsource(s) \rVert_{C^0} \lesssim E^n + \delta,
\end{equation*}
which is the content of Lemma \ref{lemma:jsource_lemma} in \ref{appendix:proof_of_error_theorem}. This is achieved by writing a recursion for the map and mass errors. The recursion reduces to a standard semi-Lagrangian estimate, and uses the $C^0$-stability assumption \eqref{eq:sl_propagator_c0_stable} of the semi-Lagrangian propagator. We next bound the error in the time derivative of $\jsource$,
\begin{equation*}
    \sup_{s \in [t_n, t_{n + 1}]} \lVert \tilde{\jsource}^n (s) - \jsource(s) \rVert_{C^0} \lesssim E^n + \delta.
\end{equation*}
This is the content of Lemma \ref{lemma:jsource_ddt_lemma}. The argument for the local orders is the same as in the pointwise case. To avoid losing a factor of $\Delta t$ in the time derivative, we can measure in the weaker spatial norm $C^{-1}$, using the assumption \eqref{eq:sl_propagator_displacement_bound}.

With the error in the action of $\jsource$ bounded, it remains to bound the error in $\tilde{\rho}^n$ using its PDE \eqref{eq:rho_tilde_pde_section4}. This is the content of Lemma \ref{lemma:error_resursion}. After using Duhamel's formula and exploiting the smoothing of the dissipative operator $\mathcal{L}$, this reduces to bounding the error in the quantity $\partial_t \tilde{\jsource}^n + \tilde{v} \cdot \nabla \tilde{\jsource}^n - \tilde{\xi} \tilde{\jsource}^n$ in a weaker spatial norm, which is achieved by the two preceding lemmas. We obtain an error recursion which is a singular discrete Grönwall inequality, and gives $E^n \lesssim \delta$. Finally, we have
\begin{equation*}
    \lVert \rho_0 \cdot \numericalGroupElement^n - \rho(t_n) \rVert_{C^0} \lesssim \lVert \numericalJacobian^n - \tilde{\jmath}^n (t_n) \rVert_{C^0} + \lVert \tilde{\rho}^n (t_n) - \rho(t_n) \rVert_{C^0}.
\end{equation*}
The first term has size $\delta$ by assumption \ref{item:error_assumption_etd} on the log-Jacobian integrator, and the second is at most $E^n \lesssim \delta$, which completes the proof.

\section{Numerical Validation and Examples}
\label{section:examples}

In this section, we validate the convergence of the proposed algorithm and perform numerical experiments illustrating its properties for a canonical set of dissipative equations. These experiments serve to validate the approach, and to illustrate the roles of the various components of the map in the group action. An in-depth study, including more diverse physical experiments and an analysis of the growing approximation space, especially for the Navier--Stokes equations, will be the subject of future work.

The experiments are performed using a GPU-enabled implementation of the algorithm and are run on an NVIDIA RTX 4070 with 12GB VRAM. The use of a GPU facilitates in particular the parallel evaluation of the spline interpolations, which is implemented via custom CUDA kernels.

\subsection{Convergence Tests}
\label{section:convergence_tests}

We first validate the error estimates of Theorem \ref{thm:error_estimate} with numerical convergence tests in space and time. We consider different dissipative evolution equations of the form $\solution_t + A(\solution) \cdot \nabla \solution = F(\solution)$ as defined in the Table \ref{fig:convergence_pdes_table}. The convergence of the proposed method will be measured using the method of manufactured solutions. Define the function
\begin{equation*}
   h(x, t) = \sin(2 y_1) \cos(3 y_2) \quad\text{where}\quad \begin{pmatrix} y_1 \\ y_2 \end{pmatrix} = \begin{pmatrix} \cos{t} & -\sin{t} \\ \sin{t} & \cos{t} \end{pmatrix} \begin{pmatrix} \sin{x_1} \\ \sin{x_2} \end{pmatrix}.
\end{equation*}
We consider the manufactured solution $\phi = h^2$ for the Allen--Cahn equation to test the case of varying total mass, and $\phi = h$ for all the other equations. Defining the function
\begin{equation*}
    g = \partial_t \solution + A(\solution) \cdot \nabla \solution - F(\solution),
\end{equation*}
it follows that $\solution$ is an exact solution of the PDE with source term
\begin{equation}
    \label{eq:manufactured_pde}
    \partial_t \solution + A(\solution) \cdot \nabla \solution = F(\solution) + g.
\end{equation}
The equation \eqref{eq:manufactured_pde} lifts to an evolution on the diffeomorphism group, in the coordinate $\rho = \solution + c$, via the lifting operator \eqref{eq:lifting_operator_def} as $(v, \xi) = \liftingOperator(\rho, t)$, and we find
\begin{equation*}
    v = A(\rho - c) + \frac{1}{\rho} \nabla \Delta^{-1} (\xi \rho - F(\rho - c) - g(x, t)).
\end{equation*}
In particular, we see that source terms require no additional treatment in our formulation. This is in contrast to traditional transport methods, where adding a source term requires an additional accumulation through a Duhamel integral \cite{yin_characteristic_2026}.

\begin{table}[h!]
\centering
\renewcommand{\arraystretch}{1.4}
\begin{tabular}{c|c|c|c|c}
\hline
\rowcolor{gray!15}
Equation & $\dot{m} = 0$ & $F(\solution)$ & $A(\solution)$ & Parameters \\ \hline
\makecell{Heat} & Yes & $\nu \Delta \solution$ & 0 & $\nu = 10^{-2}$ \\ \hline
\makecell{Allen--Cahn} & No & $\epsilon^2 \Delta \solution - \solution^3 + \solution$ & 0 & $\epsilon = 10^{-2}$ \\ \hline
\makecell{Cahn--Hilliard} & Yes & $D \Delta (\solution^3 - \solution - \epsilon^2 \Delta \solution)$ & 0 & \makecell{$D = 10^{-2}$, $\epsilon = 10^{-1}$} \\ \hline
\makecell{Navier--Stokes} & Yes & $\nu \Delta \solution$ & $-\nabla^\perp \Delta^{-1} \solution$ & $\nu = 10^{-2}$ \\ \hline
\end{tabular}
\caption{Equations and their parameters considered for the manufactured solution convergence tests.}
\label{fig:convergence_pdes_table}
\end{table}

For each convergence test, we run the method to a final time $T = 1$, with $c = 100$, using a single submap. The uniform $C^0$ error is estimated against the true solution $\solution(\cdot, 1)$ by sampling the numerical solution at the final time on a dense grid of size $2048^2$.
We take $\Delta t \in \{ 1/16, 1/32, \ldots, 1/256\}$ for the convergence test in time, requiring also a refinement $\Delta x \sim \Delta t^{1/2}$ in order to beat the CFL condition. The convergence tests in space take $\Delta x \in \{2\pi/16, 2\pi/32, \ldots, 2\pi/256\}$, this time with a fixed $\Delta t = 1/4096$. The errors are reported in Table \ref{fig:convergence_test}, along with the incremental rate, computed in time for instance as the log-difference
\begin{equation*}
    \text{rate}_i = \frac{\log(\text{error}_i) - \log(\text{error}_{i - 1})}{\log(\Delta t_i) - \log(\Delta t_{i - 1})}.
\end{equation*}
We observe the anticipated convergence at $\mathcal{O}(\Delta t^4)$ and $\mathcal{O}(\Delta x^6)$ as determined by the dominant terms in Theorem \ref{thm:error_estimate}.

\begin{table}[h!]
\centering
\begin{tabular}{c|c|c|cc|cc|cc|cc|}
\cline{2-11}
& \multicolumn{1}{c|}{\multirow{2}{*}{$\Delta t$}} & \multicolumn{1}{c|}{\multirow{2}{*}{$\Delta x$}} & \multicolumn{2}{c|}{\makecell{\rule{0pt}{2.5ex}Heat equation}} & \multicolumn{2}{c|}{\makecell{\rule{0pt}{2.5ex}Allen--Cahn}}   & \multicolumn{2}{c|}{\makecell{\rule{0pt}{2.5ex}Cahn--Hilliard}} & \multicolumn{2}{c|}{\makecell{\rule{0pt}{2.5ex}Navier--Stokes}} \\
& & & \multicolumn{1}{c:}{Error}  & Rate & \multicolumn{1}{c:}{Error} & Rate & \multicolumn{1}{c:}{Error}  & Rate & \multicolumn{1}{c:}{Error}  & Rate \\
\cline{2-11}
\multicolumn{1}{c|}{\multirow{5}{*}{\rotatebox[origin=c]{90}{Time}}} & $1/16$ & $2\pi/64$ & \multicolumn{1}{c:}{1.26e-05}    & --  & \multicolumn{1}{c:}{1.19e-03}   & --  & \multicolumn{1}{c:}{8.83e-05}    & --  & \multicolumn{1}{c:}{2.34e-04}    & --  \\
\multicolumn{1}{c|}{} & $1/32$ & $2\pi/96$ & \multicolumn{1}{c:}{8.62e-07}    & 3.87  & \multicolumn{1}{c:}{9.48e-05}   & 3.65  & \multicolumn{1}{c:}{3.74e-06}    & 4.56  & \multicolumn{1}{c:}{1.34e-05}    & 4.13  \\
\multicolumn{1}{c|}{} & $1/64$ & $2\pi/128$ & \multicolumn{1}{c:}{5.98e-08}    & 3.85  & \multicolumn{1}{c:}{6.39e-06}   & 3.89  & \multicolumn{1}{c:}{2.19e-07}    & 4.09  & \multicolumn{1}{c:}{8.05e-07}    & 4.05  \\
\multicolumn{1}{c|}{} & $1/128$ & $2\pi/192$ & \multicolumn{1}{c:}{3.92e-09}    & 3.93  & \multicolumn{1}{c:}{4.31e-07}   & 3.89  & \multicolumn{1}{c:}{1.61e-08}    & 3.76  & \multicolumn{1}{c:}{5.06e-08}    & 3.99  \\
\multicolumn{1}{c|}{} & $1/256$ & $2\pi/256$ & \multicolumn{1}{c:}{2.50e-10}    & 3.97  & \multicolumn{1}{c:}{2.79e-08}   & 3.95  & \multicolumn{1}{c:}{1.10e-09}    & 3.87  & \multicolumn{1}{c:}{3.18e-09}    & 3.99  \\
\doublehhline{~|=|=|=|=|=|=|=|=|=|=|}
\multicolumn{1}{c|}{\multirow{5}{*}{\rotatebox[origin=c]{90}{Space}}} & $1/4096$ & $2\pi/16$ & \multicolumn{1}{c:}{7.44e-03}    & --  & \multicolumn{1}{c:}{9.20e-02}   & --  & \multicolumn{1}{c:}{7.46e-03}    & --  & \multicolumn{1}{c:}{6.06e-02}    & --  \\
\multicolumn{1}{c|}{} & $1/4096$ & $2\pi/32$ & \multicolumn{1}{c:}{4.39e-05}    & 7.41  & \multicolumn{1}{c:}{4.89e-04}   & 7.55  & \multicolumn{1}{c:}{4.39e-05}    & 7.41  & \multicolumn{1}{c:}{6.03e-04}    & 6.65  \\
\multicolumn{1}{c|}{} & $1/4096$ & $2\pi/64$ & \multicolumn{1}{c:}{4.37e-07}    & 6.65  & \multicolumn{1}{c:}{3.50e-06}   & 7.13  & \multicolumn{1}{c:}{4.37e-07}    & 6.65  & \multicolumn{1}{c:}{2.91e-06}    & 7.70  \\
\multicolumn{1}{c|}{} & $1/4096$ & $2\pi/128$ & \multicolumn{1}{c:}{5.99e-09}    & 6.19  & \multicolumn{1}{c:}{4.25e-08}   & 6.36  & \multicolumn{1}{c:}{5.99e-09}    & 6.19  & \multicolumn{1}{c:}{2.74e-08}    & 6.73  \\
\multicolumn{1}{c|}{} & $1/4096$ & $2\pi/256$ & \multicolumn{1}{c:}{9.06e-11}    & 6.05  & \multicolumn{1}{c:}{6.28e-10}   & 6.08  & \multicolumn{1}{c:}{9.04e-11}    & 6.05  & \multicolumn{1}{c:}{3.86e-10}    & 6.15  \\ \cline{2-11}
\end{tabular}
\caption{Manufactured solution convergence of the method for all four PDEs in time and space.}
\label{fig:convergence_test}
\end{table}

\subsection{Heat Equation}

As a first example, we consider the purely dissipative dynamics of the heat equation $\rho_t = \Delta \rho$. The heat equation admits a formulation as an $H^{-1}$ gradient flow on $\Prob$ with constant mass $\dot{\mass} = 0$, and the lifting operator \eqref{eq:lifting_operator_def} defines a gradient vector field
\begin{equation}
    \label{eq:lift_heat_equation}
    v = \liftingOperator_1 (\rho) = -\nabla \log \rho.
\end{equation}
In particular, we find that our lift for the heat equation coincides with the Otto lift (c.f. \eqref{eq:otto_lift}).

\begin{remark}
    More generally, our lift (in the balanced case $\mass \equiv 1$) coincides with the lift defined by the Otto metric for any equation of the form $\partial_t \rho = \Delta f(\rho)$ where $f \in C^1(\mathbb{R})$. This is the case in particular for the porous medium equations $\partial_t \rho = \Delta \rho^m$, $m > 1$. To see this, note that our lift gives $v = -\frac{1}{\rho} \nabla f(\rho) = -\nabla g(\log \rho)$ where $g(x) = \int_0^x f'(e^y) dy$. Thus $v$ is a gradient field, hence equals the Otto lift.
\end{remark}

From Proposition \ref{prop:jacobian_L_plus_N_form}, the evolution equation for the submap log-Jacobian $\jacobian := \jacobian_{[t, \tau_k]}$ in this case becomes the nonlinear advection-reaction-diffusion equation
\begin{equation*}
    \partial_t \jacobian - \nabla \log \jsource \cdot \nabla \jacobian = \Delta \jacobian + |\nabla \jacobian|^2 + \Delta \log \jsource,
\end{equation*}
in terms of the source term $\jsource := \jsource_k$ (defined in \eqref{eq:j_source_term_def}). With the exponential integrator $\intJac$ defined in Section \ref{section:jacobian_time_stepper}, the leading-order linear term $\Delta$ is exponentiated exactly, removing the second-order CFL condition. The nonlinear remainder $\mathcal{N}(\jacobian, \jsource) = |\nabla \jacobian|^2 + \nabla \log \jsource \cdot \nabla \jacobian + \Delta \log\jsource$ is evolved explicitly.

\begin{figure}[h!]
    \centering
    \begin{subfigure}[t]{0.32\textwidth}
        \includegraphics[width=\textwidth]{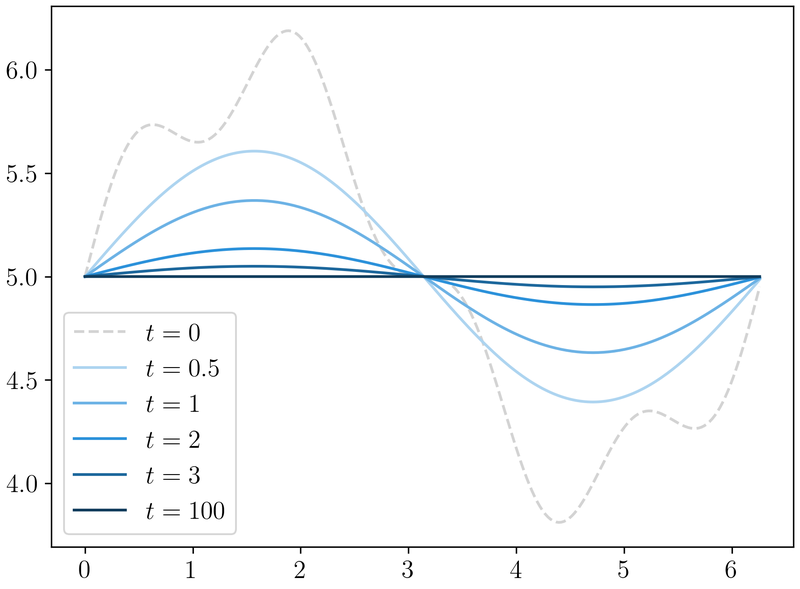}
        \caption*{$\rho(t)$}
    \end{subfigure}
    \begin{subfigure}[t]{0.32\textwidth}
        \includegraphics[width=\textwidth]{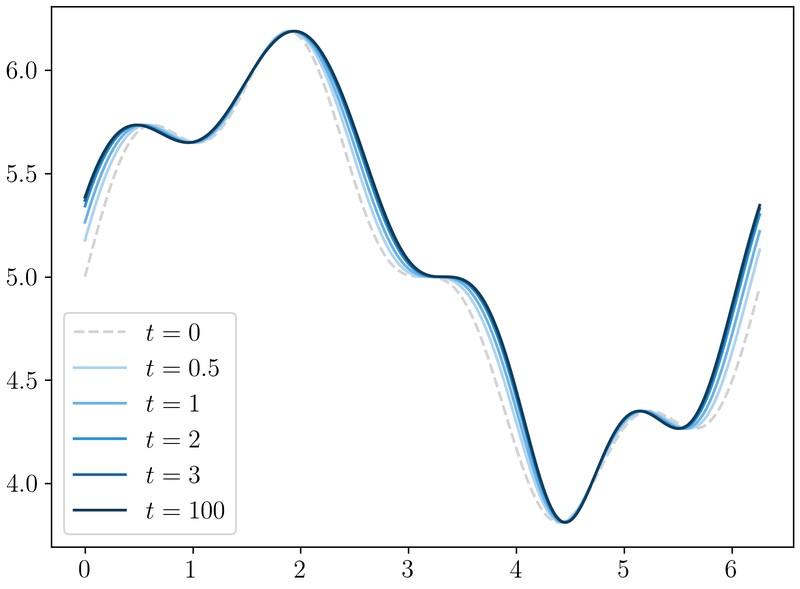}
        \caption*{$\rho_0 \circ \map(t)$}
    \end{subfigure}
    \begin{subfigure}[t]{0.32\textwidth}
        \includegraphics[width=\textwidth]{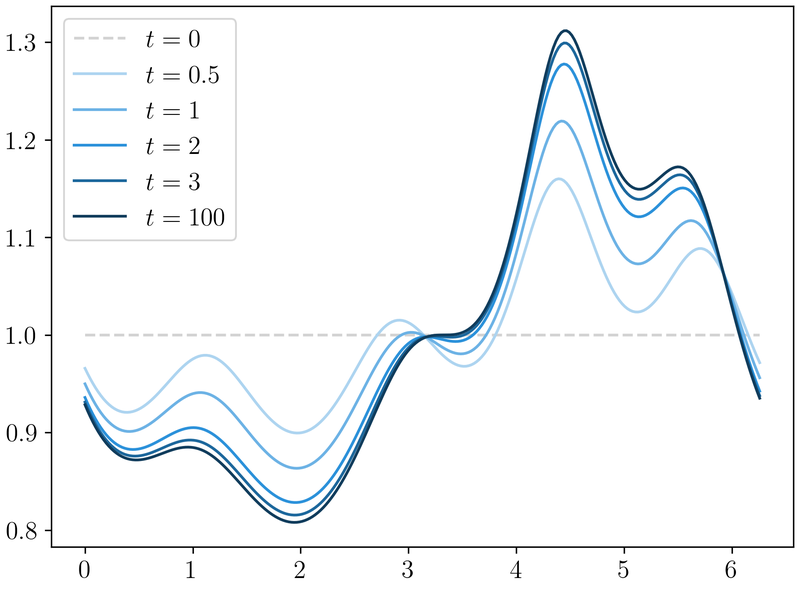}\caption*{$e^{\jacobian(t)}$}
    \end{subfigure}
    \caption{Plots of the solution (left) along with the advective (middle) and compression (right) components of the numerical solutions for the one-dimensional heat equation experiment. The dominant evolution is concentrated to the compression component of the solution.}
    \label{fig:heat_plots}
\end{figure}

To illustrate the nonlinear splitting induced by the group action in the case of pure diffusion, we solve a one-dimensional heat equation. We take the initial condition
\begin{equation*}
    \rho_0 (x) = 5 + \sin(x) + \frac{1}{4} \sin(4 x),
\end{equation*}
and run the method with $N = 32$ and $\Delta t = 1/64$ until $T = 100$ with a single submap. In Figure \ref{fig:heat_plots}, we plot at various times the solution $\rho$ (left), the transport component $\rho_0 \circ \map$ (middle), and the compression component $e^\jacobian$ (right). We see that the transport component $\rho_0 \circ \map$ contributes very little to the overall dynamics; the majority of the evolution is due to the Jacobian. Indeed, the solution of the heat equation converges towards its mean as $t \rightarrow \infty$. For the majority contribution to come from the transport component $\rho_0 \circ \map$, we would need the map to become singular as $t \rightarrow \infty$. Instead we see the map stabilize, which is expected since the exact lift $v = -\nabla \log \rho$ decays (since $\rho \to \text{const.}$) as $t \rightarrow \infty$. This example demonstrates that the compression component dominates the evolution for purely dissipative linear systems.

\subsection{Phase-Field Equations}

We next consider two nonlinear dissipative examples. These are the Allen--Cahn reaction-diffusion equation \cite{allen_microscopic_1979},
\begin{equation}
    \label{eq:ac_equation}
    \partial_t \solution = \epsilon^2 \Delta \solution  - f(\solution), \quad\quad f(\solution) = \solution^3 - \solution,
\end{equation}
and the related fourth-order Cahn--Hilliard equation \cite{cahn_free_1958},
\begin{equation}
    \label{eq:ch_equation}
    \partial_t \solution = D \Delta ( f(\solution) - \epsilon^2 \Delta \solution ).
\end{equation}
These equations are models of phase separation in a two-component system. The phase field $\solution$ evolves toward the minima of the double-well potential $\frac{1}{4} (\solution^2 - 1)^2$: both equations are gradient flows of the Ginzburg-Landau free energy functional \cite{shen_numerical_2010}
\begin{equation}
    E(\solution) = \int_{\mathbb{T}^2} \left( \frac{\epsilon^2}{2} |\nabla \solution|^2 + \frac{1}{4} (\solution^2 - 1)^2 \right).
\end{equation}
The Allen--Cahn equation is the $L^2$ gradient flow of $E$, while the Cahn--Hilliard equation is the gradient flow in $H^{-1}$ (with mobility $D$). In both equations, the evolution is characterized by a rapid initial separation of the unstable mixture into two phases, followed by a coarsening over longer time scales.
The Allen--Cahn equation has varying mass $\dot{\mass} \neq 0$, and its phase transition interfaces evolve by a curvature-driven coarsening, at an asymptotic coarsening rate of $t^{1/2}$.
On the other hand, the Cahn--Hilliard equation, being an $H^{-1}$ gradient flow, conserves mass $\dot{\mass} = 0$, and its interfaces undergo a conservative coarsening at a rate of $t^{1/3}$.

We now lift these equations in our framework, beginning with the Allen--Cahn equation. In the coordinate $\rho = \solution + c > 0$, the lift $(v, \xi) = \liftingOperator(\rho)$ reads
\begin{equation*}
    v = -\epsilon^2 \nabla \log \rho + \frac{1}{\rho} \nabla \Delta^{-1} (\xi \rho + f(\rho - c)), \quad\quad \xi = -\frac{\int_{\mathbb{T}^2} f(\rho - c)}{\int_{\mathbb{T}^2} \rho}.
\end{equation*}
From Proposition \ref{prop:jacobian_L_plus_N_form}, the log-Jacobian evolution equation reads
\begin{equation*}
    \partial_t \jacobian - \epsilon^2 \nabla \log \jsource \cdot \nabla \jacobian = \epsilon^2 \Delta \jacobian + \epsilon^2 |\nabla \jacobian|^2 + \epsilon^2 \Delta \log \jsource - \frac{1}{\jsource e^\jacobian} \left( G_0 - \nabla \log \jsource \cdot \nabla \Delta^{-1} G_0 \right),
\end{equation*}
where we have written $G_0 := f(\jsource e^\jacobian - c) + \xi \jsource e^\jacobian$. We note that the final term is zeroth order $\jacobian$, as well as nonlocal. We find the leading-order term $\mathcal{L} = \epsilon^2 \Delta$ to be exponentiated exactly, and a lower-order nonlinear remainder to be treated explicitly.

For the Cahn--Hilliard equation, we compute the lift
\begin{equation*}
    v = \frac{D}{\rho} \nabla (\epsilon^2 \Delta \rho - f(\rho - c)),
\end{equation*}
which is now local in $\rho$. In this case, we obtain an evolution equation for $\jacobian$ of the form
\begin{equation*}
    \partial_t \jacobian = -D \epsilon^2 \Delta^2 \jacobian + \mathcal{N}(\jacobian, \jsource),
\end{equation*}
where $\mathcal{N}$ is a local third-order operator in $\jacobian$.

\subsubsection{Mean Curvature Flow of a Disk}

We begin with a mean curvature flow in the Allen--Cahn equation. In this equation, the transitions between phases have a typical profile of the form $\tanh(s/\ell)$, with $s \in \mathbb{R}$ the normal coordinate and the characteristic length scale $\ell = \sqrt{2} \epsilon$ \cite{evans1992phase}. Over long time scales, the phase interfaces evolve by approximate mean curvature flow. The mean curvature flow being only approximate is a byproduct of both the Allen--Cahn equation and the fact that we are in a periodic domain. Thus, we take as initial condition a periodization of the following function,
\begin{equation*}
    \solution_0 (x) = \tanh \left( \frac{R_0 - |x|}{\ell} \right), \quad\quad R_0 = 1,
\end{equation*}
representing a disk of radius 1 with a tanh profile. We expect the disk to shrink in radius at a linear rate until extinction. We take $N = 128$, $N_\epsilon = 96$, $\Delta t = 1/16$, $c = 2$ and $\epsilon = 10^{-1}$, with uniform remapping times $\Delta \tau = 10$. The exact mean curvature flow of a disk in our time scale would have a radius over time of $R(t) = \sqrt{R_0^2 - \ell^2 t}$, giving an expected extinction time $t \approx R_0^2/\ell^2 = 50$.

\begin{figure}
    \centering
    \begin{subfigure}[c]{0.15\textwidth}
        \includegraphics[width=\textwidth]{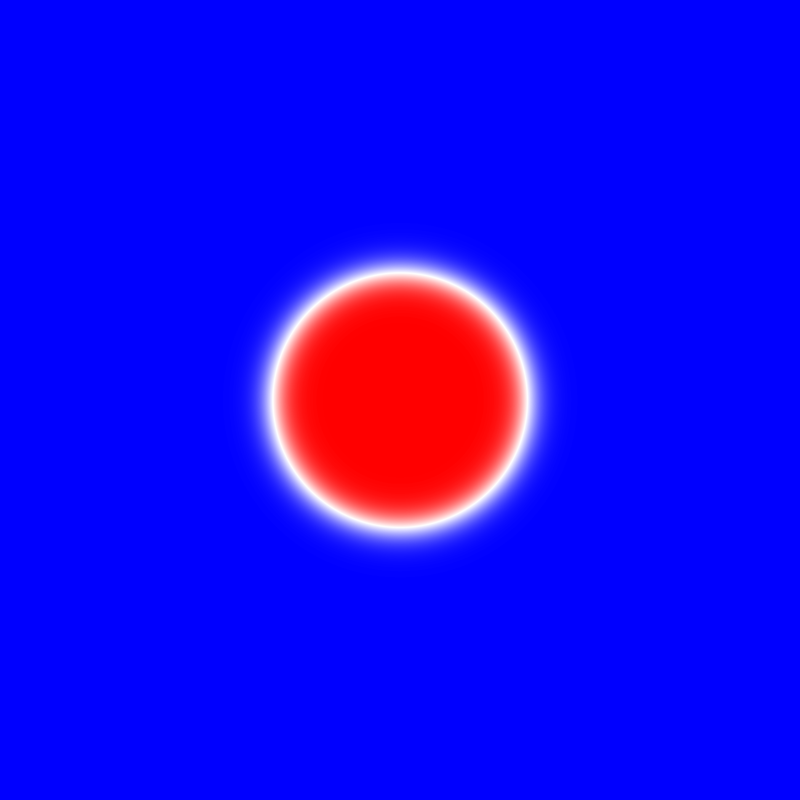}
    \end{subfigure}
    \begin{subfigure}[c]{0.15\textwidth}
        \includegraphics[width=\textwidth]{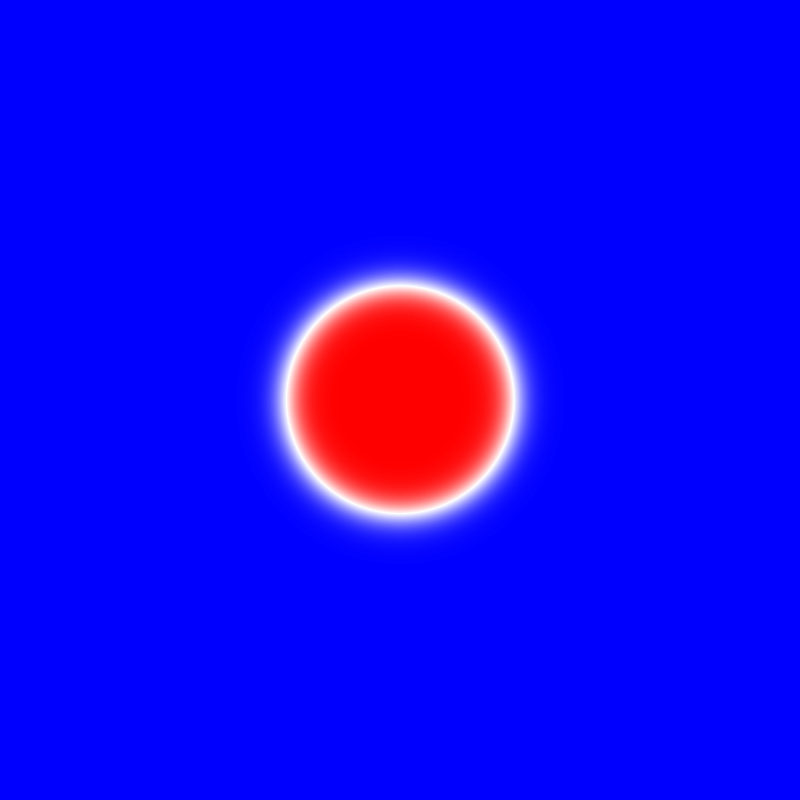}
    \end{subfigure}
    \begin{subfigure}[c]{0.15\textwidth}
        \includegraphics[width=\textwidth]{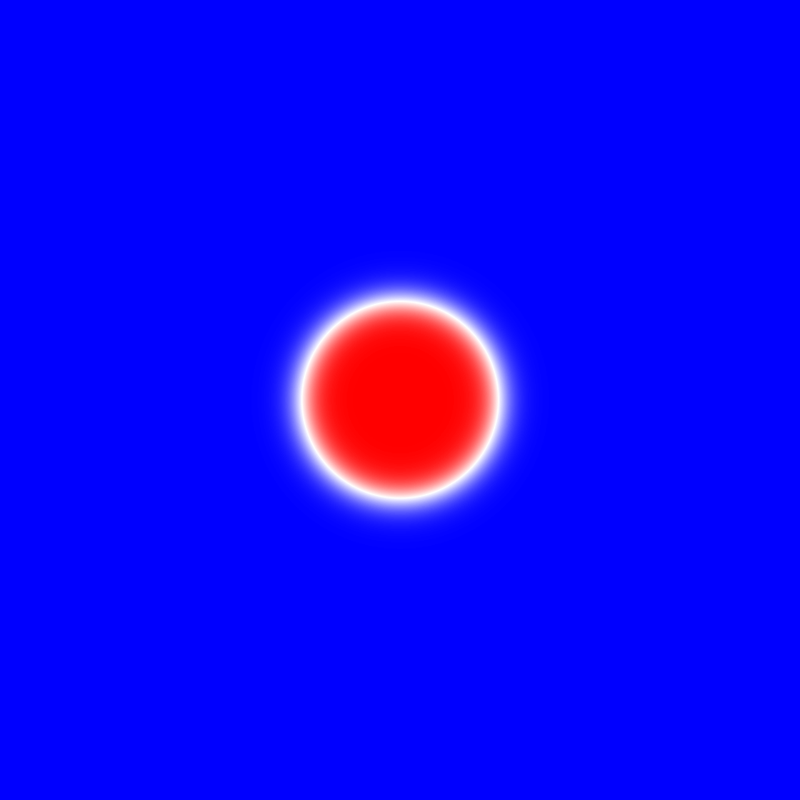}
    \end{subfigure}
    \begin{subfigure}[c]{0.15\textwidth}
        \includegraphics[width=\textwidth]{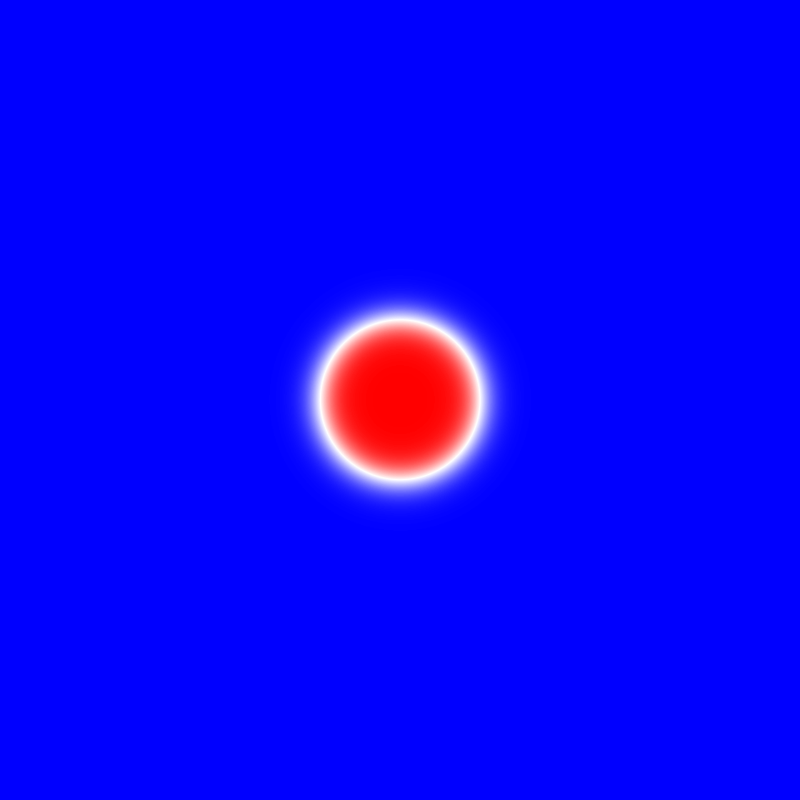}
    \end{subfigure}
    \begin{subfigure}[c]{0.15\textwidth}
        \includegraphics[width=\textwidth]{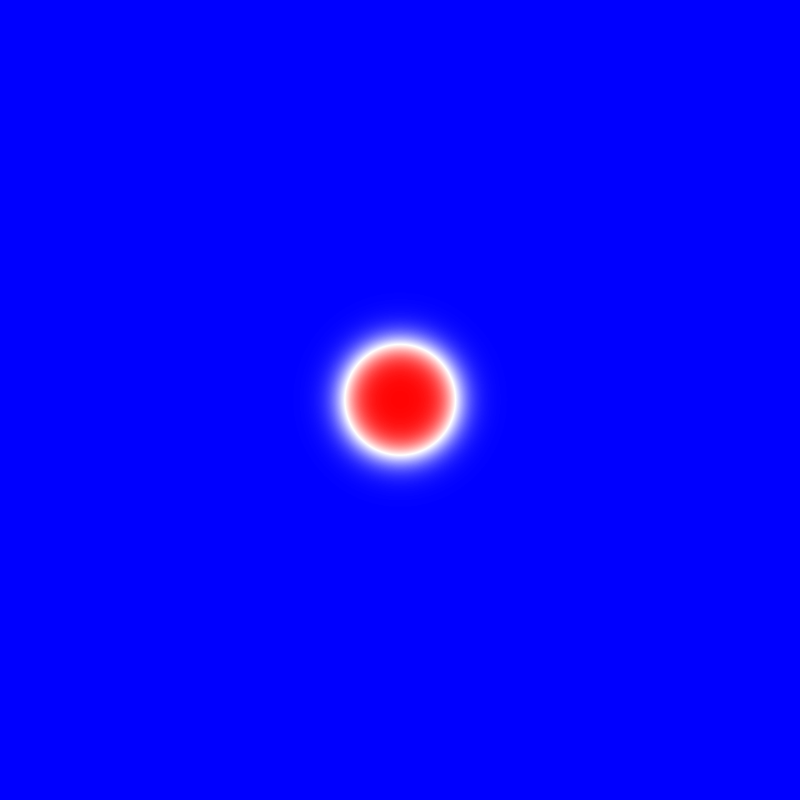}
    \end{subfigure}
    \begin{subfigure}[c]{0.15\textwidth}
        \includegraphics[width=\textwidth]{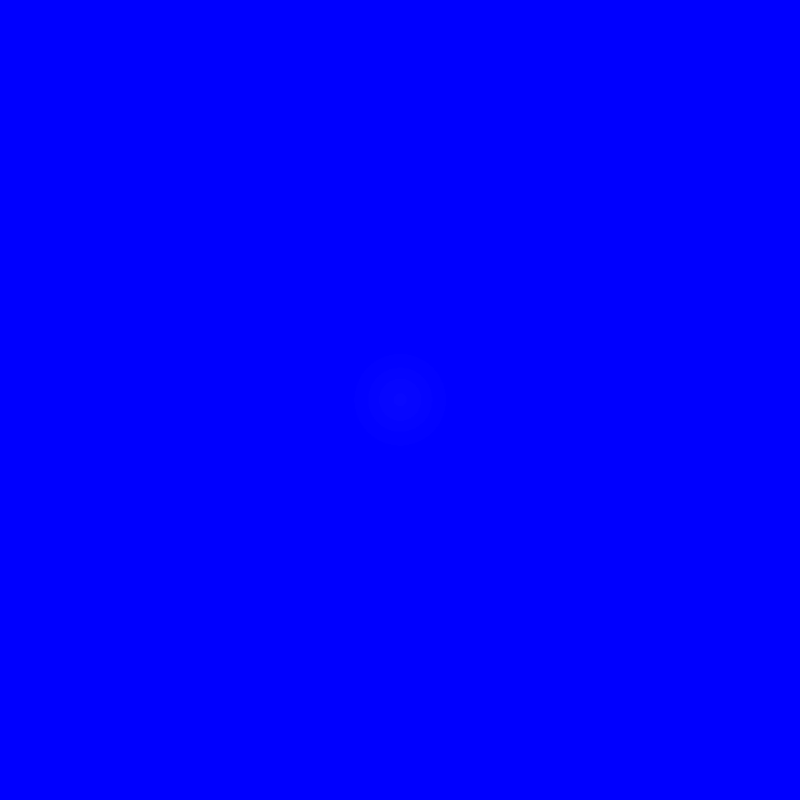}
    \end{subfigure}
    \begin{subfigure}[c]{0.02267\textwidth}
        \includegraphics[width=\textwidth]{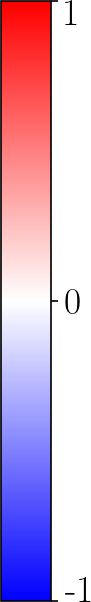}
    \end{subfigure}
    \\[0.25em]
    \begin{subfigure}[c]{0.15\textwidth}
        \includegraphics[width=\textwidth]{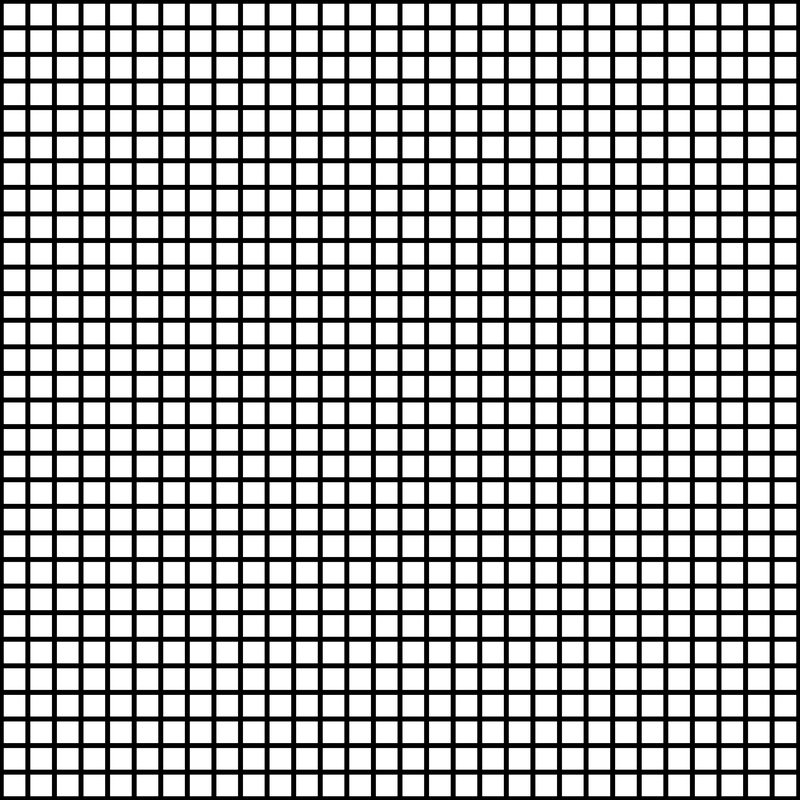}
        \caption*{$t = 0$}
    \end{subfigure}
    \begin{subfigure}[c]{0.15\textwidth}
        \includegraphics[width=\textwidth]{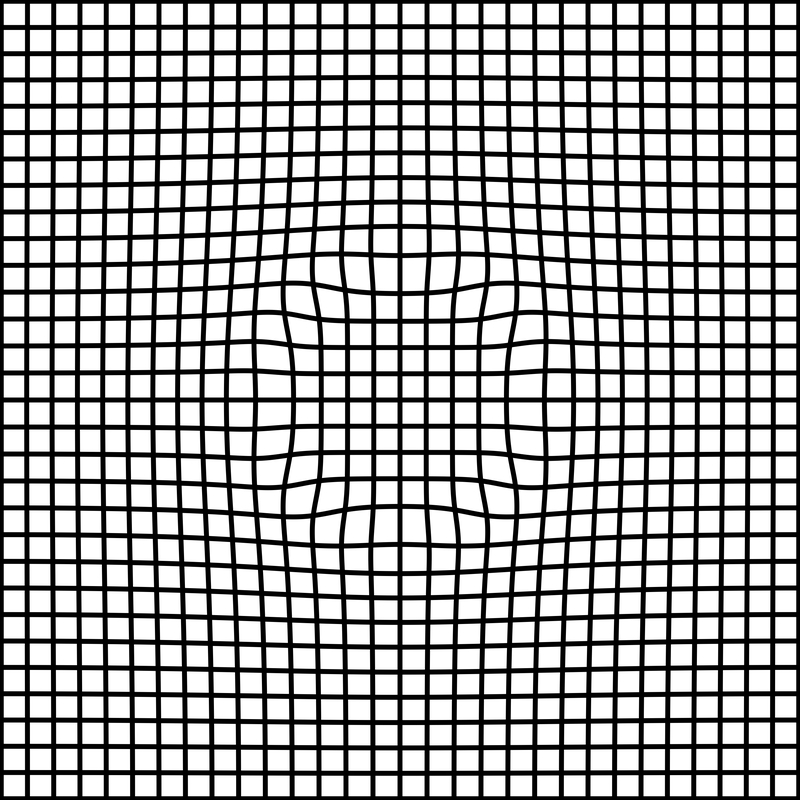}
        \caption*{$t = 10$}
    \end{subfigure}
    \begin{subfigure}[c]{0.15\textwidth}
        \includegraphics[width=\textwidth]{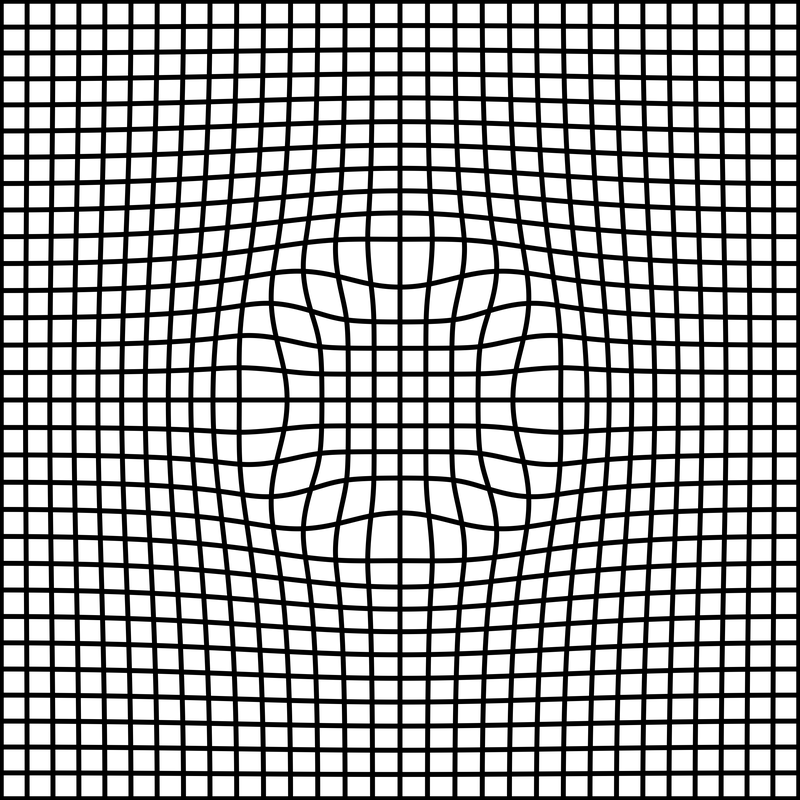}
        \caption*{$t = 20$}
    \end{subfigure}
    \begin{subfigure}[c]{0.15\textwidth}
        \includegraphics[width=\textwidth]{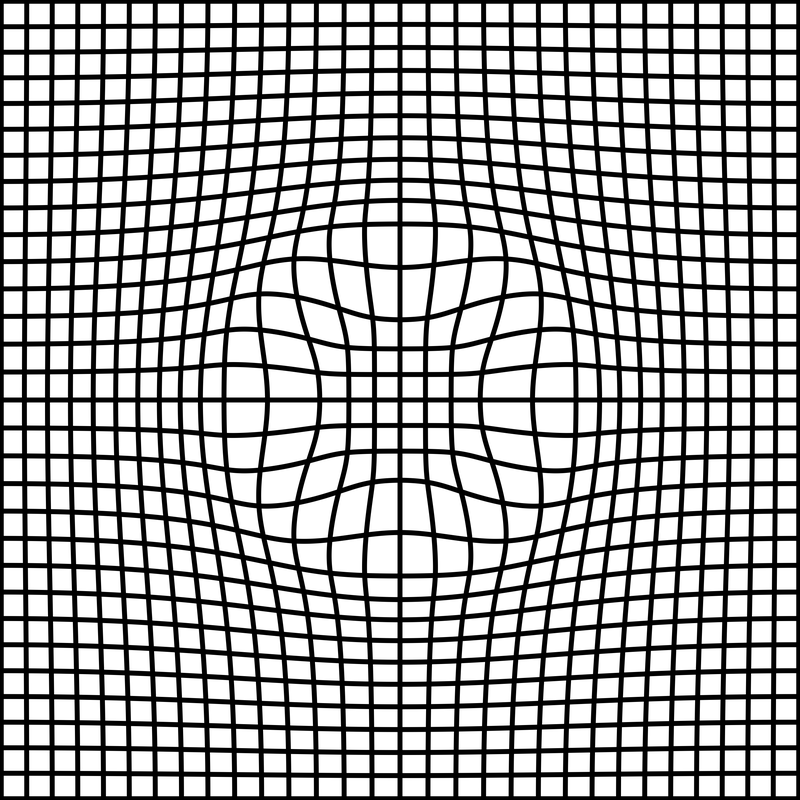}
        \caption*{$t = 30$}
    \end{subfigure}
    \begin{subfigure}[c]{0.15\textwidth}
        \includegraphics[width=\textwidth]{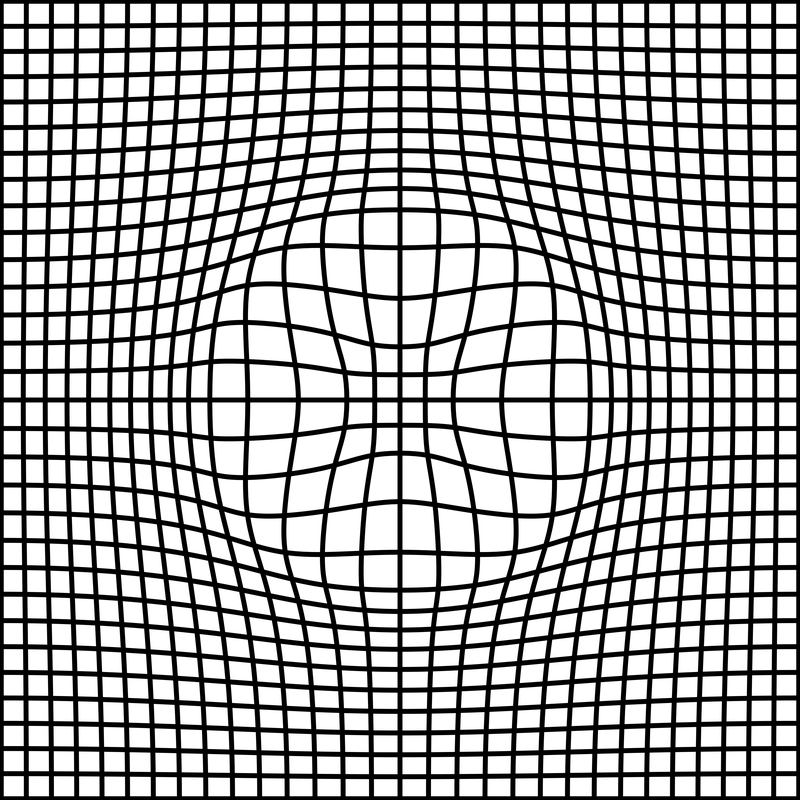}
        \caption*{$t = 40$}
    \end{subfigure}
    \begin{subfigure}[c]{0.15\textwidth}
        \includegraphics[width=\textwidth]{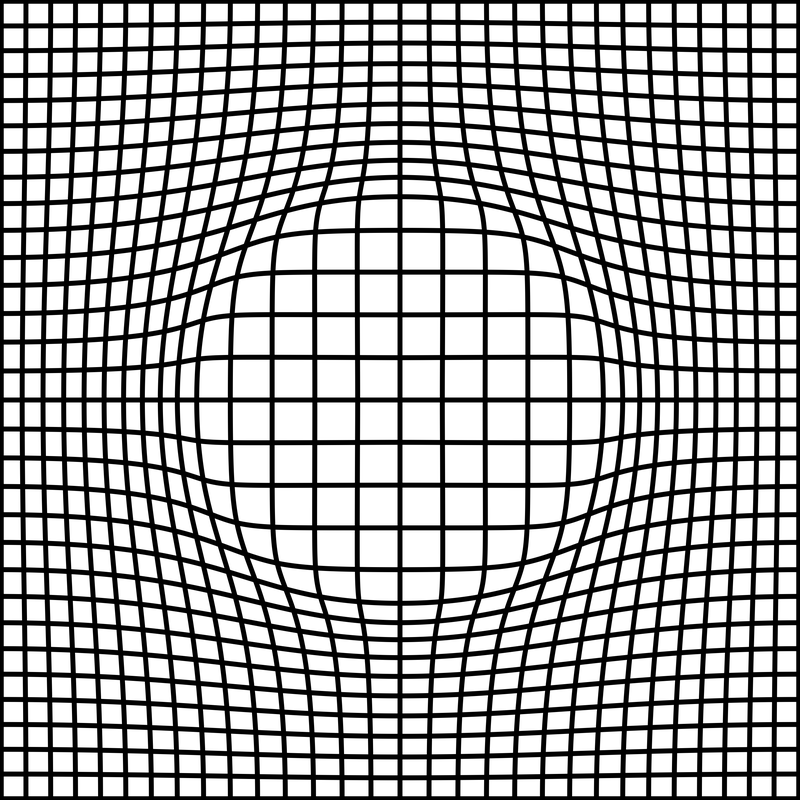}
        \caption*{$t = 50$}
    \end{subfigure}
    \hspace*{0.02267\textwidth}
    \caption{Mean curvature flow of a $\tanh$-profile disk in Allen--Cahn until extinction time. \emph{Top row:} Solution $\solution$. \emph{Bottom row:} Deformation of the space generated by the lifted evolution.}
    \label{fig:allen_cahn_mcf_components}
\end{figure}

We plot the solution $\solution$ over time until extinction in the top row of Figure \ref{fig:allen_cahn_mcf_components}. In the bottom row of Figure \ref{fig:allen_cahn_mcf_components}, we plot the deformation of the space generated by the lifted evolution. This example illustrates the need for the lift to in general be a nonlocal operator of the solution. Indeed, for the disk to shrink and eventually go extinct, we see that mass must flow outward, far past the active region of the solution. We also see that, since the solution converges to a constant, the diffusion of the disk by the density action must be offset by a varying total mass $\mass$.

\begin{figure}[h!]
    \centering
    \begin{subfigure}[c]{0.32\textwidth}
        \includegraphics[width=\textwidth]{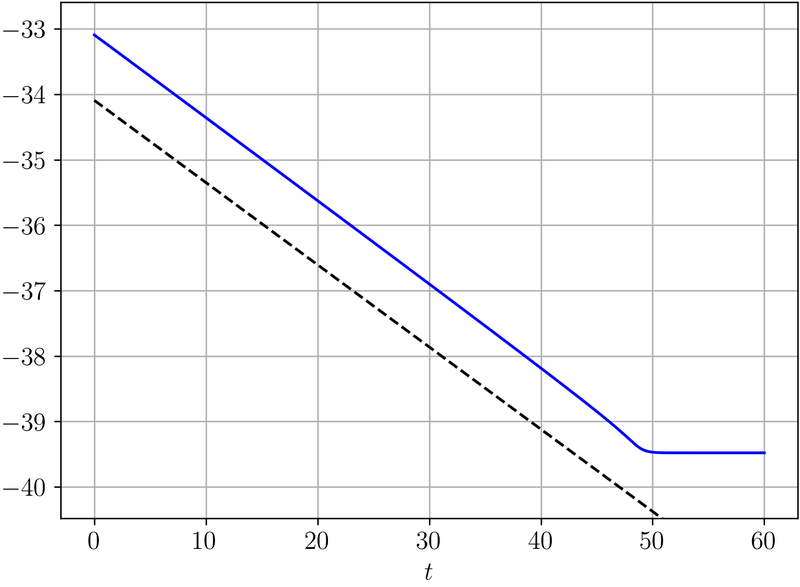}
    \end{subfigure}
    \hspace*{1em}
    \begin{subfigure}[c]{0.32\textwidth}
        \includegraphics[width=\textwidth]{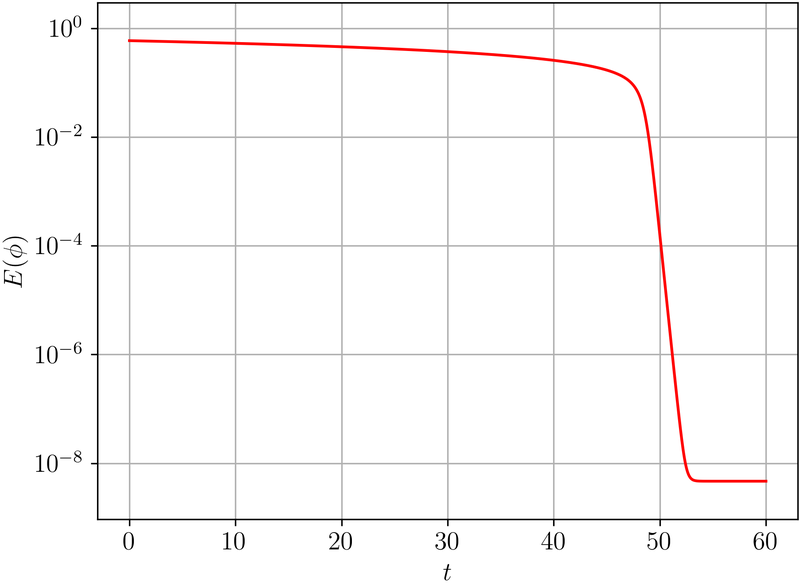}
    \end{subfigure}
    \caption{Statistics in the Allen--Cahn mean curvature disk test. \emph{Left:} Mass of the solution $\int_{\mathbb{T}^2} \solution$ over time (blue) versus the expected rate $-2\pi\ell^2$ (dashed). \emph{Right:} Energy $E(\solution)$ over time.}
    \label{fig:ac_mcf_mass_energy}
\end{figure}

In exact mean curvature flow, the radius over time $R(t) = \sqrt{R_0^2 - \ell^2 t}$ implies a constant rate of change of the area of the disk, $\frac{dA}{dt} = -\pi \ell^2$. The area is related to the integral of the solution via
\begin{equation*}
    \int_{\mathbb{T}^2} (-1) + 2 A \approx \int_{\mathbb{T}^2} \solution,
\end{equation*}
thus we expect the rate of change in mass $\frac{d}{dt} \int_{\mathbb{T}^2} \solution \approx -2\pi\ell^2$. In Figure \ref{fig:ac_mcf_mass_energy} (left), we plot the integral of the solution over time and see that it agrees with the expected rate until just before extinction, showing that the varying-mass case is correctly handled by the method.
We also plot the energy $E(\solution)$ over time in Figure \ref{fig:ac_mcf_mass_energy} (right), where we see the energy drop off significantly at the expected extinction time.

\subsubsection{Random Initial Condition}

We next perform an experiment for both the Allen--Cahn and Cahn--Hilliard equations, which consists of an identical random initial condition. The parameters for both equations are collected in Table \ref{fig:ac_ch_parameters}. For both equations, we take the initial condition $\solution_0$ to consist of 256 random Fourier modes, with phases sampled uniformly and decaying magnitude $|\hat{\solution}_0 (k)| = \mathcal{O}(|k|^{-1})$, normalized so that $\int_{\mathbb{T}^2} \solution_0 = 0$ and $\lVert \solution_0 \rVert_{L^\infty} = 1$. The method is run for both equations until a final time $T = 400$.

{
\renewcommand{\arraystretch}{1.2}
\begin{table}[h]
    \centering
    \begin{tabular}{c|c|c|c|c|c|c|c}
        & Constants & $c$ & $N$ & $N_\epsilon$ & $\Delta t$ & $\Delta \tau$ & Number of submaps \\
        \hline
        Allen--Cahn & $\epsilon = 10^{-2}$ & $2$ & $384$ & $192$ & $1/8$ & $10$ & 40 \\
        \hline
        Cahn--Hilliard & \makecell{$D = 10^{-2}$ \\ $\epsilon = 10^{-1.5}$} & $2$ & $512$ & $256$ & $1/512$ & $5$ & 80
    \end{tabular}
    \caption{Parameters in the random initial condition test for Allen--Cahn and Cahn--Hilliard.}
    \label{fig:ac_ch_parameters}
\end{table}
}

\begin{figure}
    \centering
    \begin{tabular*}{\linewidth}{@{}c@{\hspace{0.5em}}c@{\extracolsep{\fill}}c@{\extracolsep{\fill}}c@{\extracolsep{\fill}}c@{\extracolsep{\fill}}c@{\extracolsep{\fill}}c@{}}
        \multirow{2}{*}{%
            \smash{\rotatebox[origin=c]{90}{\small Allen-Cahn}}%
        }

        & \includegraphics[width=0.182\textwidth]{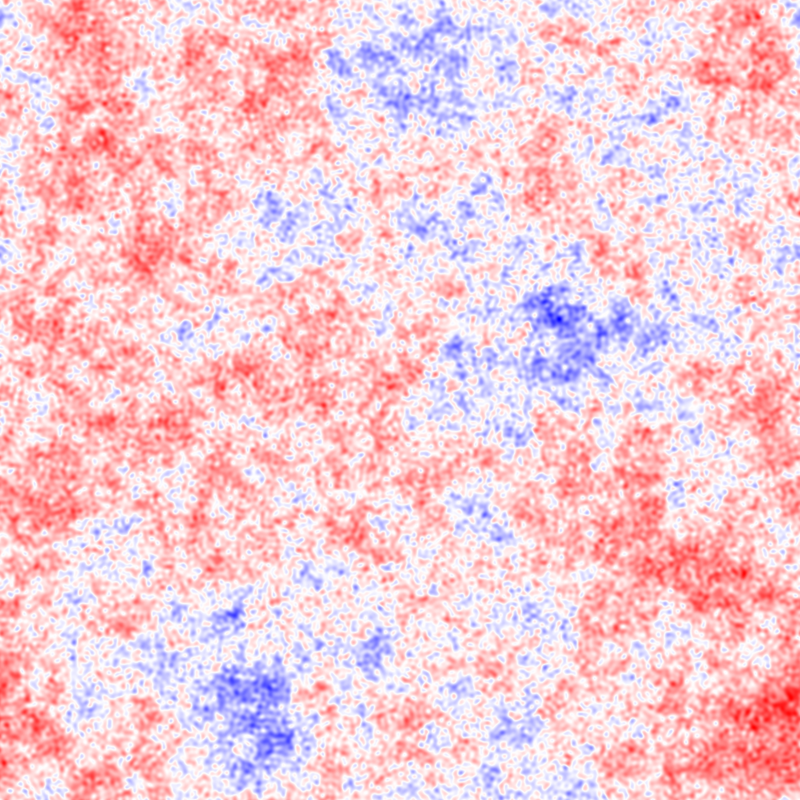}
        & \includegraphics[width=0.182\textwidth]{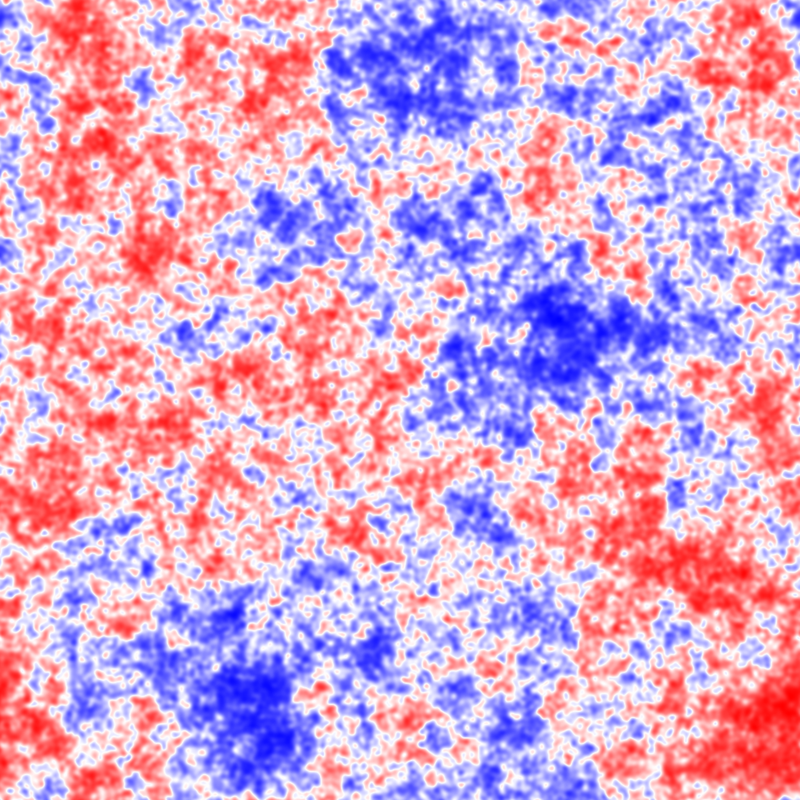}
        & \includegraphics[width=0.182\textwidth]{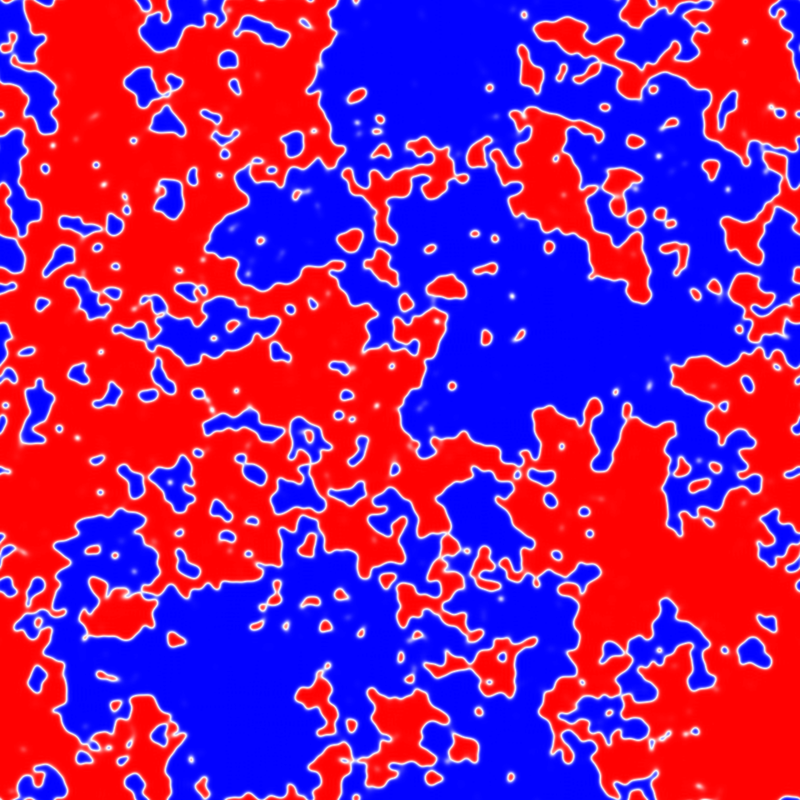}
        & \includegraphics[width=0.182\textwidth]{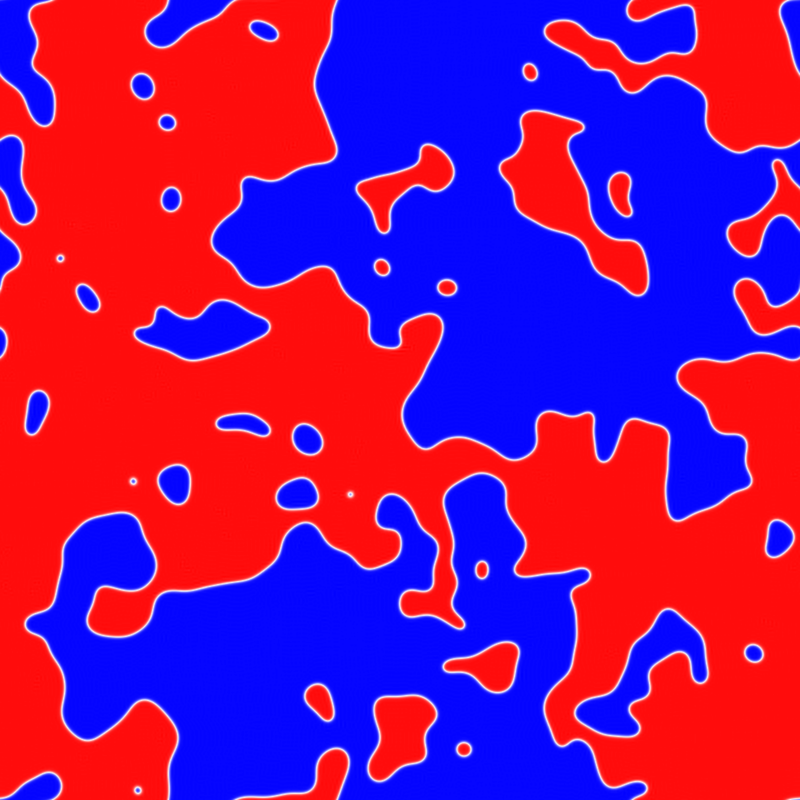}
        & \includegraphics[width=0.182\textwidth]{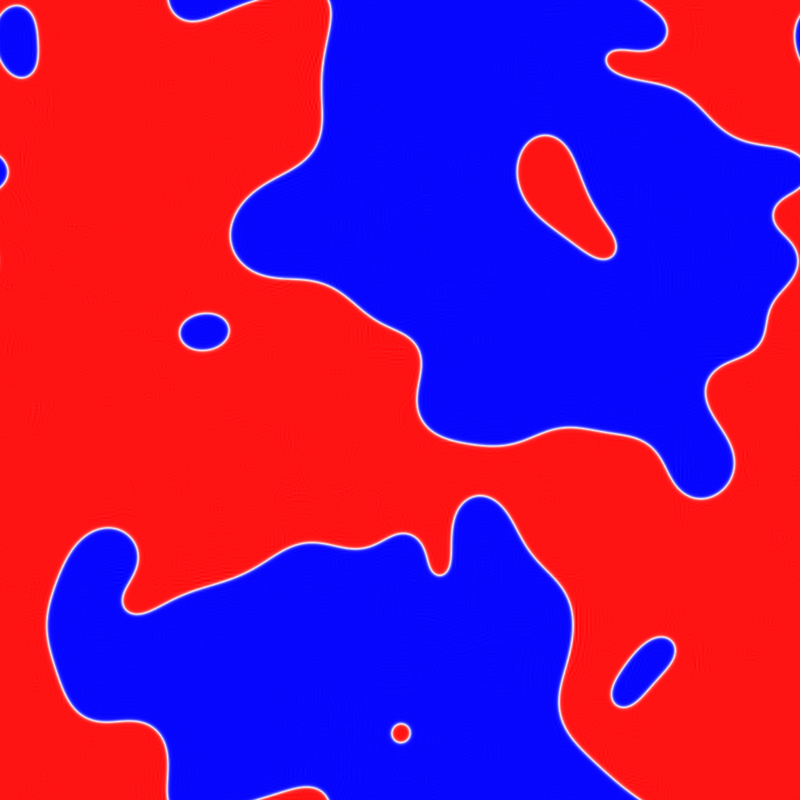}
        & \includegraphics[width=0.02751\textwidth]{allen_cahn/colorbar_thin.png}
        \\

        & \includegraphics[width=0.182\textwidth]{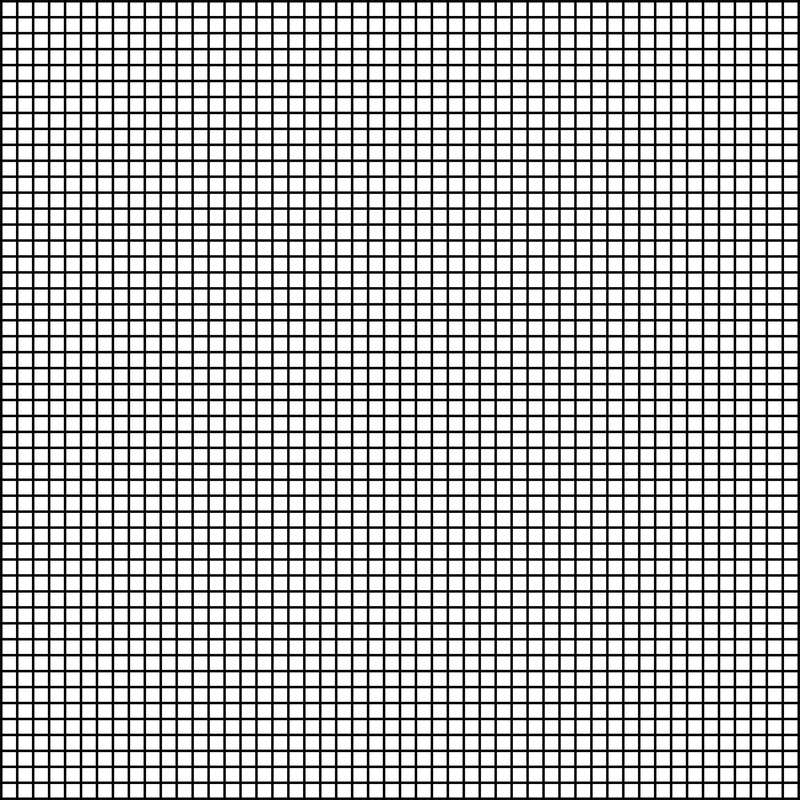}
        & \includegraphics[width=0.182\textwidth]{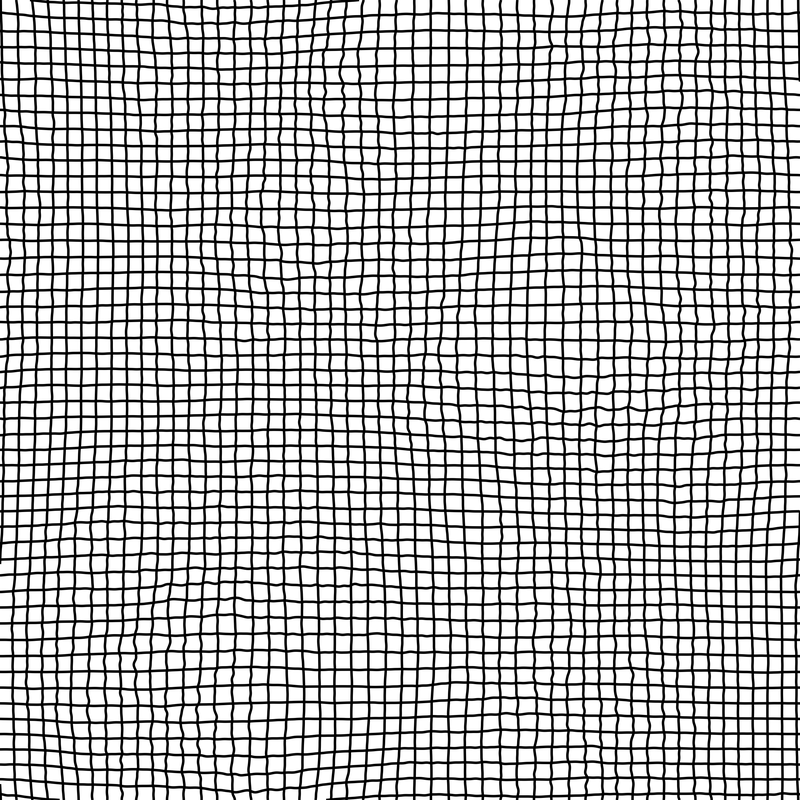}
        & \includegraphics[width=0.182\textwidth]{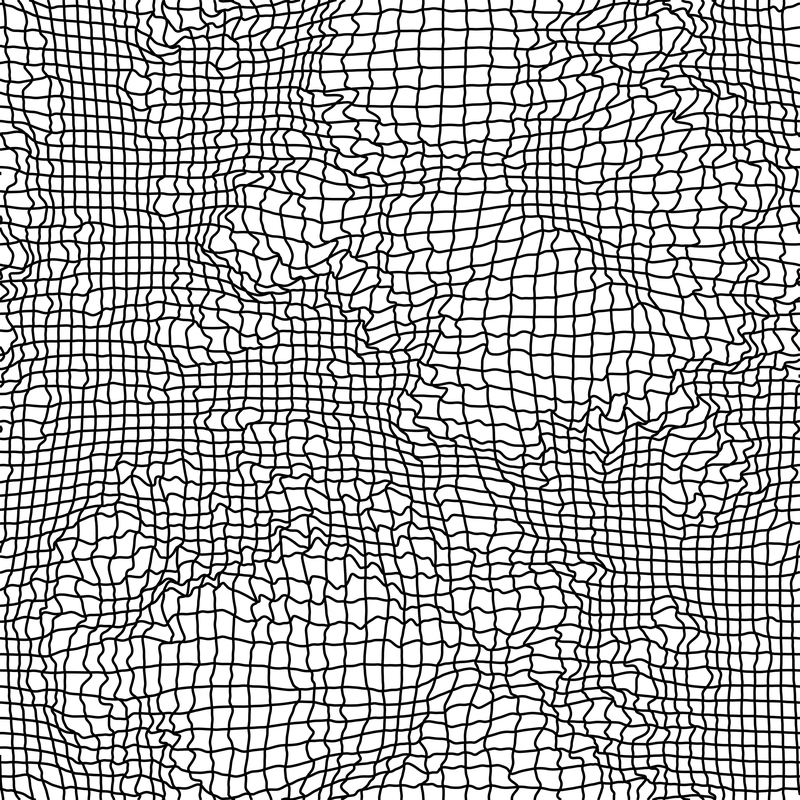}
        & \includegraphics[width=0.182\textwidth]{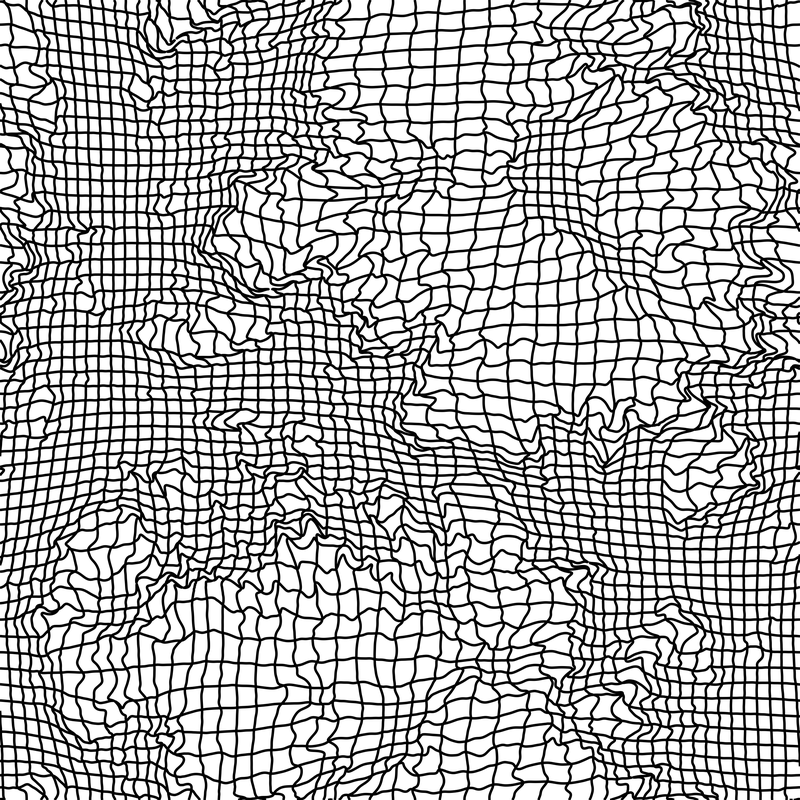}
        & \includegraphics[width=0.182\textwidth]{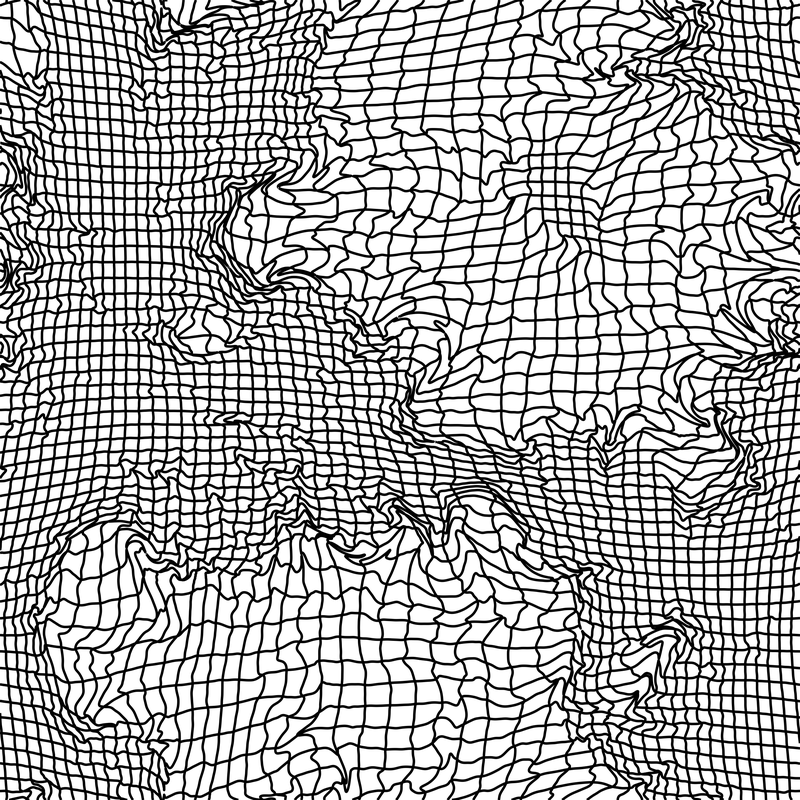}
        \\
        & {\small $t = 0$} & {\small $t = 1$} & {\small $t = 5$} & {\small $t = 50$} & {\small $t = 400$}
        \\[0.2em]

        \multirow{2}{*}{%
            \smash{\rotatebox[origin=c]{90}{\small Cahn-Hilliard}}%
        }
        
        & \includegraphics[width=0.182\textwidth]{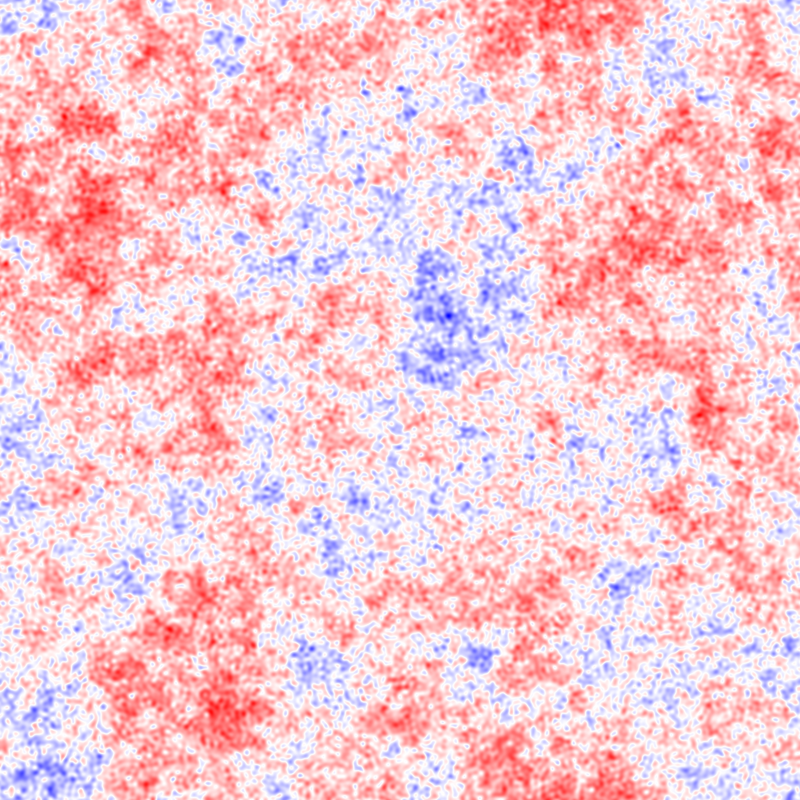}
        & \includegraphics[width=0.182\textwidth]{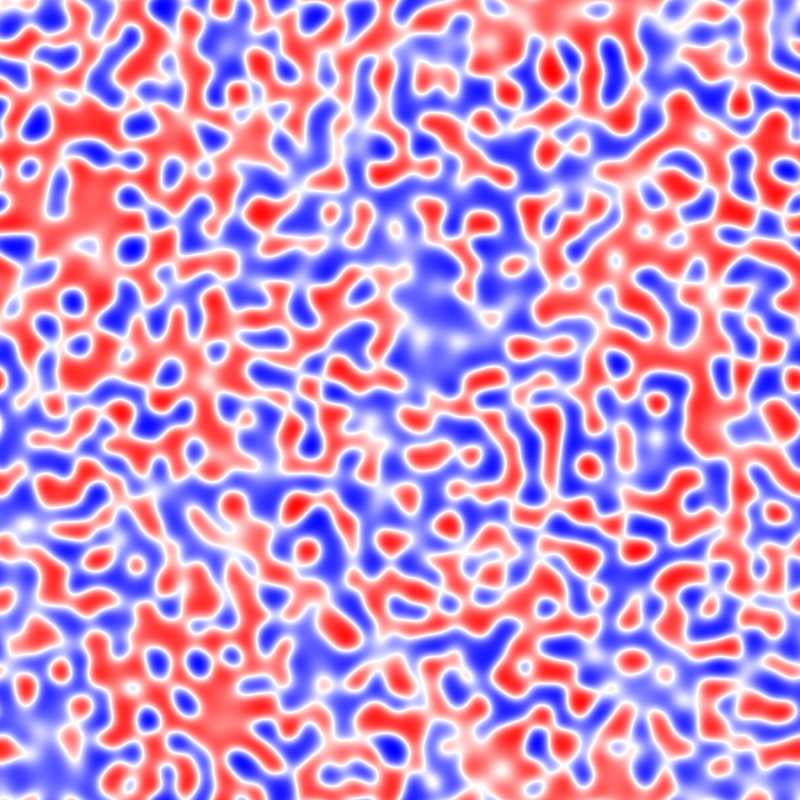}
        & \includegraphics[width=0.182\textwidth]{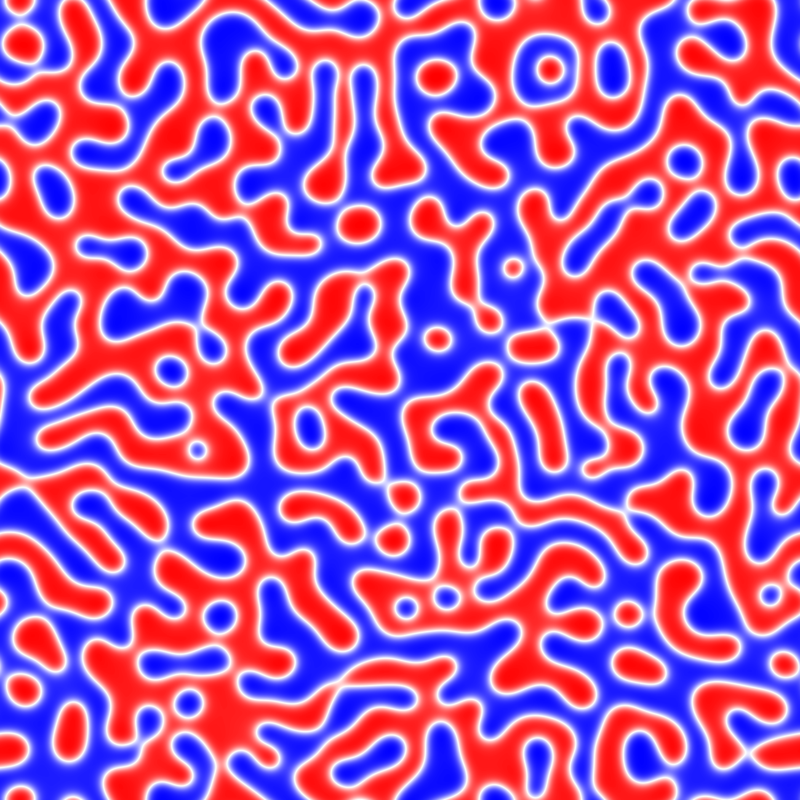}
        & \includegraphics[width=0.182\textwidth]{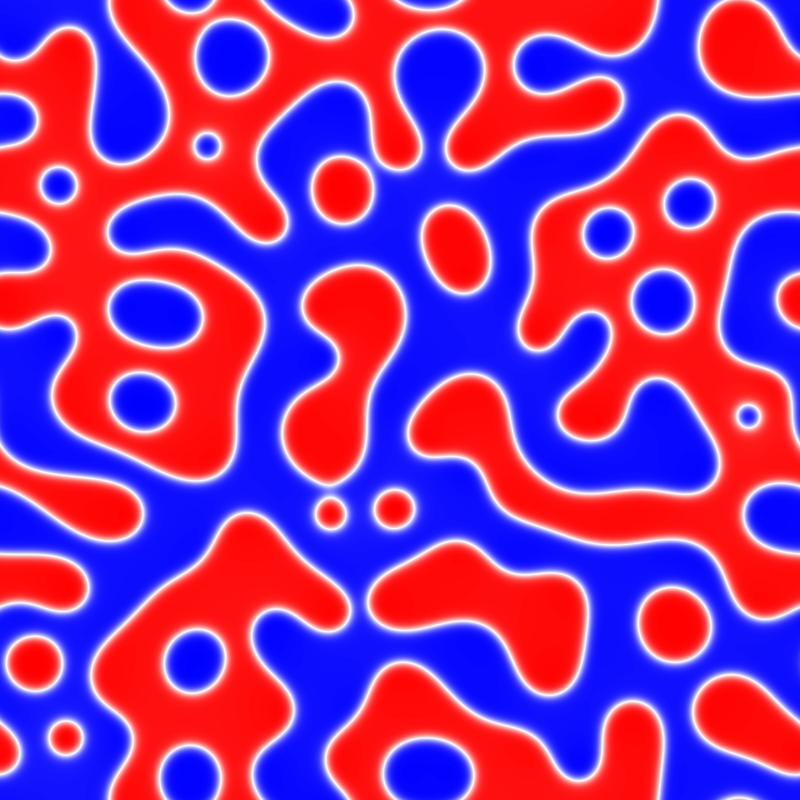}
        & \includegraphics[width=0.182\textwidth]{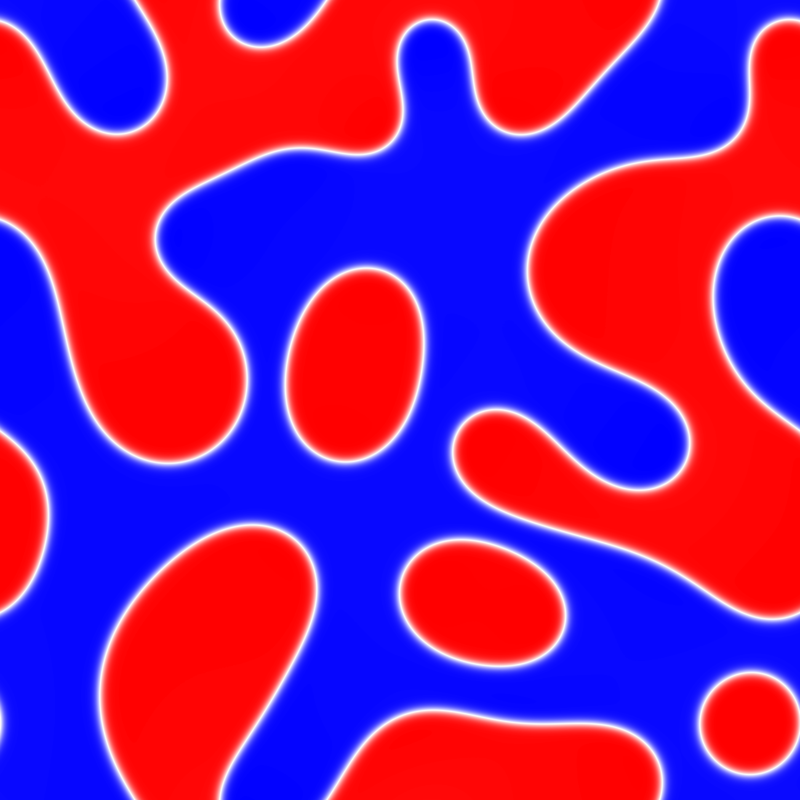}
        & \includegraphics[width=0.02751\textwidth]{allen_cahn/colorbar_thin.png}
        \\

        & \includegraphics[width=0.182\textwidth]{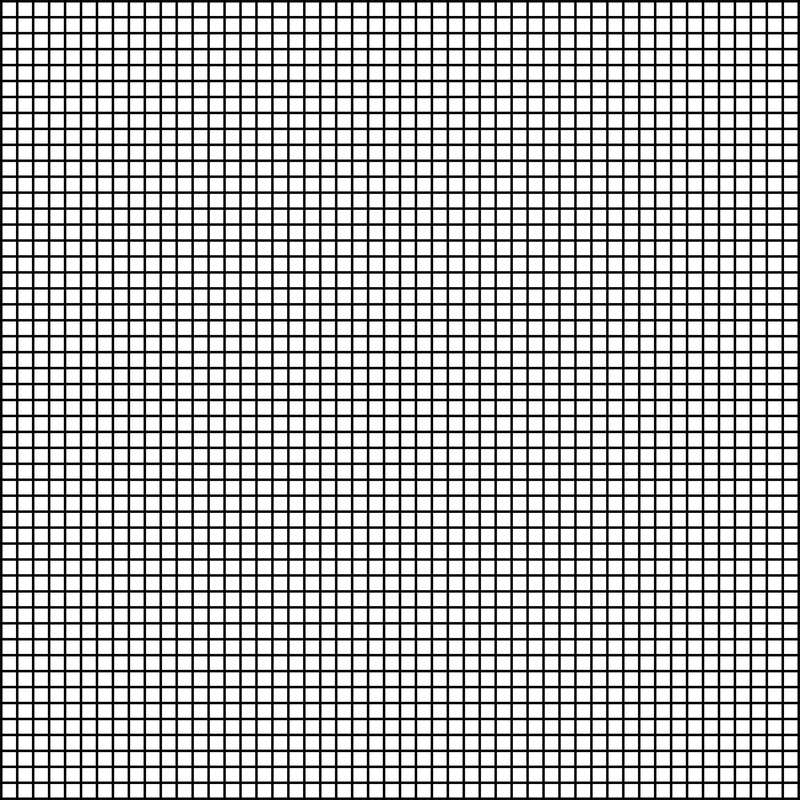}
        & \includegraphics[width=0.182\textwidth]{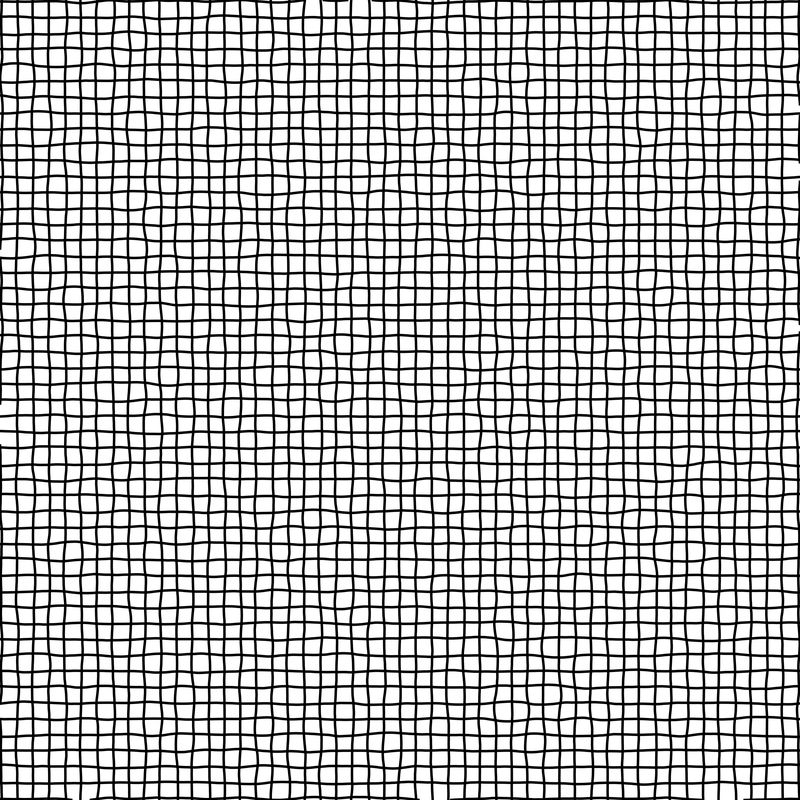}
        & \includegraphics[width=0.182\textwidth]{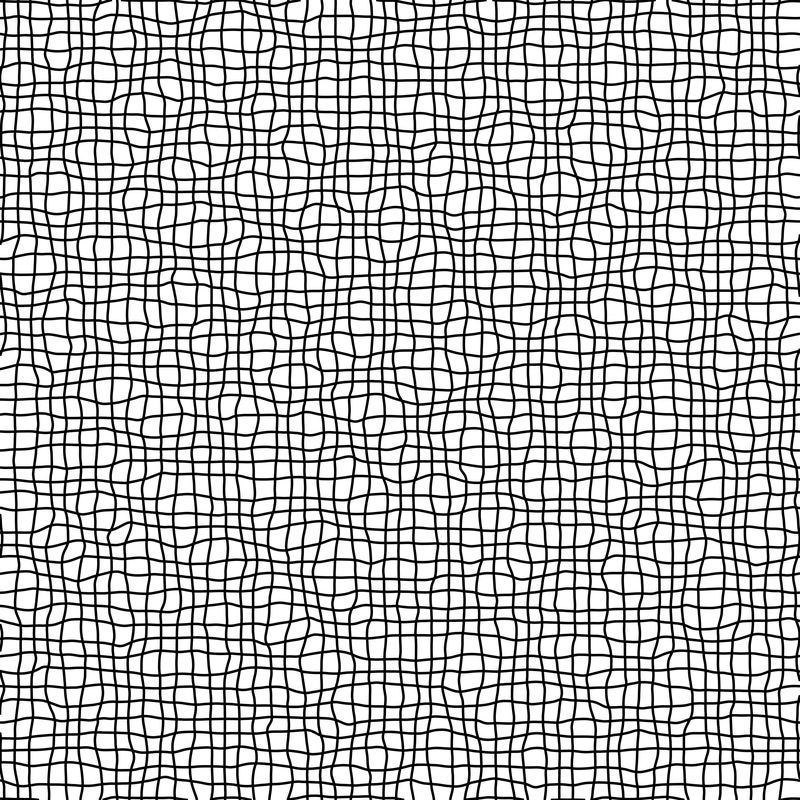}
        & \includegraphics[width=0.182\textwidth]{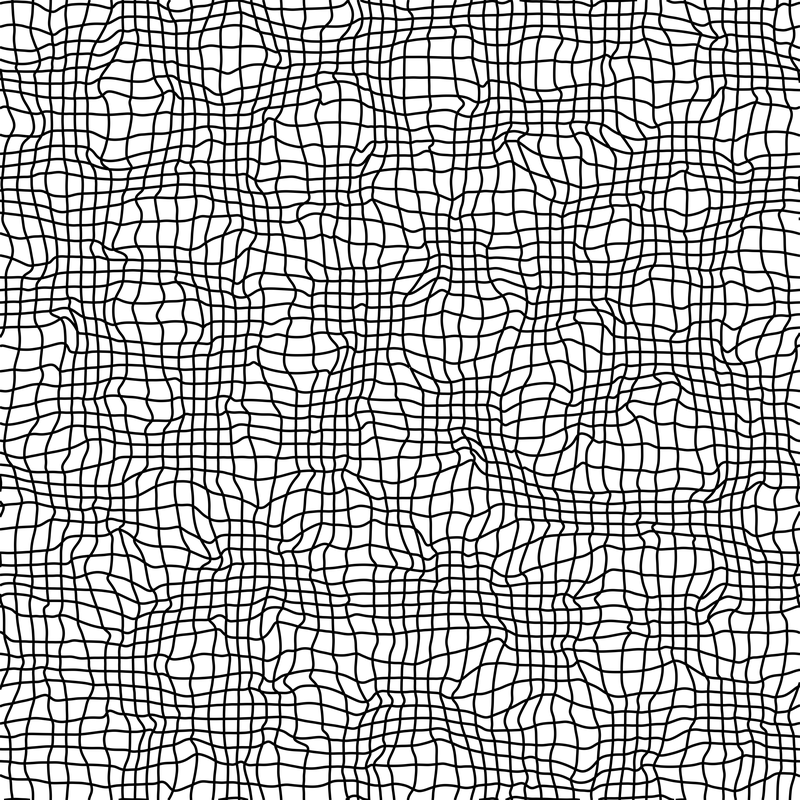}
        & \includegraphics[width=0.182\textwidth]{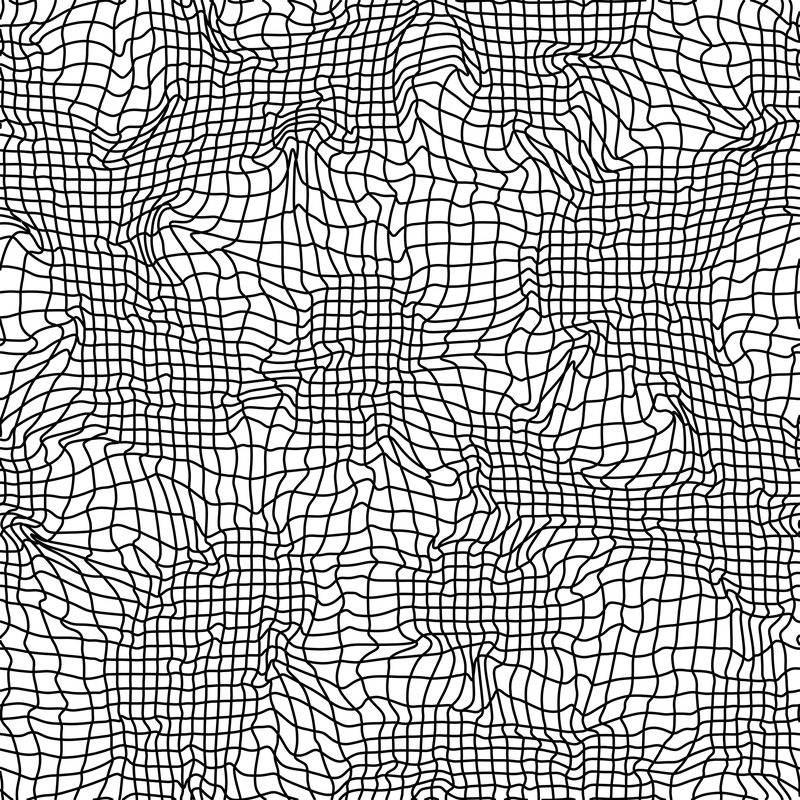}
        \\
        & {\small $t = 0$} & {\small $t = 1$} & {\small $t = 5$} & {\small $t = 50$} & {\small $t = 400$}
        \\
    \end{tabular*}
    \caption{Solution plots and deformations of the space induced by the lifted evolution over time in the random initial condition test for the Allen--Cahn and Cahn--Hilliard equations.}
    \label{fig:ac_ch_plots}
\end{figure}

\begin{figure}
    \centering
    \begin{subfigure}[c]{0.32\textwidth}
        \includegraphics[width=\textwidth]{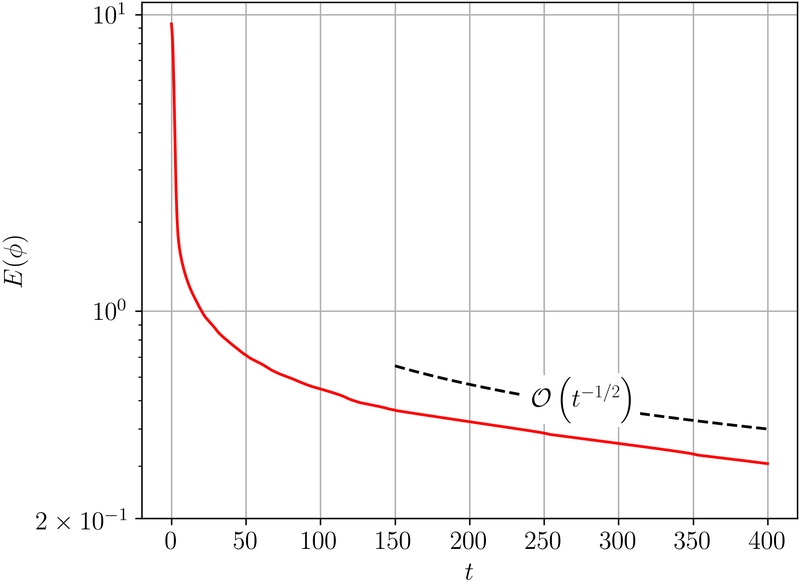}
    \end{subfigure}
    \hspace*{1em}
    \begin{subfigure}[c]{0.32\textwidth}
        \includegraphics[width=\textwidth]{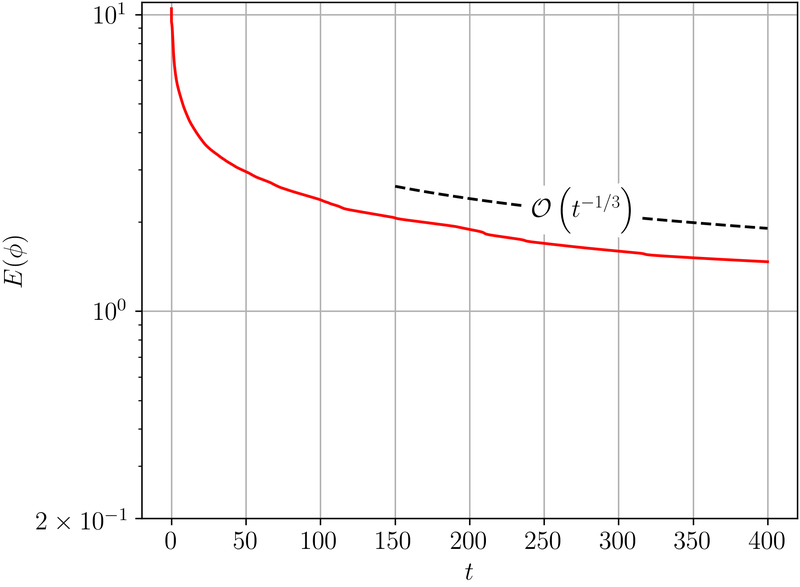}
    \end{subfigure}
    \caption{Energy $E(\solution)$ over time in the random initial condition tests in the Allen--Cahn (left) and Cahn--Hilliard (right) equations. The slight difference in initial energies is due to the different choice of $\epsilon$ for the two equations.}
    \label{fig:ac_ch_energies}
\end{figure}

For both equations, we plot the solutions, as well as the action of the lifted evolution on the space, at various times in Figure \ref{fig:ac_ch_plots}. In both systems we see the rapid initial phase separation, followed by either curvature-driven coarsening (for Allen--Cahn) or conservative coarsening (for Cahn--Hilliard). Notably, when comparing the solution plots to the induced spatial deformation at various times, we see the different roles played by the transport and compression components of the map. Indeed, the early fast diffusion does not induce a large distortion of the space; rather, it induces small, localized changes in compression, i.e. in the Jacobian. Conversely, the later coarsening induces a larger-scale mass transport, as the distortion of the space seems to align with the phase transitions. Finally, we plot the energy $E(\solution)$ over time for both equations in Figure \ref{fig:ac_ch_energies}. We see that the scheme is energy-stable, and that the energies decay at the correct asymptotic rates.

\subsection{2D Incompressible Navier--Stokes Equations}
\label{section:navier_stokes}

As a final numerical experiment we consider the incompressible Navier--Stokes equations
\begin{equation}
    \label{eq:navier_stokes}
    \partial_t u + u \cdot \nabla u = -\nabla p + \nu \Delta u, \quad \quad \operatorname{div} u = 0
\end{equation}
for a mean-zero vector field $u$ on $\mathbb{T}^2$, where $p$ is the pressure.
We solve the problem in vorticity form with $\solution = \operatorname{curl} u$, and recover the velocity field $u$ from $\solution$ by the Biot-Savart law
\begin{equation*}
    u = -\nabla^\perp \Delta^{-1} \solution =: A(\solution).
\end{equation*}
Then, taking the curl of \eqref{eq:navier_stokes}, we obtain the vorticity equation
\begin{equation*}
    \partial_t \solution + u \cdot \nabla \solution = \nu \Delta \solution.
\end{equation*}
This equation has constant mass $\dot{\mass} = 0$, and when $\nu = 0$, it reduces to the 2D incompressible Euler equations.

In the coordinate $\rho = \solution + c$, the lift reads
\begin{equation*}
    v = \liftingOperator_1 (\rho) = u - \nu \nabla \log \rho, \quad\quad u = -\nabla^\perp \Delta^{-1} (\rho - c).
\end{equation*}
That is, we lift the incompressible fluid velocity $u$ exactly, and the diffusion term lifts to a horizontal component identically to the heat equation. The evolution equation \eqref{eq:jacobian_L_plus_N_form_order2} of the log-Jacobian $\jacobian$ in the subinterval $[\tau_k, \tau_{k + 1}]$ reads
\begin{equation*}
    \partial_t \jacobian + (u - \nu \nabla \log \jsource) \cdot \nabla \jacobian = \nu \Delta \jacobian + \nu |\nabla \jacobian|^2 + \nu \Delta \log \jsource, \quad\quad u = -\nabla^\perp \Delta^{-1}(\jsource e^\jacobian - c).
\end{equation*}

\begin{figure}
    \centering
    \begin{minipage}[c]{0.02\textwidth}
        \rotatebox{90}{\makebox[0.18\textwidth][c]{$\nu = 10^{-4.5}$}}
    \end{minipage}
    \begin{subfigure}[c]{0.18\textwidth}
        \includegraphics[width=\textwidth]{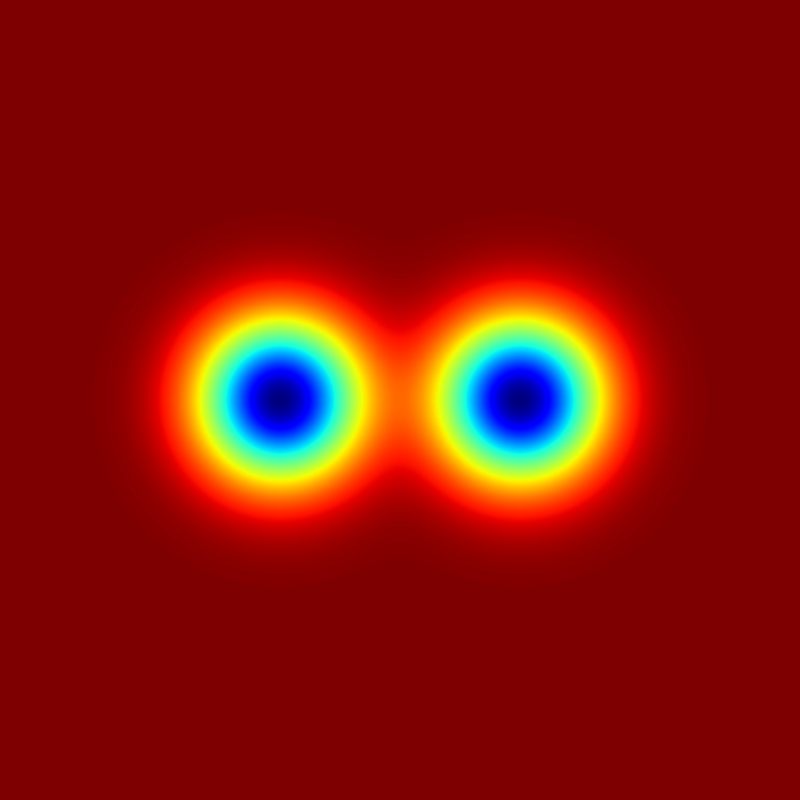}
    \end{subfigure}
    \begin{subfigure}[c]{0.18\textwidth}
        \includegraphics[width=\textwidth]{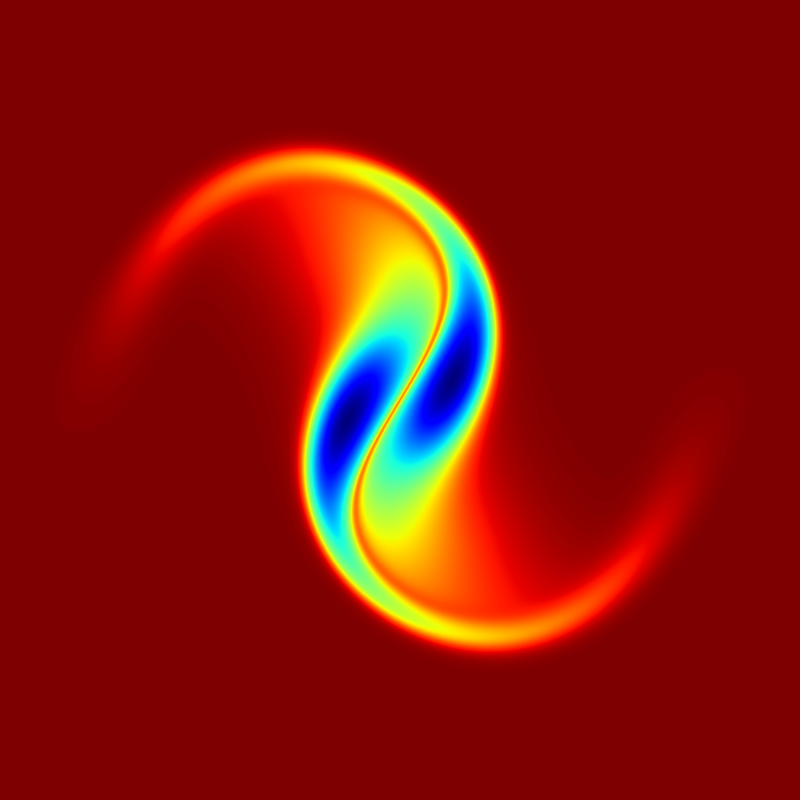}
    \end{subfigure}
    \begin{subfigure}[c]{0.18\textwidth}
        \includegraphics[width=\textwidth]{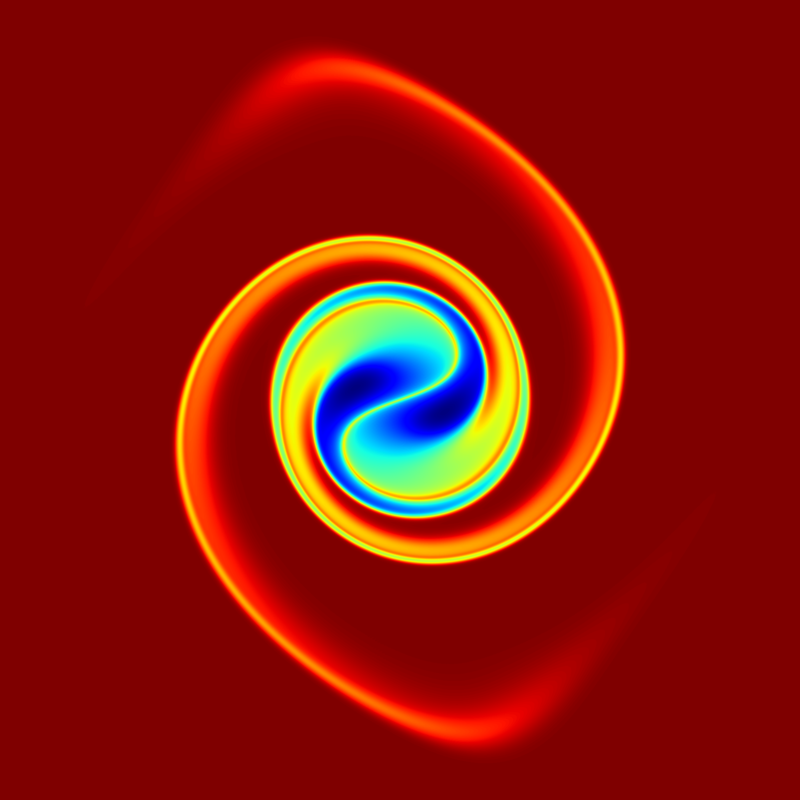}
    \end{subfigure}
    \begin{subfigure}[c]{0.18\textwidth}
        \includegraphics[width=\textwidth]{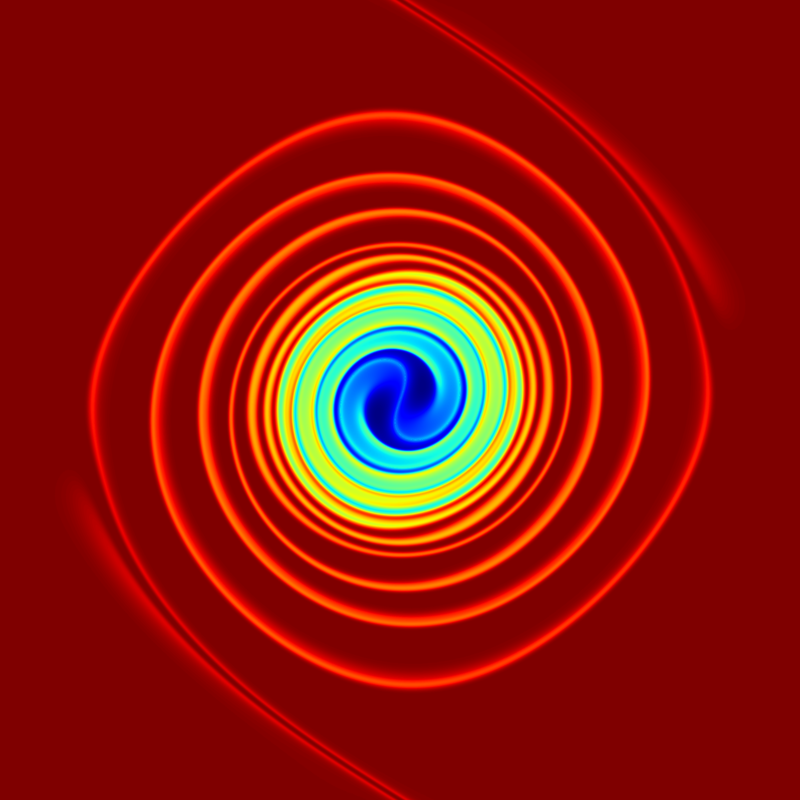}
    \end{subfigure}
    \begin{subfigure}[c]{0.18\textwidth}
        \includegraphics[width=\textwidth]{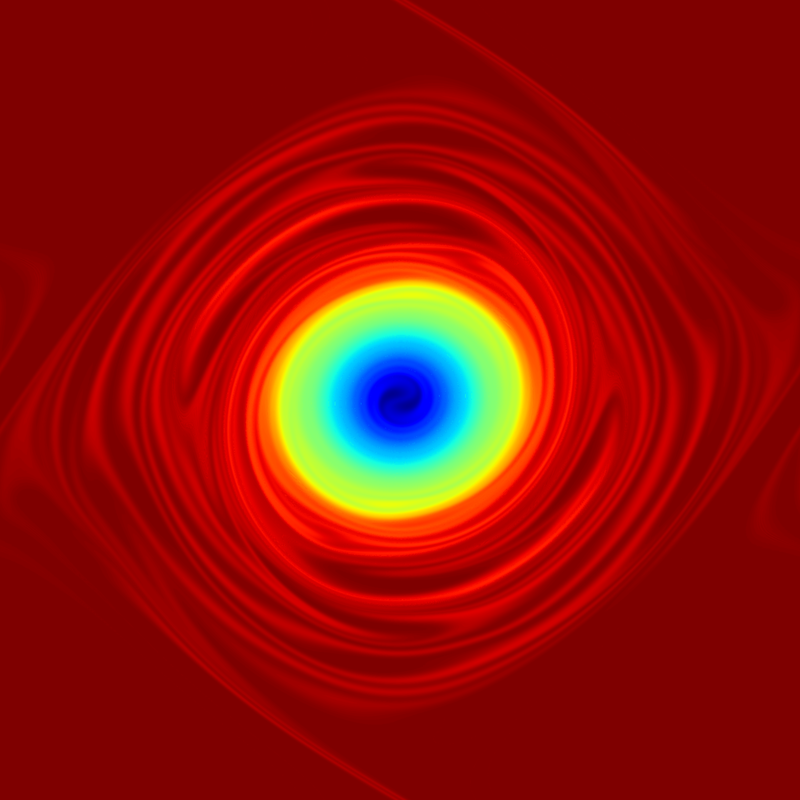}
    \end{subfigure}
    \\
    \begin{minipage}[c]{0.02\textwidth}
        \rotatebox{90}{\makebox[0.18\textwidth][c]{$\nu = 0$}}
    \end{minipage}
    \begin{subfigure}[c]{0.18\textwidth}
        \includegraphics[width=\textwidth]{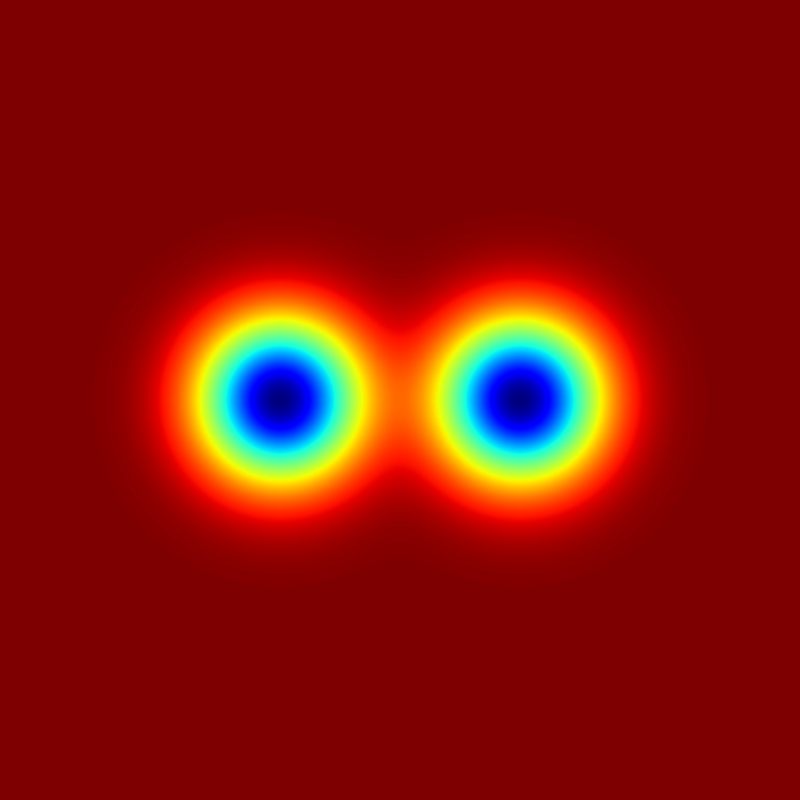}
    \end{subfigure}
    \begin{subfigure}[c]{0.18\textwidth}
        \includegraphics[width=\textwidth]{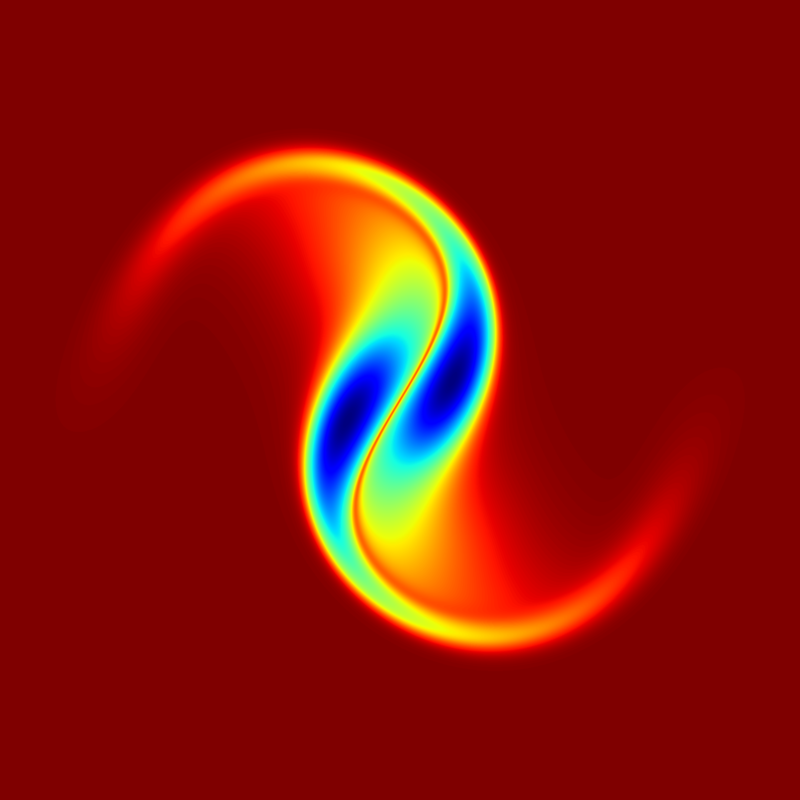}
    \end{subfigure}
    \begin{subfigure}[c]{0.18\textwidth}
        \includegraphics[width=\textwidth]{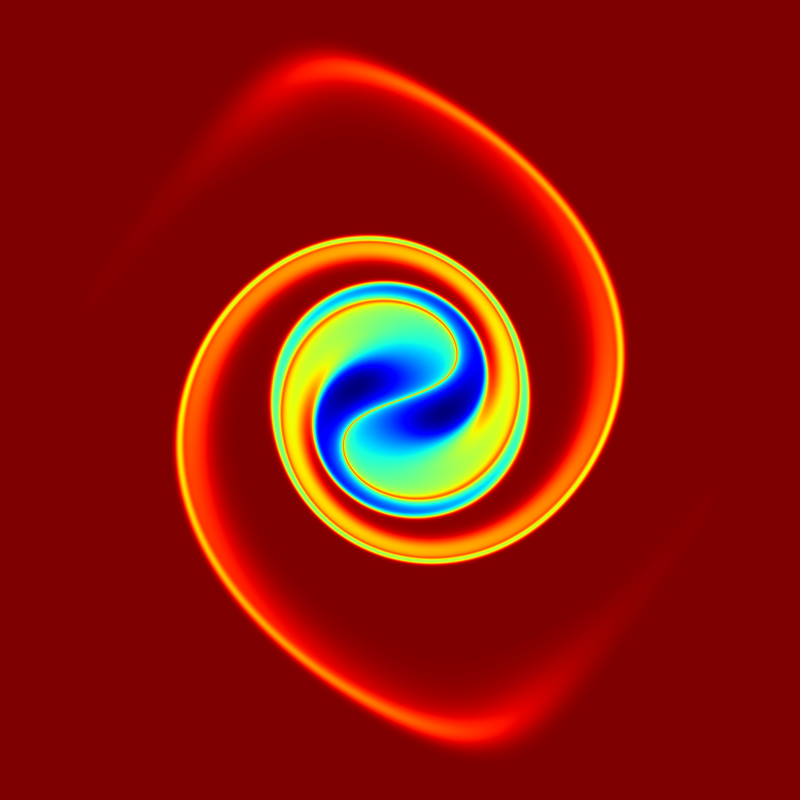}
    \end{subfigure}
    \begin{subfigure}[c]{0.18\textwidth}
        \includegraphics[width=\textwidth]{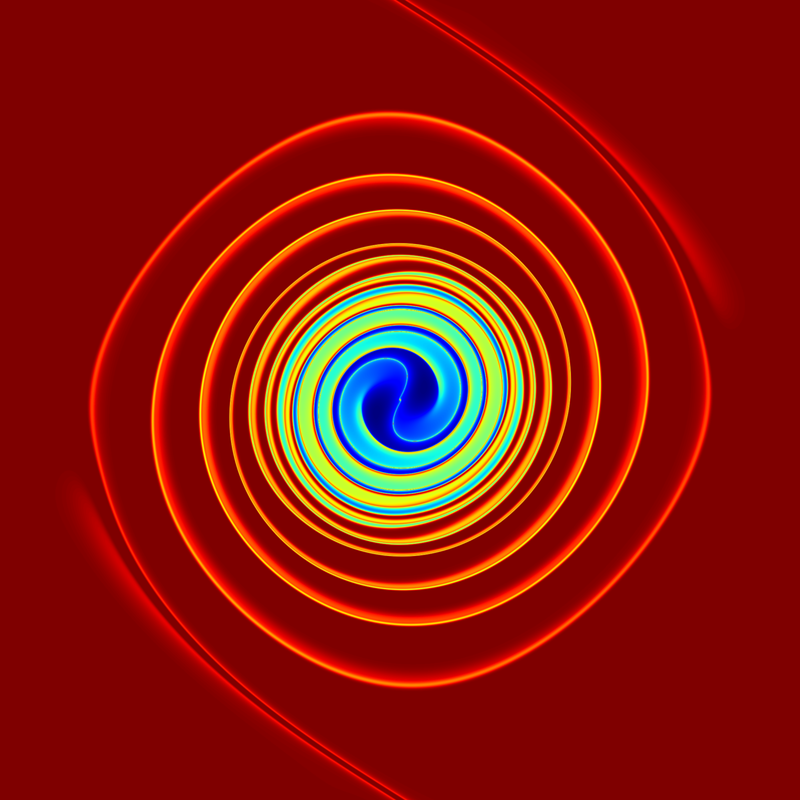}
    \end{subfigure}
    \begin{subfigure}[c]{0.18\textwidth}
        \includegraphics[width=\textwidth]{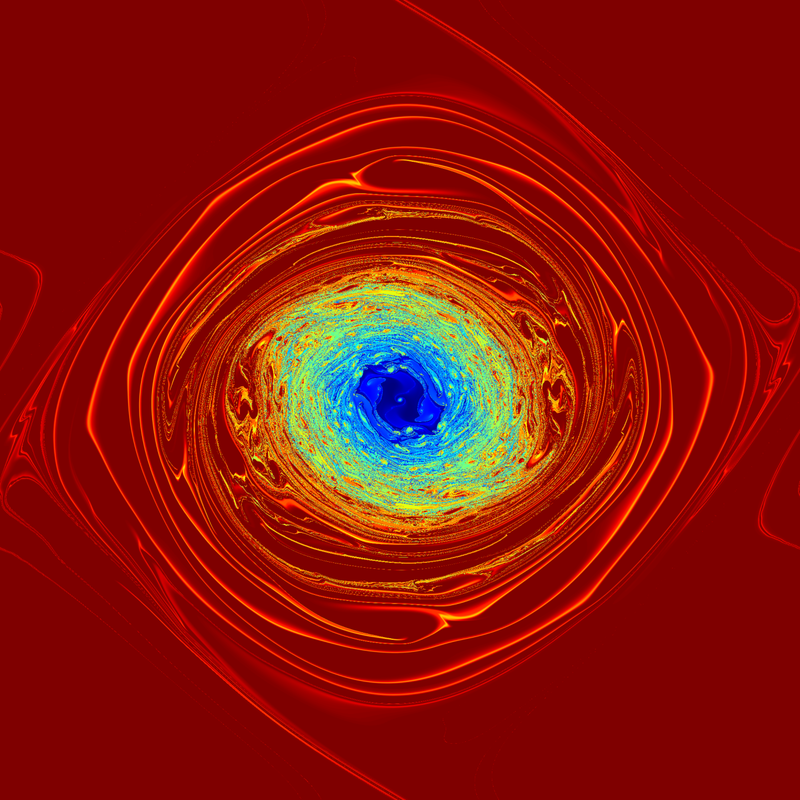}
    \end{subfigure}
    \\[0.3em]
    \hspace*{0.02\textwidth}
    \makebox[0.18\textwidth]{\parbox{0.18\textwidth}{\centering $t = 0$}}
    \makebox[0.18\textwidth]{\parbox{0.18\textwidth}{\centering $t = 2$}}
    \makebox[0.18\textwidth]{\parbox{0.18\textwidth}{\centering $t = 4$}}
    \makebox[0.18\textwidth]{\parbox{0.18\textwidth}{\centering $t = 10$}}
    \makebox[0.18\textwidth]{\parbox{0.18\textwidth}{\centering $t = 30$}}
    \caption{Vorticity plots over time in the Navier--Stokes vortex merger, for $\nu = 10^{-4.5}, 0$.}
    \label{fig:navier_stokes_summary_plots}
\end{figure}

As an experiment, we simulate the same vortex merger as was done for the CM method for 2D incompressible Euler \cite{yin_characteristic_2021}.
We take a periodized initial vorticity $\solution_0$ consisting of two negative Gaussian blobs, each with a standard deviation of $0.07 \times 2 \pi$, and separated by a distance of $0.3 \times 2\pi$.
The initial condition is shown at the left of Figure \ref{fig:navier_stokes_summary_plots}. We have $\min \solution_0 \approx -10$ and take $c = 100$. We take $\nu \in \{10^{-4.5}, 10^{-5.5}, 10^{-6}, 0\}$. These values were chosen in order to see the transition between the laminar and turbulent regimes. The case $\nu = 0$ corresponds to the incompressible Euler equations and, up to choices in discretization and integration schemes, this method is equivalent to the CM method for 2D incompressible Euler. For all values of $\nu$, we take $N = N_\epsilon = 512$, and $\Delta t = 1 / 512$. We use 8 maps per unit of time, and run the experiment to $t = 30$, for a total of 240 submaps.

We plot the vorticity at various times in Figure \ref{fig:navier_stokes_summary_plots} for $\nu = 10^{-4.5}$ and $\nu = 0$. The vortex blobs spin around each other. In the inviscid case $\nu = 0$, the blobs will never merge, and will continue spinning around each other, creating instabilities and finer scales as time goes on. In the dissipative regime $\nu > 0$ this is not the case, and the creation of fine scales is balanced by diffusion. However, we see that more instabilities are created and sustained as $\nu \rightarrow 0$, as expected.

The vorticity at the final time $t = 30$ is shown for all values of $\nu$ in the leftmost column of Figure \ref{fig:navier_stokes_zoom}. Further, since the submaps can be evaluated everywhere in the domain by interpolation, we are able to zoom into the solution, as in previous iterations of the CM method \cite{yin_characteristic_2021}. This is shown in the rest of Figure \ref{fig:navier_stokes_zoom}, where we gradually zoom in around the point $(0.6, 0.5) \in [0, 1]^2$, each time by a factor of $1 / 2$. Here we again see the increase in instability and the finer scales generated by the flow as $\nu \rightarrow 0$. In particular, we see that we retain the resolution benefits of the CM method.

\newcommand{\navierstokeszoomimagewidthpercent}{0.18}
\begin{figure}
    \centering
    \makebox[1em]{\rotatebox{90}{\parbox{\navierstokeszoomimagewidthpercent\textwidth}{\centering $\nu = 10^{-4.5}$}}}
    \begin{subfigure}{\navierstokeszoomimagewidthpercent\textwidth}
        \includegraphics[width=\textwidth]{navier_stokes_v2/nu_e-4.5/30.png}
    \end{subfigure}
    \begin{subfigure}{\navierstokeszoomimagewidthpercent\textwidth}
        \includegraphics[width=\textwidth]{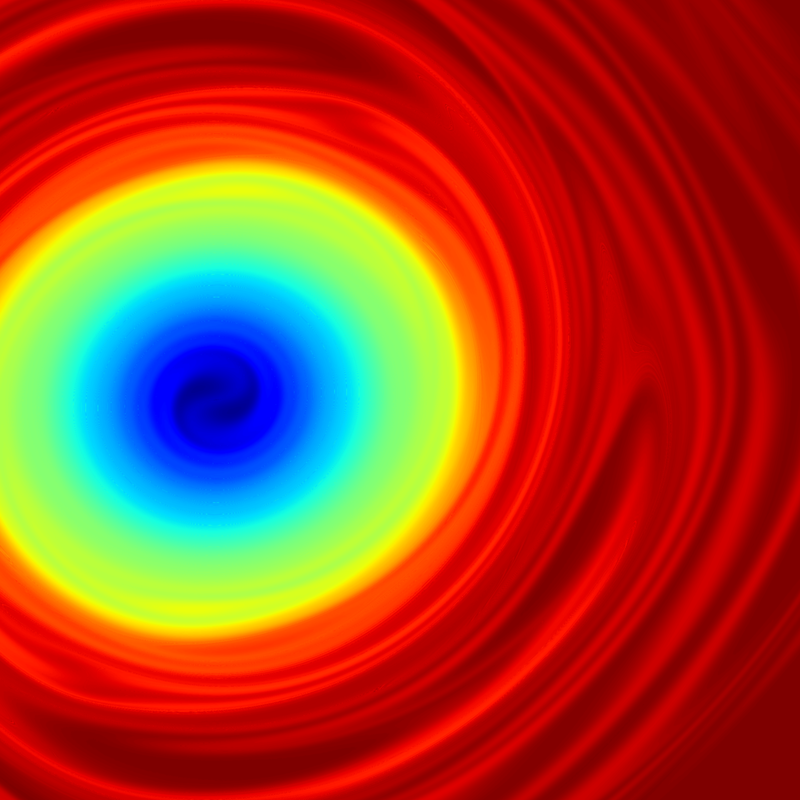}
    \end{subfigure}
    \begin{subfigure}{\navierstokeszoomimagewidthpercent\textwidth}
        \includegraphics[width=\textwidth]{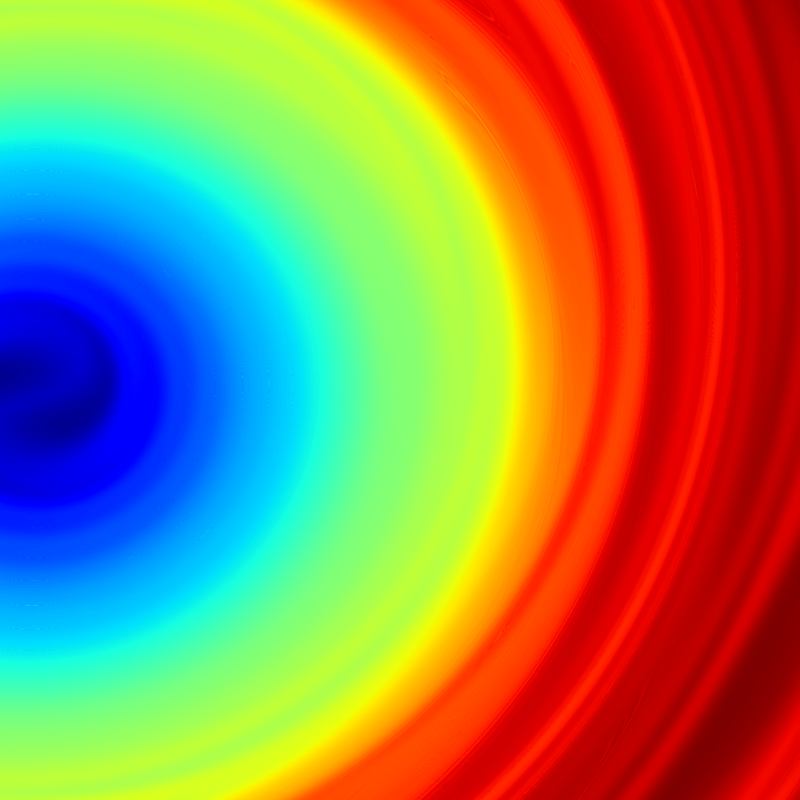}
    \end{subfigure}
    \begin{subfigure}{\navierstokeszoomimagewidthpercent\textwidth}
        \includegraphics[width=\textwidth]{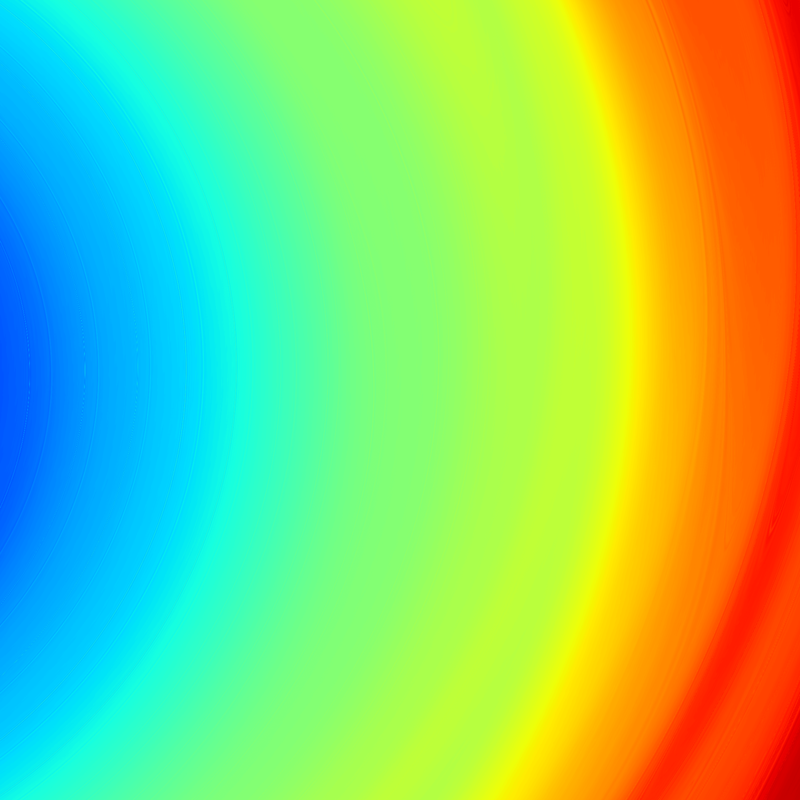}
    \end{subfigure}
    \begin{subfigure}{\navierstokeszoomimagewidthpercent\textwidth}
        \includegraphics[width=\textwidth]{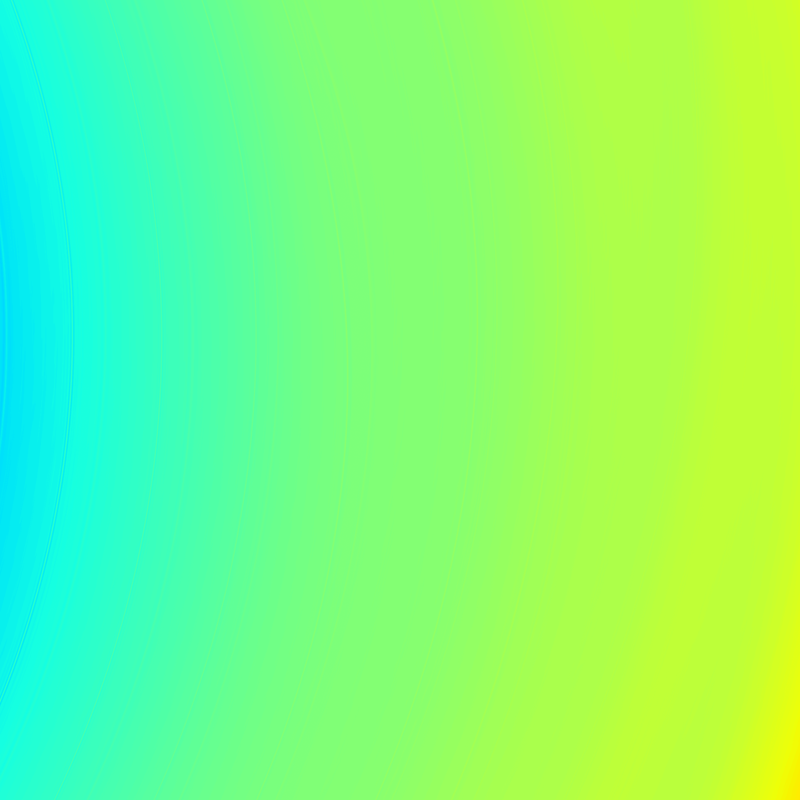}
    \end{subfigure}
    \\
    \makebox[1em]{\rotatebox{90}{\parbox{\navierstokeszoomimagewidthpercent\textwidth}{\centering $\nu = 10^{-5.5}$}}}
    \begin{subfigure}{\navierstokeszoomimagewidthpercent\textwidth}
        \includegraphics[width=\textwidth]{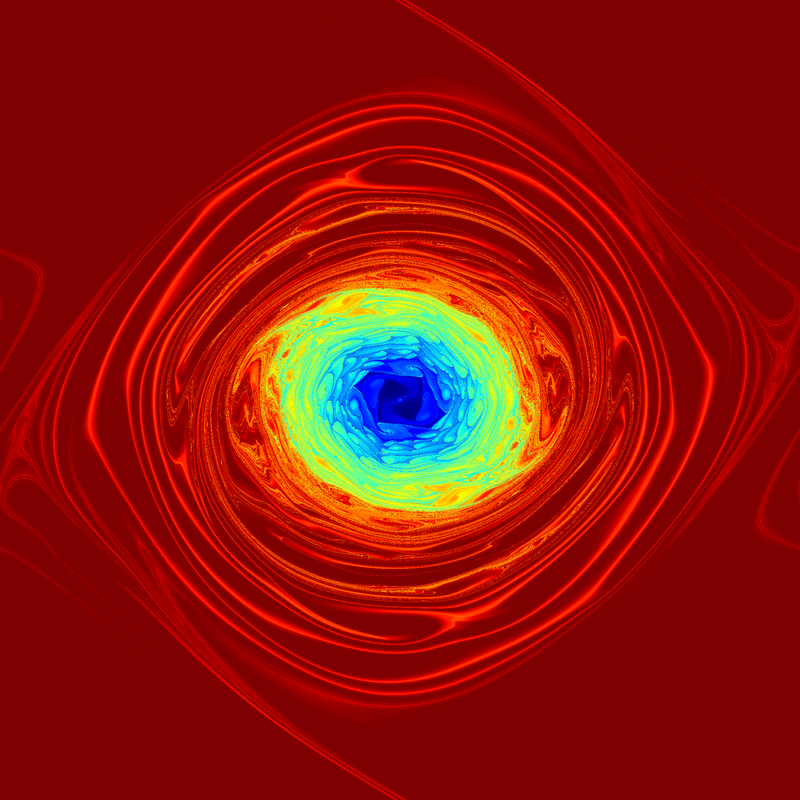}
    \end{subfigure}
    \begin{subfigure}{\navierstokeszoomimagewidthpercent\textwidth}
        \includegraphics[width=\textwidth]{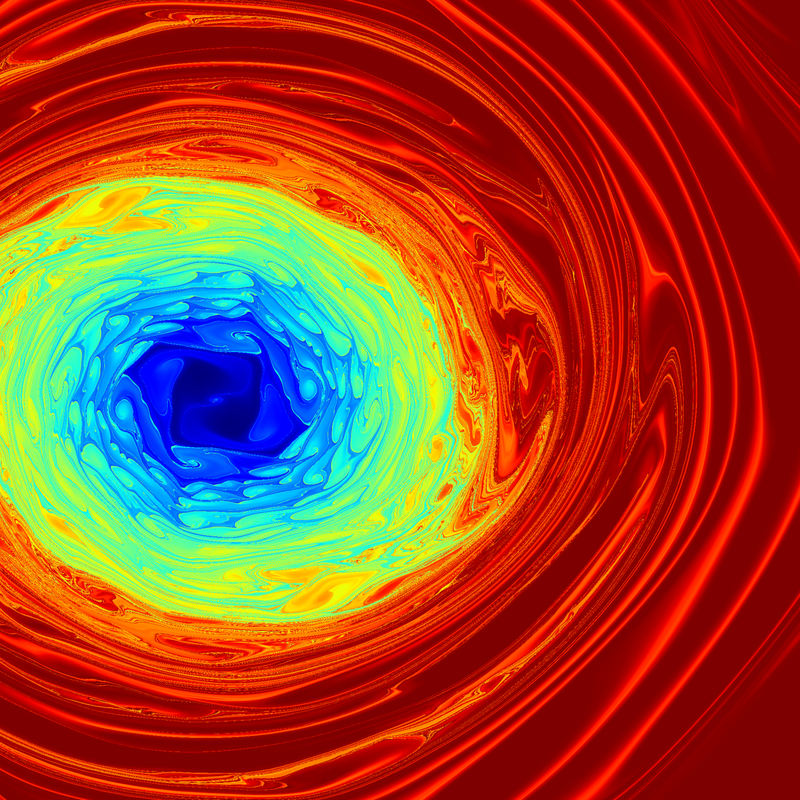}
    \end{subfigure}
    \begin{subfigure}{\navierstokeszoomimagewidthpercent\textwidth}
        \includegraphics[width=\textwidth]{navier_stokes_v2/zooms_3/nu_e-5.5_z2.png}
    \end{subfigure}
    \begin{subfigure}{\navierstokeszoomimagewidthpercent\textwidth}
        \includegraphics[width=\textwidth]{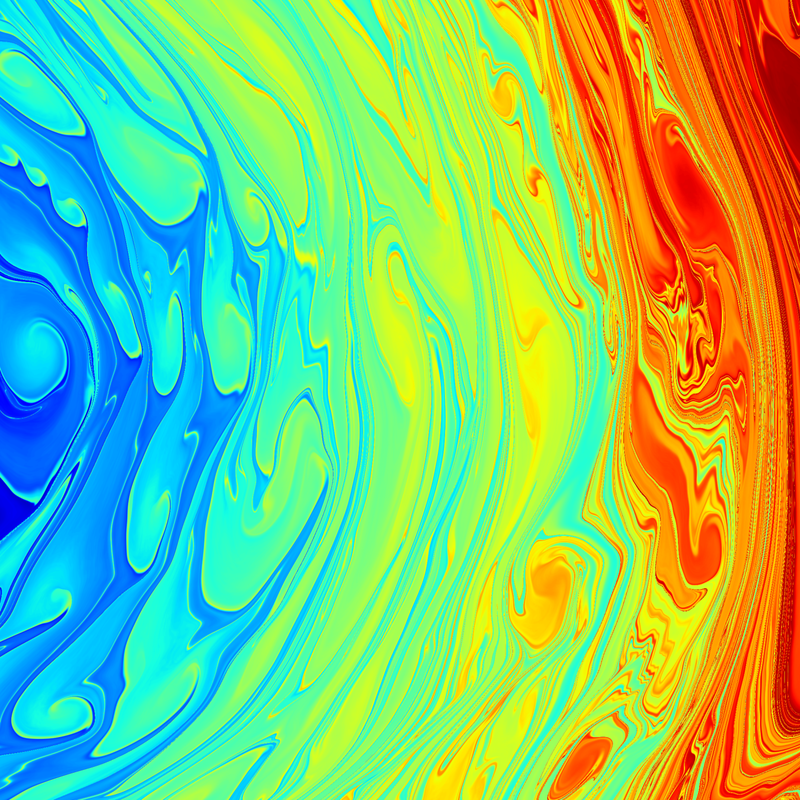}
    \end{subfigure}
    \begin{subfigure}{\navierstokeszoomimagewidthpercent\textwidth}
        \includegraphics[width=\textwidth]{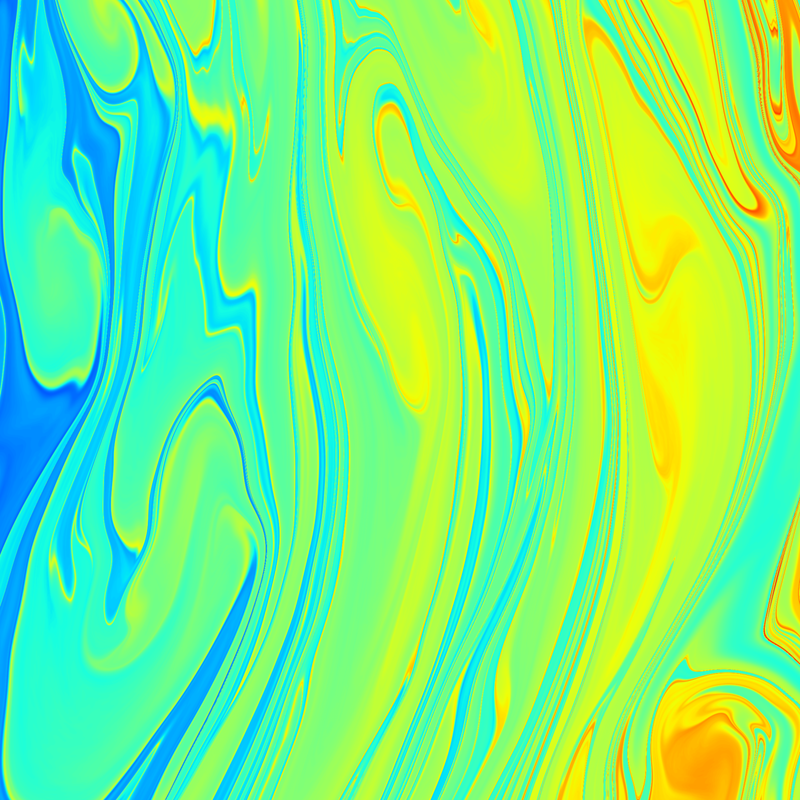}
    \end{subfigure}
    \\
    \makebox[1em]{\rotatebox{90}{\parbox{\navierstokeszoomimagewidthpercent\textwidth}{\centering $\nu = 10^{-6}$}}}
    \begin{subfigure}{\navierstokeszoomimagewidthpercent\textwidth}
        \includegraphics[width=\textwidth]{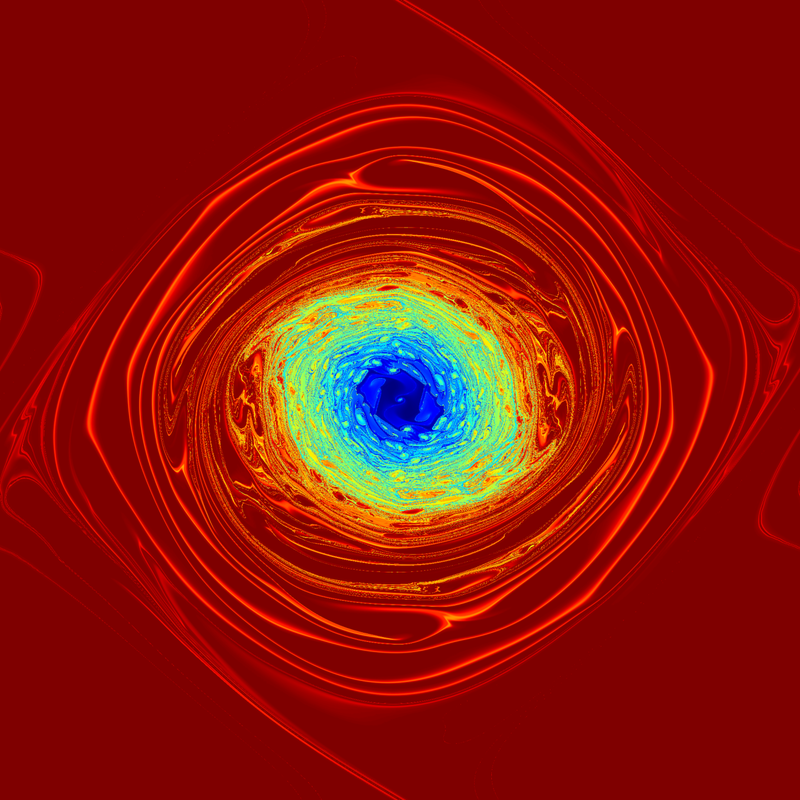}
    \end{subfigure}
    \begin{subfigure}{\navierstokeszoomimagewidthpercent\textwidth}
        \includegraphics[width=\textwidth]{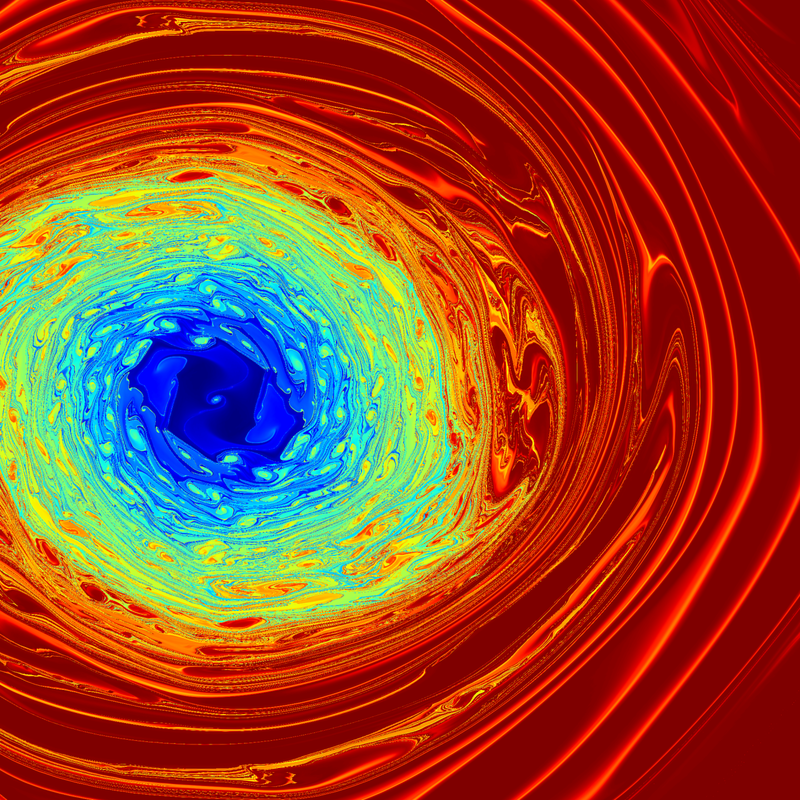}
    \end{subfigure}
    \begin{subfigure}{\navierstokeszoomimagewidthpercent\textwidth}
        \includegraphics[width=\textwidth]{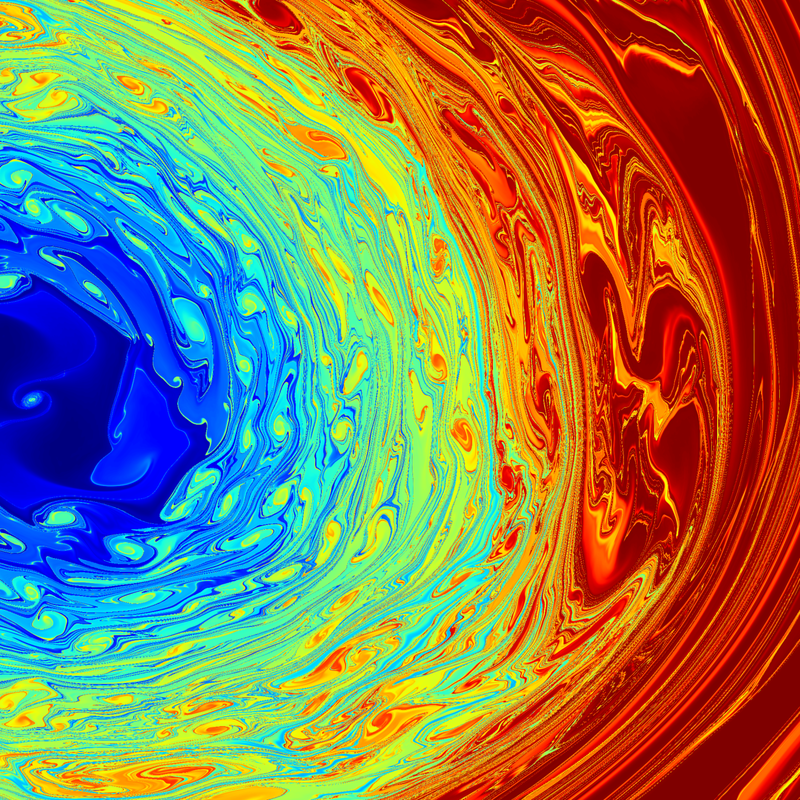}
    \end{subfigure}
    \begin{subfigure}{\navierstokeszoomimagewidthpercent\textwidth}
        \includegraphics[width=\textwidth]{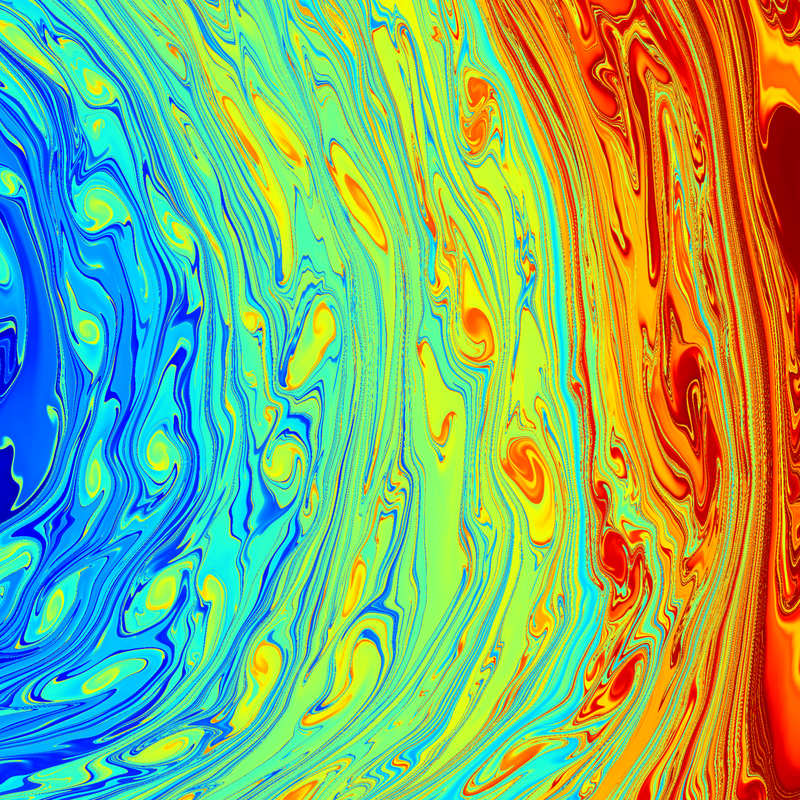}
    \end{subfigure}
    \begin{subfigure}{\navierstokeszoomimagewidthpercent\textwidth}
        \includegraphics[width=\textwidth]{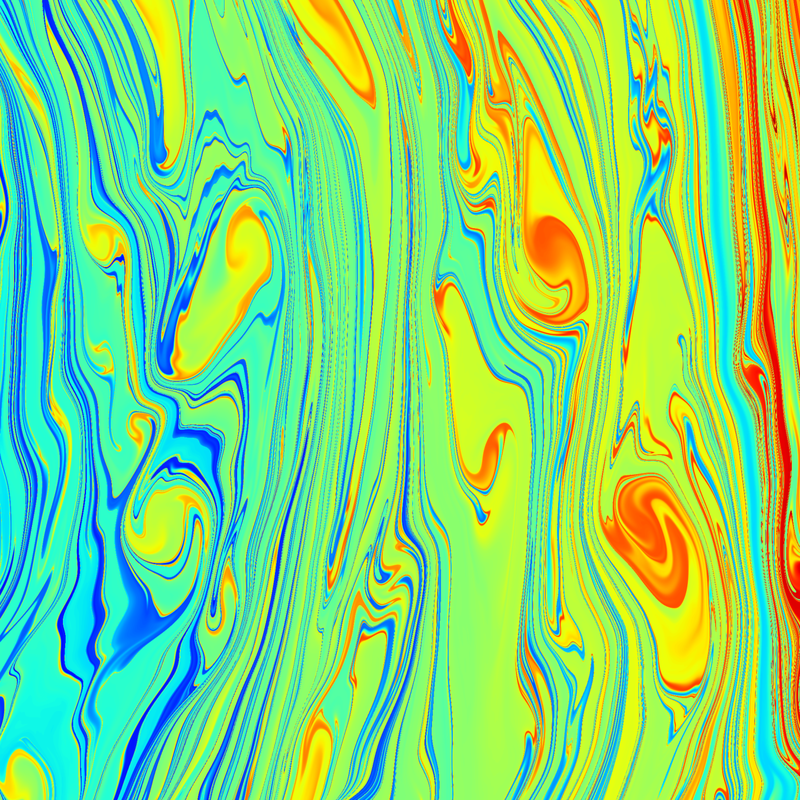}
    \end{subfigure}
    \\
    \makebox[1em]{\rotatebox{90}{\parbox{\navierstokeszoomimagewidthpercent\textwidth}{\centering $\nu = 0$}}}
    \begin{subfigure}{\navierstokeszoomimagewidthpercent\textwidth}
        \includegraphics[width=\textwidth]{navier_stokes_v2/nu_0/30.png}
    \end{subfigure}
    \begin{subfigure}{\navierstokeszoomimagewidthpercent\textwidth}
        \includegraphics[width=\textwidth]{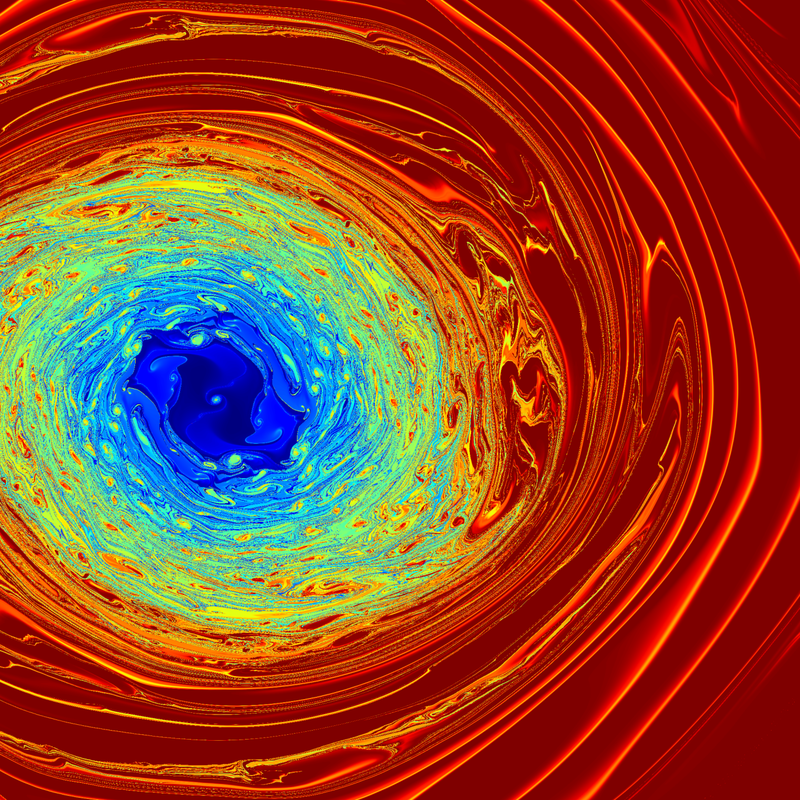}
    \end{subfigure}
    \begin{subfigure}{\navierstokeszoomimagewidthpercent\textwidth}
        \includegraphics[width=\textwidth]{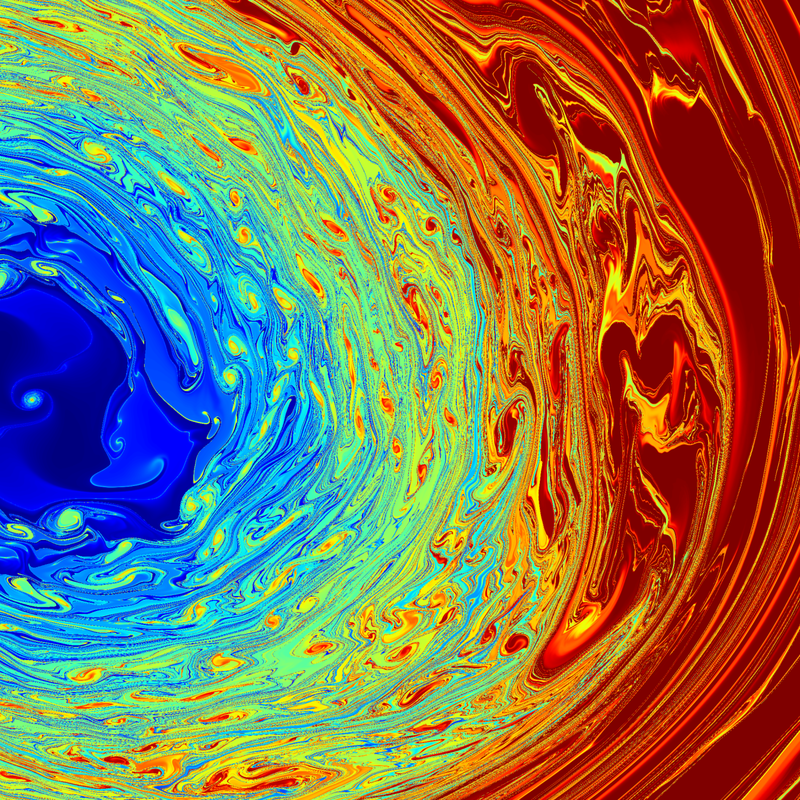}
    \end{subfigure}
    \begin{subfigure}{\navierstokeszoomimagewidthpercent\textwidth}
        \includegraphics[width=\textwidth]{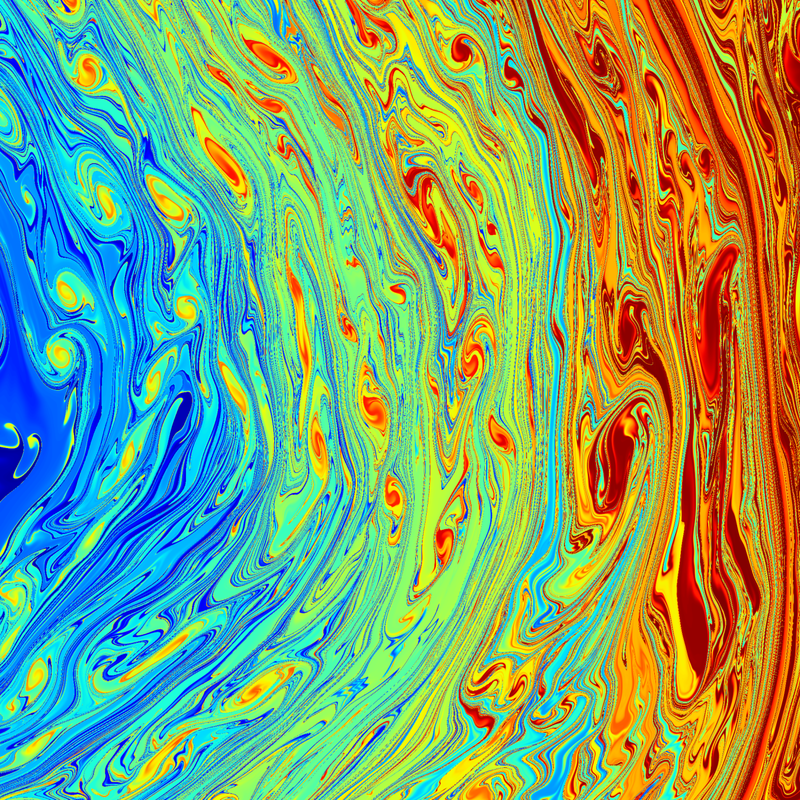}
    \end{subfigure}
    \begin{subfigure}{\navierstokeszoomimagewidthpercent\textwidth}
        \includegraphics[width=\textwidth]{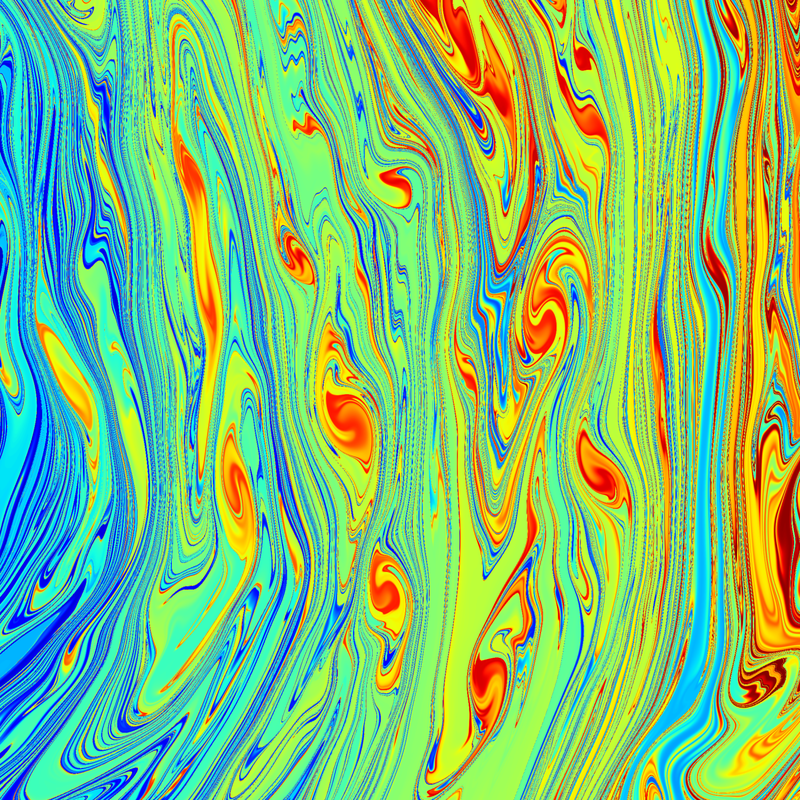}
    \end{subfigure}
    \\
    \makebox[1em]{\rotatebox{90}{\parbox{0\textwidth}{\centering}}}
    \makebox[\navierstokeszoomimagewidthpercent\textwidth]{\parbox{\navierstokeszoomimagewidthpercent\textwidth}{\centering 1x zoom}}
    \makebox[\navierstokeszoomimagewidthpercent\textwidth]{\parbox{\navierstokeszoomimagewidthpercent\textwidth}{\centering 2x zoom}}
    \makebox[\navierstokeszoomimagewidthpercent\textwidth]{\parbox{\navierstokeszoomimagewidthpercent\textwidth}{\centering 4x zoom}}
    \makebox[\navierstokeszoomimagewidthpercent\textwidth]{\parbox{\navierstokeszoomimagewidthpercent\textwidth}{\centering 8x zoom}}
    \makebox[\navierstokeszoomimagewidthpercent\textwidth]{\parbox{\navierstokeszoomimagewidthpercent\textwidth}{\centering 16x zoom}}
    \caption{Zooms into the vorticity plots in the Navier--Stokes vortex merger at $t = 30$ across all values of $\nu$.}
    \label{fig:navier_stokes_zoom}
\end{figure}

Finally, in Figure \ref{fig:navier_stokes_components}, we plot the transport component $\rho_0 \circ \map$, the compression component $e^\jacobian$, and the solution $\rho = (\rho_0 \circ \map) e^\jacobian$ at $t = 30$ for the positive values of $\nu$. Here we see the contribution of the Jacobian shrink as $\nu \rightarrow 0$, as expected. We also see, from the pure-transport component $\rho_0 \circ \map$, that the flow itself is far less violent for larger $\nu$. This illustrates how the method fundamentally differs from existing methods which track a compression component along particle trajectories. In our formulation, the diffusion term lifts to a gradient component in the map velocity, altering the computed flow itself. The compression component $e^\jacobian$, being also advected by the flow through the semidirect product action, then naturally lives in a growing approximation space as well; this can be seen from the fine scales in the Jacobian in the middle panel of Figure \ref{fig:navier_stokes_components}. These results indicate a state of the art ability of the method to capture turbulent fluid dynamics at arbitrary Reynolds number without relying on an increasingly fine scale spatial resolution. The nonlinear decomposition of the solution into transport and compression is a central enabling factor for these results.  

\newcommand{\navierstokescomponentsimagewidthpercent}{0.28}
\begin{figure}[H]
    \centering
    \hspace*{0.02\textwidth}
    \begin{subfigure}[c]{\navierstokescomponentsimagewidthpercent\textwidth}
        \includegraphics[width=\textwidth]{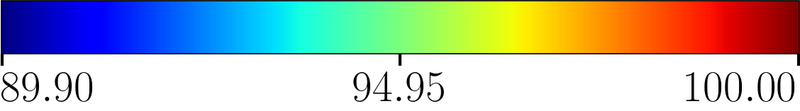}
    \end{subfigure}
    \hspace*{0.045\textwidth}
    \begin{subfigure}[c]{\navierstokescomponentsimagewidthpercent\textwidth}
        \includegraphics[width=\textwidth]{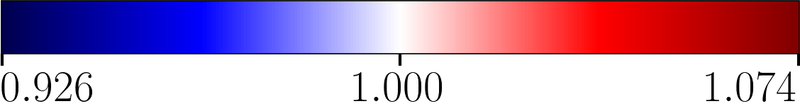}
    \end{subfigure}
    \hspace*{0.045\textwidth}
    \begin{subfigure}[c]{\navierstokescomponentsimagewidthpercent\textwidth}
        \includegraphics[width=\textwidth]{navier_stokes_v2/vorticity_colorbar.png}
    \end{subfigure}
    \\
    \begin{minipage}[c]{0.02\textwidth}
        \rotatebox{90}{\makebox[\navierstokescomponentsimagewidthpercent\textwidth][c]{$\nu = 10^{-4.5}$}}
    \end{minipage}
    \begin{subfigure}[c]{\navierstokescomponentsimagewidthpercent\textwidth}
        \includegraphics[width=\textwidth]{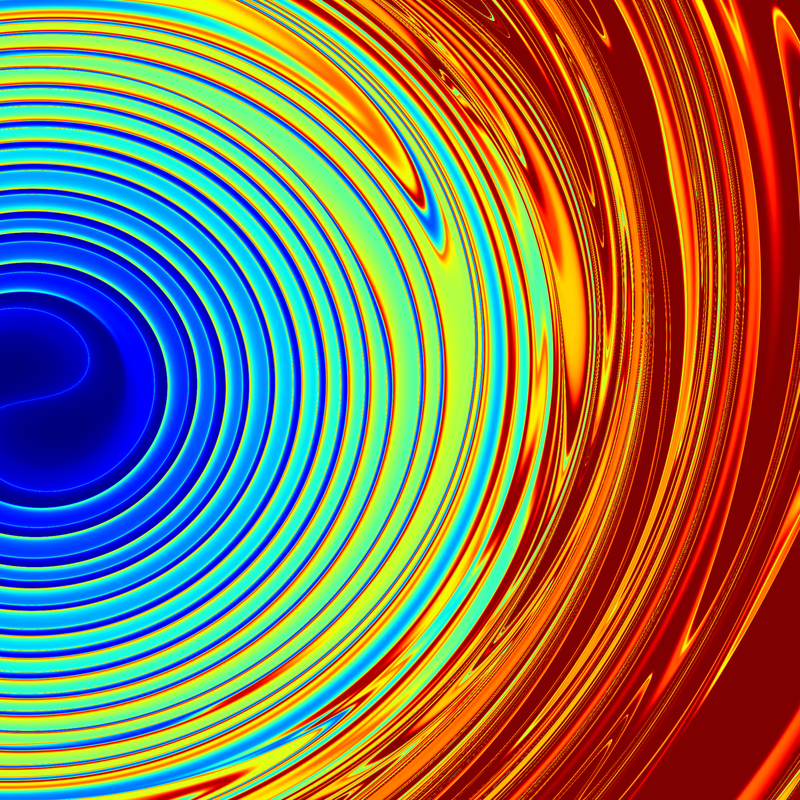}
    \end{subfigure}
    \begin{minipage}[c]{0.045\textwidth}
        \centering
        \Huge $\times$
    \end{minipage}
    \begin{subfigure}[c]{\navierstokescomponentsimagewidthpercent\textwidth}
        \includegraphics[width=\textwidth]{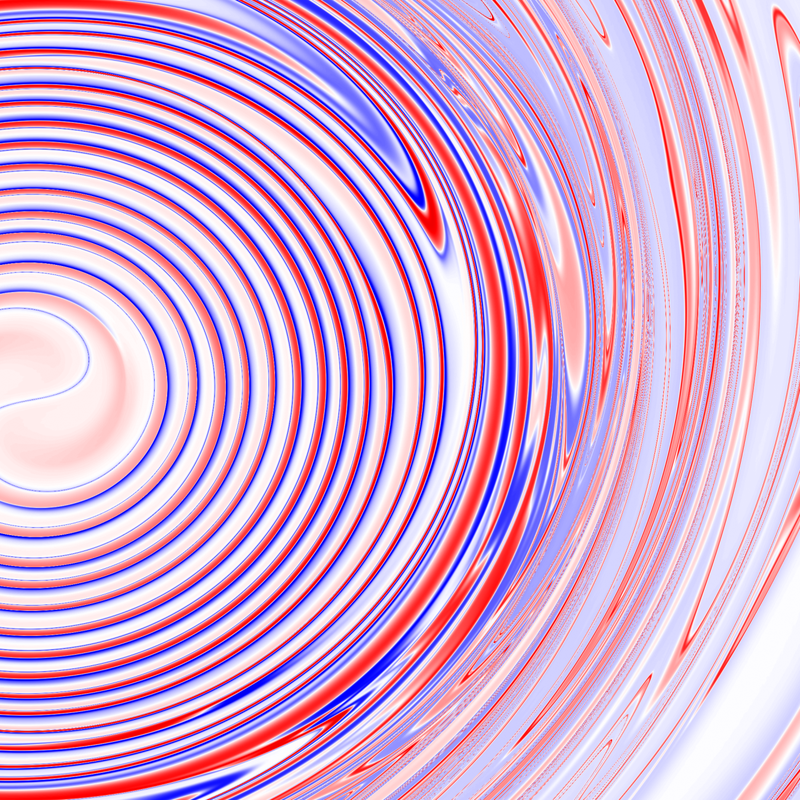}
    \end{subfigure}
    \begin{minipage}[c]{0.045\textwidth}
        \centering
        \Huge $=$
    \end{minipage}
    \begin{subfigure}[c]{\navierstokescomponentsimagewidthpercent\textwidth}
        \includegraphics[width=\textwidth]{navier_stokes_v2/zooms_3/nu_e-4.5_z2.png}
    \end{subfigure}
    \\
    \begin{minipage}[c]{0.02\textwidth}
        \rotatebox{90}{\makebox[\navierstokescomponentsimagewidthpercent\textwidth][c]{$\nu = 10^{-5.5}$}}
    \end{minipage}
    \begin{subfigure}[c]{\navierstokescomponentsimagewidthpercent\textwidth}
        \includegraphics[width=\textwidth]{navier_stokes_v2/components/nu_e-5.5_z2_transport.png}
    \end{subfigure}
    \begin{minipage}[c]{0.045\textwidth}
        \centering
        \Huge $\times$
    \end{minipage}
    \begin{subfigure}[c]{\navierstokescomponentsimagewidthpercent\textwidth}
        \includegraphics[width=\textwidth]{navier_stokes_v2/components/nu_e-5.5_z2_jacobian_seismic.png}
    \end{subfigure}
    \begin{minipage}[c]{0.045\textwidth}
        \centering
        \Huge $=$
    \end{minipage}
    \begin{subfigure}[c]{\navierstokescomponentsimagewidthpercent\textwidth}
        \includegraphics[width=\textwidth]{navier_stokes_v2/zooms_3/nu_e-5.5_z2.png}
    \end{subfigure}
    \\
    \begin{minipage}[c]{0.02\textwidth}
        \rotatebox{90}{\makebox[\navierstokescomponentsimagewidthpercent\textwidth][c]{$\nu = 10^{-6}$}}
    \end{minipage}
    \begin{subfigure}[c]{\navierstokescomponentsimagewidthpercent\textwidth}
        \includegraphics[width=\textwidth]{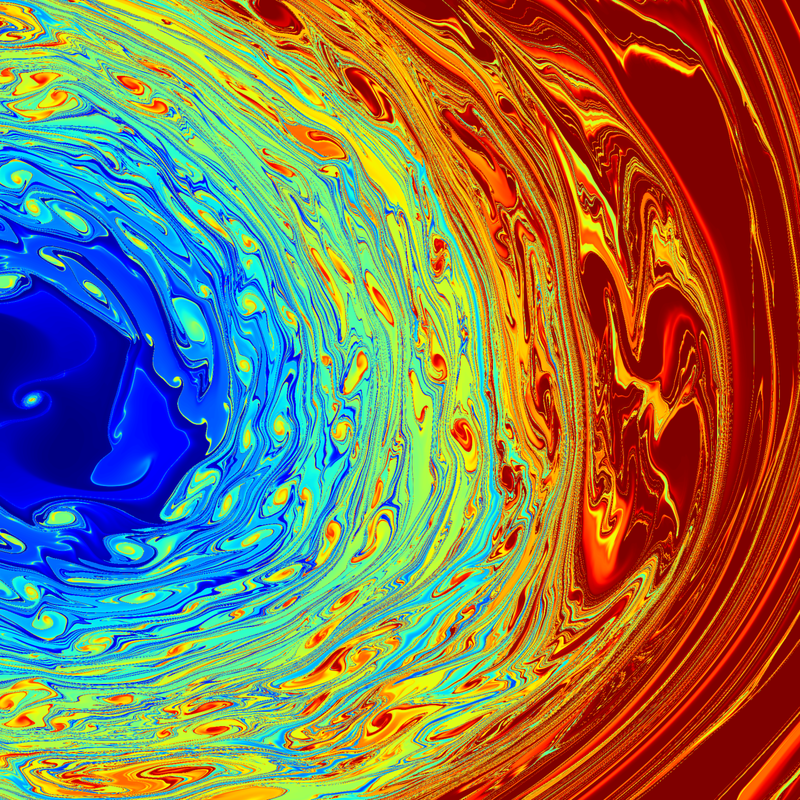}
    \end{subfigure}
    \begin{minipage}[c]{0.045\textwidth}
        \centering
        \Huge $\times$
    \end{minipage}
    \begin{subfigure}[c]{\navierstokescomponentsimagewidthpercent\textwidth}
        \includegraphics[width=\textwidth]{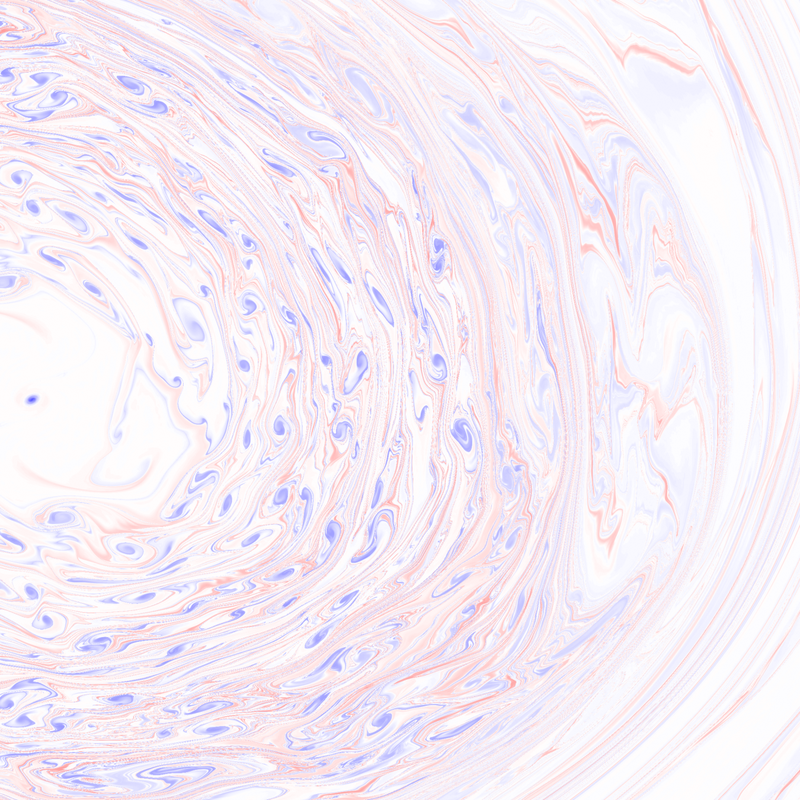}
    \end{subfigure}
    \begin{minipage}[c]{0.045\textwidth}
        \centering
        \Huge $=$
    \end{minipage}
    \begin{subfigure}[c]{\navierstokescomponentsimagewidthpercent\textwidth}
        \includegraphics[width=\textwidth]{navier_stokes_v2/zooms_3/nu_e-6_z2.png}
    \end{subfigure}
    \\[0.3em]
    \hspace*{0.02\textwidth}
    \makebox[\navierstokescomponentsimagewidthpercent\textwidth]{\parbox{\navierstokescomponentsimagewidthpercent\textwidth}{\centering $\rho_0 \circ \map$}}
    \hspace*{0.045\textwidth}
    \makebox[\navierstokescomponentsimagewidthpercent\textwidth]{\parbox{\navierstokescomponentsimagewidthpercent\textwidth}{\centering $\jacobianOperator(\map)$}}
    \hspace*{0.045\textwidth}
    \makebox[\navierstokescomponentsimagewidthpercent\textwidth]{\parbox{\navierstokescomponentsimagewidthpercent\textwidth}{\centering $\rho$}}
    \caption{Decomposition of the solution $\rho$ in the Navier--Stokes vortex merger at $t = 30$ into the different map components: the transport component $\rho_0 \circ \map$ and the compression component $\jacobianOperator(\map)$.}
    \label{fig:navier_stokes_components}
\end{figure}

\section{Conclusion}

In this work, we have presented an extension of the CM method to the solution of general scalar evolution equations. We showed that an evolution equation can be rewritten into the form of a continuity equation, after accounting for a varying total mass. This defines a lift of the dynamics to an extended space of diffeomorphisms, giving access to a semigroup decomposition, and retaining many of the resolution benefits of the CM method for transport equations. We defined a specific, numerically-motivated, construction of this lifting operation and provided a geometric interpretation within the theory of unbalanced optimal transport.\par

Based on this lifted evolution equation in the extended diffeomorphism group, we devised a semidirect product evolution algorithm to explicitly evolve the log-Jacobian of the map. This construction mitigated for a primary source of instability arising from computing this operation analytically. We constructed a numerical implementation on the torus, which is designed for the stable, high-order integration for a wide class of PDEs, and is directly assembled from standard numerical techniques. We analyzed and validated the order of the method on a range of nonlinear dissipative evolution equations. The experiments demonstrated that the scheme recovers the correct qualitative behaviour of the solution across different regimes. For the Navier--Stokes equations specifically, we showed that the extension retains the resolution benefits of the CM method and yields a uniquely promising ability to simulate flows across all Reynolds numbers.

There are many avenues to extend upon and apply the numerical method presented in this work. We are currently working on generalizing these approaches to systems of equations using recent developments in the theory of vector optimal transport \cite{khesin2025universalvectormatrixoptimal}. Another generalization would consider the case of domains with boundary and the implementation of boundary conditions. On the numerical side, the possibility of constructing an unconditionally energy-stable integration scheme for the Jacobian is warranted. Finally, as mentioned previously, a future work will involve a deeper study of the Navier--Stokes equations, with a focus on the quality of the growing approximation spaces and the inviscid limit $\nu \rightarrow 0$.

\section*{Acknowledgments}

W.H.B. is grateful to Nathaniel and Ciara Giroux for many helpful discussions.
S.T. gratefully acknowledges that this research was supported in part by Pacific Institute for the Mathematical Sciences through its Collaborative Research Group on Structure-Preserving Discretizations and their Applications, the Natural Sciences and Engineering Research Council of Canada through its Alliance Advantage program, and Environment and Climate Change Canada through its Grants \& Contributions program.
J.C.N. acknowledges financial support from the NSERC Discovery Grant program of Canada and from the French ANR under projects CM2E and CROKE.   

\section*{\normalsize Statement of Author Contributions}

\textbf{William Holman--Bissegger}: Conceptualization, methodology, formal analysis, software, validation, visualization, writing -- original draft, writing -- review \& editing.
\textbf{Seth Taylor}: Conceptualization, methodology, formal analysis, writing -- original draft, writing -- review \& editing.
\textbf{Jean--Christophe Nave}: Conceptualization, methodology, supervision, writing -- review \& editing, project administration, funding acquisition.

\appendix

\section{Proof of Theorem \ref{thm:error_estimate}}
\label{appendix:proof_of_error_theorem}

In this appendix, we establish the proof of the error estimate, Theorem \ref{thm:error_estimate}. We begin by recalling some notation from our outline of the proof in Section \ref{section:error_estimate}. We denote the target error as $\delta := \mathcal{O}(\Delta t^4 + \Delta x^6 + \Delta t \Delta x^4 + \epsilon)$.
The log-Jacobian data term $\jsource = \mass \cdot (\rho_0 \circ \map)$ satisfies $\partial_t \jsource + v \cdot \nabla \jsource = \xi \jsource$, which is checked by differentiating the action on $\jsource$. At each time step we have the numerical approximations $\jsource^n$, as well as their local time interpolations $\tilde{\jsource}^n$ over $[t_n, t_{n + 1}]$. Recall from assumption \ref{item:error_assumption_etd} that the log-Jacobian integrator approximates the quantity $\tilde{\jmath}^n$, which satisfies the PDE $\partial_t \tilde{\jmath}^n = \mathcal{L} \tilde{\jmath}^n + \mathcal{N}(\tilde{\jmath}^n, \tilde{\jsource}^n, t)$.

We consider the quantities $\tilde{\rho}^n := \tilde{\jsource}^n e^{\tilde{\jmath}^n}$ in $[t_n, t_{n + 1}]$. We note that in general these quantities jump at nodes, $\tilde{\rho}^n (t_{n + 1}) \neq \tilde{\rho}^{n + 1} (t_{n + 1})$. We define their errors relative to $\rho$,
\begin{equation*}
    E^n := \sup_{s \in [t_n, t_{n + 1}]} \lVert \tilde{\rho}^n (s) - \rho(s) \rVert_{C^0}.
\end{equation*}
We have
\begin{equation}
    \label{eq:total_error_is_j_plus_rhotilde}
    \lVert \rho_0 \cdot \numericalGroupElement^n - \rho(t_n) \rVert_{C^0} \lesssim \lVert \numericalJacobian^n - \tilde{\jmath}^n (t_n) \rVert_{C^0} + \lVert \tilde{\rho}^n (t_n) - \rho(t_n) \rVert_{C^0}.
\end{equation}
The first term has size $\delta$ by assumption \ref{item:error_assumption_etd} on the log-Jacobian integrator, and the second one is bounded by $E^n$. Thus, it remains to show $E^n \lesssim \delta$. One checks that in time interval $[t_n, t_{n + 1}]$, $\tilde{\rho}^n$ satisfies
\begin{equation}
    \label{eq:rho_tilde_pde_appendix}
    \partial_t \tilde{\rho}^n + \tilde{u} \cdot \nabla \tilde{\rho}^n = \mathcal{L} \tilde{\rho}^n + G(\tilde{\rho}^n - c, t) + e^{\tilde{\jmath}^n} \left( \partial_t \tilde{\jsource}^n + \tilde{v} \cdot \nabla \tilde{\jsource}^n - \tilde{\xi} \tilde{\jsource}^n \right)
\end{equation}
where we have defined $\tilde{u} := A(\tilde{\rho}^n - c, t)$ and $(\tilde{v}, \tilde{\xi}) := \liftingOperator(\tilde{\rho}^n, t)$. Thus, $\tilde{\rho}^n$ solves the original PDE \eqref{eq:pde_with_highest_order_linear_term}, up to an error in the action of $\jsource$. We proceed by first bounding the error in the action of $\jsource$, then bounding the error in $\tilde{\rho}^n$.

The following lemma bounds the pointwise error in $\jsource$ by writing a recursion in the global map and mass errors. Write these respectively by
\begin{equation*}
    E_\map^n = \numericalMap^n - \map(t_n) \quad\quad\text{and}\quad\quad E_\mass^n = |\numericalMass^n - \mass(t_n)|.
\end{equation*}
For the map $\numericalMap$ we estimate the error on the grid in the norm $\lVert \cdot \rVert_{\ell^\infty(\Gamma)} := \max_\Gamma |\cdot|$, as opposed to over the whole space. This is because in our regime $\Delta t \lesssim \Delta x$, the departure of a semi-Lagrangian method remains close to the grid points. Thus, estimating only the grid error allows us to exploit properties such as the superconvergence of uniform periodic spline derivatives at the nodes \cite{behforooz_overconvergence_1988}. Conversely, iterating the map error over the entire space requires estimating errors over entire cells, which is wasteful and incurs the final error factor of $1/\Delta t$ typically found in semi-Lagrangian methods. The map error over the whole space may then be incorporated as a final interpolation error $\mathcal{O}(\Delta x^6)$ at the end without entering the iteration.

\begin{lemma}
    \label{lemma:jsource_lemma}
    The error in $\tilde{\jsource}^n$ satisfies
    \begin{equation*}
        \sup_{s \in [t_n, t_{n + 1}]} \lVert \tilde{\jsource}^n (s) - \jsource(s) \rVert_{C^0} \lesssim E^n + \delta.
    \end{equation*}
\end{lemma}

\begin{proof}
    We first note that since $\rho^n = T_{N_\epsilon} [\rho_0 \cdot \numericalGroupElement^n]$ we have
    \begin{equation}
        \label{eq:rho_Cp_error}
        \begin{aligned}
            \lVert \rho^n - \rho(t_n) \rVert_{C^p} &\leq \lVert T_{N_\epsilon} [\rho_0 \cdot \numericalGroupElement^n - \rho(t_n)] \rVert_{C^p} + \lVert T_{N_\epsilon} [\rho(t_n)] - \rho(t_n) \rVert_{C^p} \\
            &\lesssim \lVert \rho_0 \cdot \numericalGroupElement^n - \rho(t_n) \rVert_{C^0} + \mathcal{O}(\epsilon) \leq E^n + \delta
        \end{aligned}
    \end{equation}
    using \eqref{eq:total_error_is_j_plus_rhotilde}. Here we have also used $\lVert T_{N_\epsilon} \rVert_{C^0 \rightarrow C^p} = \mathcal{O}(1)$, i.e. we have exploited bandlimitedness to upgrade regularity. Next we note that the first component of the lifting operator $\liftingOperator_1$ is a locally Lipschitz map $C^p \rightarrow C^0$. This follows from the differential orders of $F$ and $A$, and the fact that $\Delta^{-1} : C^{k + 1} \rightarrow C^k$ is bounded on the mean-zero subspace for all $k$. Hence the error in velocity at $t_n$ can be bounded as
    \begin{equation}
        \label{eq:velocity_error}
        \lVert v^n - v(t_n) \rVert_{C^0} = \lVert \liftingOperator(\rho^n, t_n) - \liftingOperator(\rho(t_n), t_n) \rVert_{C^0} \lesssim \lVert \rho^n - \rho(t_n) \rVert_{C^p} \lesssim E^n + \delta
    \end{equation}
    using \eqref{eq:rho_Cp_error}. For the same reason, the error in the mass velocity satisfies $|\xi^n - \xi(t_n)| \lesssim E^n + \delta$. Interpolating $\xi^i$ in time, we have
    \begin{equation*}
        \sup_{s \in [t_n, t_{n + 1}]} |\tilde{\xi}^n (s) - \xi(s)| \lesssim |\xi^n - \xi(t_n)| + \mathcal{O}(\Delta t^4) \lesssim E^n + \delta.
    \end{equation*}
    Then Grönwall's lemma applied to the time integration of the mass gives
    \begin{equation}
        \label{eq:mass_error_iteration}
        E_\mass^{n + 1} \leq E_\mass^n + \mathcal{O} \bigg( \Delta t \sup_{s \in [t_n, t_{n + 1}]} |\tilde{\xi}^n (s) - \xi(s)| + \Delta t^5 \bigg) \leq E_\mass^{n + 1} + \mathcal{O}(\Delta t E^n + \Delta t \, \delta) 
    \end{equation}
    which iterates to $E_\mass^n \lesssim E^n + \delta$.

    For the map, one checks that the error $E_\map^n = \numericalMap^n - \map(t_n)$ satisfies the recursion
    \begin{equation}
        \label{eq:map_error_recursion}
        E_\map^{n + 1} = \slPropagator^n (E_\map^n) + \tau^n
    \end{equation}
    where the local truncation error
    \begin{equation*}
        \tau^n := \intMap^n (\map(t_n), \tilde{v}^n) - \map(t_{n + 1}) = \interp[\map(t_n)] \circ \numericalMap_{[t_{n + 1}, t_n]} - \map(t_{n + 1}).
    \end{equation*}
    Iterating \eqref{eq:map_error_recursion} gives
    \begin{equation*}
        \lVert E_\map^n \rVert_{\ell^\infty(\Gamma)} \leq \sum_{i = 0}^{n - 1} \lVert \slPropagator^{n - 1} \circ \cdots \circ \slPropagator^{i + 1} (\tau^i) \rVert_{\ell^\infty (\Gamma)} \lesssim \sum_{i = 0}^{n - 1} \lVert \tau^i \rVert_{\ell^\infty (\Gamma)}
    \end{equation*}
    by the stability assumption \ref{item:error_assumption_sl}, and so it suffices to measure the local truncation errors. We split $\tau^n$ as
    \begin{equation*}
        \lVert \tau^n \rVert_{\ell^\infty (\Gamma)} \leq \lVert \underbrace{\interp[\map(t_n)] \circ \numericalMap_{[t_{n + 1}, t_n]} - \map(t_n) \circ \numericalMap_{[t_{n + 1}, t_n]}}_{\mathrm{(I)}} \rVert_{\ell^\infty (\Gamma)} + \lVert \underbrace{\map(t_n) \circ \numericalMap_{[t_{n + 1}, t_n]} - \map(t_{n + 1})}_{\mathrm{(II)}} \rVert_{\ell^\infty (\Gamma)}.
    \end{equation*}

    \paragraph{Bounding (I)} Write $\numericalMap_{[t_{n + 1}, t_n]} = \operatorname{id} + \, \eta$ as a deformation of the identity. In particular, we have the size of the displacement $\lVert \eta \rVert_{C^0} = \mathcal{O}(\Delta t)$. Next, write the interpolation error $r = \interp[\map_{[t_n, 0]}] - \map_{[t_n, 0]}$. Now fix a grid point $x \in \Gamma$. Taylor expanding $r$ about $x$ gives
    \begin{equation}
        \label{eq:map_error_decomposition_taylor_term}
        r(\numericalMap_{[t_{n + 1}, t_n]}(x)) = r(x + \eta(x)) = r(x) + \nabla r(x) \cdot \eta(x) + \frac{1}{2} \eta(x)^T \nabla^2 r(\zeta) \eta(x)
    \end{equation}
    for some $\zeta$. The first term $r(x)$ in \eqref{eq:map_error_decomposition_taylor_term} is zero since $x$ is a grid point. For the second term in \eqref{eq:map_error_decomposition_taylor_term}, we can exploit the superconvergence of the derivative at the nodes for uniform spline interpolation \cite{behforooz_overconvergence_1988}: since $x$ is a grid point, we have $|r(x)| = \mathcal{O}(\Delta x^6)$, and we do not lose one order. Thus the second term of \eqref{eq:map_error_decomposition_taylor_term} is bounded by $\mathcal{O}(\Delta t \Delta x^6)$.
    The third term in \eqref{eq:map_error_decomposition_taylor_term} has size
    \begin{equation*}
        \lVert \nabla^2 r \rVert_{C^0} |\eta(x)|^2 = \mathcal{O}(\Delta t^2 \Delta x^4).
    \end{equation*}
    Here we have lost two orders in $\Delta x$ due to the second derivative. Altogether we have
    \begin{equation*}
        \lVert \mathrm{(I)} \rVert_{\ell^\infty (\Gamma)} = \mathcal{O}(\Delta t \Delta x^6 + \Delta t^2 \Delta x^4) = \mathcal{O}(\Delta t \, \delta).
    \end{equation*}

    \paragraph{Bounding (II)} Interpolating the map velocity in space and time gives
    \begin{equation*}
        \sup_{t \in [t_n, t_{n + 1}]} \lVert \tilde{v}(t) - v(t) \rVert_{C^0} \lesssim \lVert \interp[v^n] - v(t_n) \rVert_{C^0} + \mathcal{O}(\Delta t^4) \lesssim \lVert v^n - v(t_n) \rVert_{C^0} + \mathcal{O}(\Delta t^4 + \Delta x^6) \lesssim E^n + \delta
    \end{equation*}
    using \eqref{eq:velocity_error}, and where we have also used the stability of uniform spline interpolation \cite{schurer_note_1968}. Using a Lipschitz constant for $\map_{[t_n, 0]}$ and then applying Grönwall's lemma, we obtain
    \begin{equation*}
        \lVert \mathrm{(II)} \rVert_{C^0} \lesssim \lVert \numericalMap_{[t_{n + 1}, t_n]} - \map_{[t_{n + 1}, t_n]} \rVert_{C^0} = \mathcal{O} \bigg( \Delta t \sup_{t \in [t_n, t_{n + 1}]} \lVert \tilde{v}(t) - v(t) \rVert_{C^0} + \Delta t^5 \bigg) = \mathcal{O}(\Delta t E^n + \Delta t \, \delta)
    \end{equation*}
    where the fourth-order error in the integration is guaranteed since $\tilde{v}^n \in C^4$.
    
    Altogether, we have $\lVert E_\map^n \rVert_{\ell^\infty(\Gamma)} \lesssim E^n + \delta$. Replacing the error on the grid by the error over the entire space by interpolation at the final time, and interpolating in time, gives
    \begin{equation*}
        \sup_{s \in [t_n, t_{n + 1}]} \lVert \tilde{\jsource}^n (s) - \jsource(s) \rVert_{C^0} \lesssim \lVert E_\map^n \rVert_{\ell^\infty(\Gamma)} + E_\mass^n + \mathcal{O}(\Delta t^4 + \Delta x^6) \lesssim E^n + \delta.
    \end{equation*}
\end{proof}

The following lemma bounds the error in $\partial_t \jsource$, which has the same error rate as the pointwise values when measured in the weaker spatial norm $C^{-1}$.

\begin{lemma}
    \label{lemma:jsource_ddt_lemma}
    The error in the time derivative of $\tilde{\jsource}^n$ satisfies
    \begin{equation*}
        \sup_{s \in [t_n, t_{n + 1}]} \lVert \partial_t \tilde{\jsource}^n (s) - \partial_t \jsource(s) \rVert_{C^{-1}} \lesssim E^n + \delta.
    \end{equation*}
\end{lemma}

\begin{proof}
    The proof to establish the local orders is the same as the previous lemma; thus we only illustrate the arguments pertaining to the time derivative and the weaker norm. We need to bound the error in the time derivatives of the map and mass. For the mass, we note that its propagator is linear, and thus advances the time derivative $\partial_t \mass$ at the same local error rate,
    \begin{equation*}
        |\partial_t \mass(t_1) - \partial_t \intMass(\mass(t_0))| = |\partial_t \mass(t_1) - \intMass(\partial_t \mass(t_0))| = \mathcal{O}(\Delta t^5),
    \end{equation*}
    giving the same local truncation error of $\mathcal{O}(\Delta t E^n + \Delta t \, \delta)$ as estimated in \eqref{eq:mass_error_iteration}.

    Now for the map, recall from Lemma \ref{lemma:jsource_lemma} the error $E_\map^n = \numericalMap^n - \map(t_n)$ and its recursion \eqref{eq:map_error_recursion}
    \begin{equation*}
        E_\map^{n + 1} = \slPropagator^n (E_\map^n) + \tau^n
    \end{equation*}
    with the local truncation error $\tau^n$, which satisfies $\lVert \tau^n \rVert_{C^{-1}} \lesssim \lVert \tau^n \rVert_{\ell^\infty(\Gamma)} = \mathcal{O}(\Delta t E^n + \Delta t \, \delta)$. We then have
    \begin{equation*}
        E_\map^{n + 1} - E_\map^n = (\slPropagator^n - I) E_\map^n + \tau^n.
    \end{equation*}
    By assumption \ref{item:error_assumption_sl}, this gives
    \begin{equation*}
        \lVert E_\map^{n + 1} - E_\map^n \rVert_{C^{-1}} \lesssim \Delta t \lVert E_\map^n \rVert_{C^0} + \lVert \tau^n \rVert_{C^{-1}} = \mathcal{O}(\Delta t E^n + \Delta t \, \delta).
    \end{equation*}
    Then, if $\tilde{E}_\map^n$ denotes the interpolation in time of $E_\map^n, \ldots, E_\map^{n - 3}$, the derivative $\partial_t \tilde{E}_\map^n$ depends only on the differences $E_\map^{i + 1} - E_\map^i$ and adds a factor $1/\Delta t$,
    \begin{equation*}
        \lVert \partial_t \tilde{E}_\map^n \rVert_{C^{-1}} \lesssim \frac{1}{\Delta t} \lVert E_\map^{n + 1} - E_\map^n \rVert_{C^{-1}} \lesssim E^n + \delta.
    \end{equation*}
    The same argument applies when working with $q^n$ directly. In any case, we see that we can avoid the usual $1/\Delta t$ loss in the derivative of a time interpolation, at the cost of measuring in a weaker spatial norm.
\end{proof}

Finally, the following lemma bounds the error in $\tilde{\rho}^n$. It is essentially a Duhamel (with jumps) and singular Grönwall argument, and it completes the proof of Theorem \ref{thm:error_estimate}.

\begin{lemma}
    \label{lemma:error_resursion}
    The errors $E^n$ satisfy $E^n \lesssim \delta$.
\end{lemma}

\begin{proof}
    Recall equation \eqref{eq:rho_tilde_pde_appendix}, which says that in $[t_n, t_{n + 1}]$, $\tilde{\rho}^n$ satisfies
    \begin{equation*}
        \partial_t \tilde{\rho}^n = \mathcal{L} \tilde{\rho}^n + H(\tilde{\rho}^n, t) + e^{\tilde{\jmath}^n} h^n,
    \end{equation*}
    where we have defined $H(\rho, t) = G(\rho - c, t) - A(\rho - c, t) \cdot \nabla \rho$ and
    \begin{equation*}
        h^n = \partial_t \tilde{\jsource}^n + \liftingOperator_1 (\tilde{\rho}^n) \cdot \nabla \tilde{\jsource}^n - \liftingOperator_2 (\tilde{\rho}^n) \tilde{\jsource}^n.
    \end{equation*}
    Then on the same interval, we have that $w^n = \tilde{\rho}^n - \rho$ satisfies
    \begin{equation*}
        \partial_t w^n = \mathcal{L} w^n + H(\tilde{\rho}^n, t) - H(\rho, t) + e^{\tilde{\jmath}^n} h^n.
    \end{equation*}
    Applying Duhamel's formula to $w^i$ on each interval $[t_i, t_{i + 1}]$ and combining then gives the Duhamel with jumps formula,
    \begin{equation}
        \label{eq:duhamel_w}
        \begin{aligned}
            w^n (t_{n + 1}) &= \sum_{0 \leq i \leq n} \int_{t_i}^{t_{i + 1}} e^{(t_{n + 1} - s) \mathcal{L}} \big( \underbrace{H(\tilde{\rho}^i (s), s) - H(\rho(s), s)}_{\mathrm{(I)}} + \underbrace{e^{\tilde{\jmath}^i (s)} h^i (s)}_{\mathrm{(II)}} \big) \, ds \\
            &\quad\quad + \sum_{1 \leq i \leq n} e^{(t_{n + 1} - t_i) \mathcal{L}} \big[ \underbrace{w^i(t_i) - w^{i - 1}(t_i)}_{\mathrm{(III)}} \big]
        \end{aligned}
    \end{equation}
    Crucially, by integrating over the entire interval $[t_0, t_{n + 1}]$, we obtain the full smoothing $e^{(t_{n + 1} - s) \mathcal{L}}$. If we were to bound the errors over each subinterval separately, we would obtain a Grönwall inequality which is unbounded as $\Delta t \rightarrow 0$.

    We now bound \eqref{eq:duhamel_w} in $C^0$.

    \paragraph{Bounding $\mathrm{(I)}$}
    We have that $H$ is a locally Lipschitz map $C^{p - 1} \rightarrow C^0$. Conversely, the operator $e^{s \mathcal{L}}$ provides smoothing at the cost of a time factor, $\lVert e^{s \mathcal{L}} \rVert_{C^{p - 1} \rightarrow C^0} \lesssim s^{-(p - 1)/p}$. Combining these gives
    \begin{equation*}
        \lVert e^{(t_{n + 1} - s) \mathcal{L}} (H(\tilde{\rho}^i (s), s) - H(\rho(s), s)) \rVert_{C^0} \lesssim (t_{n + 1} - s)^{-(p - 1)/p} \lVert w^i(s) \rVert_{C^0} \leq (t_{n + 1} - s)^{-(p - 1)/p} E^i.
    \end{equation*}

    \paragraph{Bounding $\mathrm{(II)}$}
    Adding $\partial_t \jsource + v \cdot \nabla \jsource - \xi \jsource = 0$ to the definition of $h^i$ gives
    \begin{equation*}
        h^i = \partial_t \tilde{\jsource}^i - \partial_t \jsource + \liftingOperator_1 (\tilde{\rho}^i) \cdot \nabla \tilde{\jsource}^i - \liftingOperator_1 (\rho) \cdot \nabla \jsource - (\liftingOperator_2 (\tilde{\rho}^i) \tilde{\jsource}^i - \liftingOperator_2 (\rho) \jsource).
    \end{equation*}
    This is order $p - 1$ in $\tilde{\rho}^i - \rho$ and order 1 in $\tilde{\jsource}^i - \jsource$. Thus the smoothing $\lVert e^{s \mathcal{L}} \rVert_{C^{p - 1} \rightarrow C^0} \lesssim s^{-(p - 1)/p}$ gives
    \begin{equation*}
        \begin{aligned}
            \lVert e^{(t_{n + 1} - s) \mathcal{L}} e^{\tilde{\jmath}^i} h^i \rVert_{C^0} &\lesssim (t_{n + 1} - s)^{-(p - 1)/p} \left( \lVert \partial_t \tilde{\jsource}^i - \partial_t \jsource \rVert_{C^{-1}} + \lVert \tilde{\jsource}^i - \jsource \rVert_{C^0} + \lVert w^i \rVert_{C^0} \right) \\
            &\lesssim (t_{n + 1} - s)^{-(p - 1)/p} \left( E^i + \delta \right)
        \end{aligned}
    \end{equation*}
    using Lemmas \ref{lemma:jsource_lemma} and \ref{lemma:jsource_ddt_lemma}.

    \paragraph{Bounding $\mathrm{(III)}$}
    The jumps at $t_i$ are $w^i (t_i) - w^{i - 1} (t_i) = e^{\tilde{\jmath}^i (t_i)} \left( \tilde{\jsource}^i (t_i) - \tilde{\jsource}^{i - 1} (t_i) \right)$, so the smoothing $\lVert e^{s \mathcal{L}} \rVert_{C^1 \rightarrow C^0} \lesssim s^{-1/p}$ gives
    \begin{equation*}
        \lVert e^{(t_{n + 1} - s) \mathcal{L}} \left[ w^i (t_i) - w^{i - 1} (t_i) \right] \rVert_{C^0} \lesssim (t_{n + 1} - t_i)^{-1/p} \lVert \tilde{\jsource}^i (t_i) - \tilde{\jsource}^{i - 1} (t_i) \rVert_{C^{-1}}.
    \end{equation*}
    We Taylor expand $\tilde{\jsource}^i (s) - \tilde{\jsource}^{i - 1} (s)$ about $s = t_{i - 1}$, where $\tilde{\jsource}^i$ and $\tilde{\jsource}^{i - 1}$ agree, to get
    \begin{equation*}
        \tilde{\jsource}^i (t_i) - \tilde{\jsource}^{i - 1} (t_i) = (t_i - t_{i - 1}) \left( \partial_t \tilde{\jsource}^i (\zeta) - \partial_t \tilde{\jsource}^{i - 1} (\zeta) \right),
    \end{equation*}
    for some $\zeta \in [t_{i - 1}, t_i]$. Thus
    \begin{equation*}
        \lVert \tilde{\jsource}^i (t_i) - \tilde{\jsource}^{i - 1} (t_i) \rVert_{C^{-1}} \lesssim \Delta t \left( \lVert \partial_t \tilde{\jsource}^i (\zeta) - \partial_t \jsource(\zeta) \rVert_{C^{-1}} + \lVert \partial_t \jsource(\zeta) - \partial_t \tilde{\jsource}^{i - 1} (\zeta) \rVert_{C^{-1}} \right) = \mathcal{O}(\Delta t E^i + \Delta t \, \delta)
    \end{equation*}
    by Lemma \ref{lemma:jsource_ddt_lemma}, giving $\lVert e^{(t_{n + 1} - s) \mathcal{L}} \left[ w^i (t_i) - w^{i - 1} (t_i) \right] \rVert_{C^0} \lesssim (t_{n + 1} - t_i)^{-1/p} \cdot \Delta t (E^i + \delta)$.

    Finally, assuming that $w^n (s)$ is increasing for simplicity, with these bounds \eqref{eq:duhamel_w} becomes
    \begin{equation*}
        E^n \lesssim \sum_{0 \leq i \leq n} \int_{t_i}^{t_{i + 1}} (t_{n + 1} - s)^{-(p - 1)/p} \left( E^i + \delta \right) ds + \Delta t \sum_{1 \leq i \leq n} (t_{n + 1} - t_i)^{-1/p} \left( E^i + \delta \right).
    \end{equation*}
    Routine integral computations transform this recursion into the form
    \begin{equation}
        \label{eq:weakly_singular_discrete_gronwall}
        E^n \leq C \Delta t^{1/p} E^n + C \Delta t \sum_{0 \leq i \leq n - 1} (t_n - t_i)^{-(p - 1)/p} E^i + C \delta
    \end{equation}
    for some constant $C$. Taking $\Delta t$ small enough that $C \Delta t^{1/p} \leq \frac{1}{2}$, we can move $C \Delta t^{1/p} E^n$ to the left. Then \eqref{eq:weakly_singular_discrete_gronwall} becomes a weakly singular discrete Grönwall inequality, and gives $E^n \lesssim \delta$ uniformly in $\Delta t$, see \cite[Thm. 6.1]{dixon_weakly_1986}.
\end{proof}

\small
\bibliographystyle{unsrt}  
\bibliography{Bibliography}

\end{document}